\documentclass{amsart}

\usepackage{hyperref}
\usepackage[utf8]{inputenc}
\usepackage{color}
\usepackage{amsmath}
\usepackage{amsfonts}
\usepackage{amsthm}
\usepackage{amssymb}
\usepackage{graphicx}
\usepackage{enumerate}
\usepackage{tikz-cd}
\usepackage{tikz}
\usepackage{circuitikz}
\usepackage{caption}
\usepackage{makecell}

\newtheorem{theorem}{Theorem}[section]
\newtheorem{lemma}[theorem]{Lemma}
\newtheorem{proposition}[theorem]{Proposition}
\newtheorem{corollary}[theorem]{Corollary}
\newtheorem{question}[theorem]{Question}

\theoremstyle{definition}
\newtheorem{definition}[theorem]{Definition}
\newtheorem{example}[theorem]{Example}

\theoremstyle{remark}
\newtheorem{remark}[theorem]{Remark}

\numberwithin{equation}{section}

\begin{document}

\title[The deformation spaces of orderable Coxeter 3-polytopes]{The smoothness of the real projective deformation spaces of orderable Coxeter 3-polytopes}

\author{Suhyoung Choi}
\address{Department of Mathematical Sciences, KAIST, 291 Daehak-ro Yuseong-gu Daejeon, South Korea}
\email{schoi@math.kaist.ac.kr}
\thanks{The first and second authors were supported by NRF Grant 2022R1A2C300316213.}

\author{Seungyeol Park}
\address{Department of Mathematical Sciences, KAIST, 291 Daehak-ro Yuseong-gu Daejeon, South Korea}
\email{sypark14@kaist.ac.kr}


\date{}

\dedicatory{}


\begin{abstract}
A Coxeter polytope is a convex polytope in a real projective space equipped with linear reflections in its facets, such that the orbits of the polytope under the action of the group generated by the linear reflections tessellate a convex subset in the real projective space.
Vinberg proved that the group generated by these reflections acts properly discontinuously on the interior of this convex subset, thus inducing a natural orbifold structure on the polytope.

In this paper, we consider labeled combinatorial polytopes $\mathcal{G}$ associated to such orbifolds, and study the deformation space $\mathcal{C} (\mathcal{G})$ of Coxeter polytopes realizing $\mathcal{G}$.
We prove that if $\mathcal G$ is orderable and of normal type, and its underlying combinatorial polytope is not a cone over a polygon, then the deformation space $\mathcal C(\mathcal G)$ is a smooth manifold.
This result is obtained by analyzing a natural map from $\mathcal C(\mathcal G)$ to a smooth manifold called the realization space.
\end{abstract}

\maketitle

\begingroup 
\hypersetup{hidelinks}
\tableofcontents
\endgroup


\section{Introduction}\label{introduction}

Let $\mathbb{S}^d = (\mathbb{R}^{d+1} \setminus \{0\}) / \sim$ be the \textit{real projective} $d$-sphere consisting of the open rays of $\mathbb{R}^{d+1}$ from the origin. 
The group $\text{SL}_\pm (d+1, \mathbb{R})$ can then be viewed as the group of real projective automorphisms of $\mathbb{S}^d$. 
A \textit{real projective structure} on a smooth $d$-orbifold $\mathcal{O}$ is an atlas of coordinate charts on $\mathcal{O}$ valued in $\mathbb{S}^d$ such that the changes of coordinates locally lie in $\text{SL}_\pm (d+1, \mathbb{R})$.
The \textit{deformation space} of real projective structures on $\mathcal{O}$ is the space of real projective structures modulo isotopy.

Many real projective $d$-orbifolds arise as quotients $\Omega / \Gamma$, where $\Omega \subset \mathbb{S}^d$ is a convex open domain and $\Gamma \subset \text{SL}_\pm (d + 1, \mathbb{R})$ is a discrete subgroup acting properly on $\Omega$. 
Such quotients are called \textit{convex real projective orbifolds}.

A large class of discrete groups arising in this setting is provided by \textit{linear reflection groups}, which are generated by reflections in the facets of convex polytopes.
Throughout this paper, a convex polytope in $\mathbb{S}^d$ means a properly convex polytope, namely one contained in an affine chart whose closure is convex and bounded in that chart. 
Vinberg characterized linear reflection groups in terms of the linear data defining their reflections.

Linear reflection groups naturally give rise to Coxeter orbifolds. 
Let $P \subset \mathbb{S}^d$ be a convex polytope and let $\Gamma \subset \text{SL}_\pm (d + 1, \mathbb{R})$ be generated by reflections in its facets. 
Vinberg's theorem \cite{MR0302779} implies that the union of the images of $P$ under the elements of $\Gamma$ is convex, that its interior $\Omega$ is convex, and that $\Omega/\Gamma$ carries a natural orbifold structure.
The stabilizer of an interior point of a codimension-one face is isomorphic to $\mathbb Z/2\mathbb Z$, whereas the stabilizer of an interior point of a codimension-two face $e$ is the dihedral group of order $2m(e)$, where $m(e)\geq2$ is the label of $e$.
It follows that only faces of codimension at least $3$ may be removed.
In the usual terminology, every Coxeter polytope considered in this paper is $2$-perfect, since every face of codimension at most $2$ has a finite stabilizer.
We call the resulting orbifolds \textit{reflection orbifolds} or \textit{Coxeter orbifolds}.

The aim of this paper is to study deformation spaces of real projective structures on $3$-dimensional Coxeter orbifolds, where the removed faces are some vertices of the underlying polytope.
To analyze the topology of these spaces without working directly with orbifold terminology, we encode each Coxeter orbifold by a \textit{labeled combinatorial polytope} $\mathcal G$ and study Coxeter polytopes realizing $\mathcal G$.

We first recall the underlying combinatorial object.
A \textit{combinatorial $d$-polytope} is an isomorphism class of finite regular CW-complexes arising from convex $d$-polytopes, where the open cells are the relative interiors of the faces.
Let $\mathcal G$ be a combinatorial $d$-polytope, and denote its facets by $\mathcal G_1,\ldots,\mathcal G_f$.
Two facets are called \textit{adjacent} when their intersection is a ridge.
A \textit{labeled combinatorial polytope} is obtained from $\mathcal G$ by assigning an integer $m_{i,j}\geq 2$ to each ridge $\mathcal G_i\cap\mathcal G_j$ determined by a pair of adjacent facets.

A \textit{Coxeter $d$-polytope realizing $\mathcal G$} consists of a convex $d$-polytope $Q\subset\mathbb S^d$, whose facets $Q_1,\ldots,Q_f$ correspond to $\mathcal G_1,\ldots,\mathcal G_f$, together with \textit{linear reflections} $r_1,\ldots,r_f$ in the projective hyperplanes supporting these facets.
We further require that, for each pair of adjacent facets, the product $r_ir_j$ has the prescribed order $m_{i,j}$, and that the interiors of the images of $Q$ under distinct elements of the group generated by $r_1,\ldots,r_f$ are disjoint.
Vinberg's results \cite{MR0302779} imply that these images tessellate a convex subset of $\mathbb S^d$, and hence each Coxeter polytope realizing $\mathcal G$ determines a convex real projective structure on the Coxeter orbifold associated with $\mathcal G$.
	
The \textit{deformation space} $\mathcal C(\mathcal G)$ is the space of Coxeter polytopes realizing $\mathcal G$ modulo projective equivalence, where a projective transformation acts simultaneously on the polytope and, by conjugation, on its linear reflections.
Using Vinberg's results, we will identify $\mathcal C(\mathcal G)$ with the solution space of certain polynomial equalities and inequalities \cite{MR0302779}.
Precise definitions of combinatorial and labeled combinatorial polytopes, Coxeter polytopes realizing $\mathcal G$, and the deformation space $\mathcal C(\mathcal G)$ are given in Sections~\ref{Section 2.1} and~\ref{Section 2.2}.

\begin{figure}[!ht]
	\centering
	\resizebox{0.6\textwidth}{!}{%
		\begin{circuitikz}
			\tikzstyle{every node}=[font=\normalsize]
			\draw [short] (11.25,19.75) -- (9,16);
			\draw [short] (9,16) -- (12,14.5);
			\draw [short] (12,14.5) -- (15.5,16);
			\draw [short] (11.25,19.75) -- (15.5,16);
			\draw [short] (11.25,19.75) -- (12,14.5);
			\draw [dashed] (9,16) -- (15.5,16);
			
			\node [font=\LARGE] at (18,17) {$\to$};
			
			\draw [short] (13.5,17.75) -- (13.25,15);
			\draw [dashed] (13.5,17.75) -- (14.25,16);
			\draw [dashed] (13.25,15) -- (14.25,16);
			\draw [short] (22.25,19.5) -- (20,15.75);
			\draw [short] (20,15.75) -- (23,14.25);
			\draw [short] (23,14.25) -- (24.25,14.75);
			
			\draw [short] (22.25,19.5) -- (24.5,17.5);
			
			\draw [short] (22.25,19.5) -- (23,14.25);
			
			\draw [dashed] (20,15.75) -- (25.25,15.75);
			
			\draw [short] (24.5,17.5) -- (24.25,14.75);
			
			\draw [short] (24.5,17.5) -- (25.25,15.75);
			
			\draw [short] (24.25,14.75) -- (25.25,15.75);
			\draw [short] (12.25,20.25) -- (12.5,13);
			\draw [short] (12.25,20.25) -- (14.5,21.75);
			\draw [short] (14.5,21.75) -- (14.75,14.5);
			\draw [short] (12.5,13) -- (14.75,14.5);
			\node [font=\LARGE] at (13.5,19.5) {$\Pi$};
			\node [font=\LARGE] at (16,16) {$v$};
			\node [font=\LARGE] at (11,16.5) {$P$};
			\node [font=\LARGE] at (23.5,13.5) {$P^{\dagger v}$};
		\end{circuitikz}
	}%
	\caption{A truncation of a convex polytope $P$ by a hyperplane $\Pi \subset \mathbb{S}^3$.}
	\label{test}
\end{figure}

In related work, Marquis \cite{MR2660566} studied labeled combinatorial polytopes whose underlying polytopes are \textit{truncation polytopes}.
We first recall the class of truncation polytopes.
Let $\mathcal G$ be a combinatorial $d$-polytope and let $v$ be one of its vertices.
The truncation of $\mathcal G$ at $v$ is obtained by replacing $v$ with a new facet whose face poset is the link of $v$, with all other incidences inherited naturally.
A \textit{truncation $d$-polytope} is obtained from a $d$-simplex by successive vertex truncations.
Geometrically, if a convex $d$-polytope $P$ realizes $\mathcal G$ and $p$ is the vertex of $P$ corresponding to $v$, this operation can be realized by choosing a projective hyperplane $\Pi$ sufficiently close to $p$ so that it separates $p$ from all other vertices of $P$ and intersects precisely the edges incident to $p$, and then intersecting $P$ with the closed half-space bounded by $\Pi$ that contains the remaining vertices.
The new facet cut out by $\Pi$ has face poset equal to the link of $v$, so the resulting polytope realizes the combinatorial truncation of $\mathcal G$ at $v$.
In dimension $3$, a truncation polytope is therefore obtained from a tetrahedron by successively truncating vertices in this geometric sense, as illustrated in Figure \ref{test}.

Marquis proved that, when the underlying polytope is a truncation $3$-polytope, $\mathcal C(\mathcal G)$ is a disjoint union of copies of $\mathbb{R}^k$ for some $k\geq0$ \cite{MR2660566}.
An analogous result holds for higher-dimensional truncation polytopes \cite{MR4524212}.
These results motivate the following question.

\begin{question}
	What can be determined about the local and global structures of the deformation space of Coxeter 3-orbifolds under suitable combinatorial conditions?
\end{question}

Marquis's result gives a particularly simple global description of the deformation space in the truncation-polytope case. 
For more general Coxeter orbifolds, however, even the local topology of the deformation space can be more complicated. 
A natural first step toward Question~1.1 is to determine when a Coxeter deformation space is a smooth manifold.

Previous work gives both positive and negative results in this direction. Choi--Hodgson--Lee and Choi--Lee found classes of Coxeter 3-orbifolds for which the point corresponding to the unique hyperbolic structure has a neighborhood homeomorphic to an open cell \cite{MR2944525,MR3375519}. 
By contrast, Choi--Lee exhibited a Coxeter 4-orbifold whose deformation space is the union of two lines meeting at a single point \cite[Section~5.3]{MR3375519}.

In this paper, we provide a sufficient condition on a labeled combinatorial $3$-polytope $\mathcal G$ for the space $\mathcal C(\mathcal G)$ to be a smooth manifold.
The main condition is \emph{orderability}, whose precise definition is given in Definition~\ref{orderable}.
It provides an indexing $\mathcal G_1,\ldots,\mathcal G_f$ of the facets such that, for each $i$, at most three edges of $\mathcal G_i$ either have order $2$ or are shared with a previously indexed facet $\mathcal G_j$ with $j<i$.
In the construction of local charts in Section~\ref{Section 6.3}, this indexing allows the vectors defining the facet reflections to be determined successively by solving systems of four linear equations.

In addition to orderability, we assume that $\mathcal G$ is of normal type; the precise exceptional cases excluded by this condition are listed in Definition~\ref{normal type}.
We impose one further combinatorial restriction in the main theorem: the underlying combinatorial polytope is not a cone over a polygon.
Under this additional restriction, the realization-space and principal-bundle arguments used in Section~\ref{proof of the main theorem} apply.

Choi proved that, under the orderability and normal-type hypotheses, every nonempty restricted deformation space obtained by fixing the projective equivalence class of the fundamental polytope is a smooth manifold \cite{MR2247648}.
Building on this result, our main theorem establishes the smoothness of the full deformation space under the additional combinatorial hypothesis above.

\begin{theorem}\label{main theorem}
	Let $\mathcal G$ be a labeled combinatorial $3$-polytope.
	Let $f$ and $e_2$ be the numbers of facets and edges of order $2$ of $\mathcal G$, respectively.
	Suppose that $\mathcal G$ is orderable and of normal type.
	Suppose further that the underlying combinatorial polytope of $\mathcal G$ is not a cone over a polygon.
	Then $\mathcal C(\mathcal G)$ is a smooth manifold.
	If it is nonempty, its dimension is $3f-e_2-9$.
\end{theorem}

Although the tetrahedral case is excluded from Theorem~\ref{main theorem}, it can be handled directly under the remaining hypotheses.
Since all realizations of a tetrahedron are projectively equivalent as realizations of $\mathcal G$, the full deformation space, when nonempty, is a restricted deformation space; hence Theorem~\ref{restricted}, together with Proposition~\ref{prop:stabilizer-dimension}, proves its smoothness and gives $\dim\mathcal C(\mathcal G)=3-e_2$ without using the argument of Section~\ref{proof of the main theorem}.

The proof is based on an idea of Hodgson formulated in Corollary~2 of \cite{MR2247648}.
We study the natural forgetful map $\kappa:\mathcal C(\mathcal G)\to\mathcal{RS}(\mathcal G)$ which forgets the reflection data and records only the projective equivalence class of the convex polytope, together with the correspondence between its faces and those of $\mathcal G$.
Here $\mathcal{RS}(\mathcal G)$, defined precisely in Definition~\ref{projective realization}, is the \textit{projective realization space} of $\mathcal G$; roughly speaking, it is the space of convex $3$-polytopes realizing the underlying combinatorial polytope of $\mathcal G$, with the correspondence between their faces and those of $\mathcal G$ retained, modulo projective transformations.
See also Remark~\ref{face correspondence of realizations} for this marking convention.
Filpo Molina proved that the realization space is smooth when the underlying combinatorial polytope is not a cone over a polygon \cite{FilpoMolina}, while Choi's result gives smoothness of every nonempty fiber of $\kappa$ \cite{MR2247648}.
The map $\kappa$ is not a fiber bundle in general; Section~\ref{counterexample} gives an explicit example in which it is neither surjective nor a fiber bundle.
The main step to prove Theorem \ref{main theorem} is to construct compatible local product charts on the lifted deformation space and then pass to the quotient by $\text{SL}_\pm(4,\mathbb{R})\times\mathbb{R}_+^f$.

Theorem~\ref{main theorem} is compatible with Marquis's stronger description in the \textit{simple} (i.e., each vertex of $\mathcal{G}$ is joined to exactly three edges) case.
Every simple orderable combinatorial $3$-polytope is a truncation polytope \cite[Remark~3.20]{MR2660566}.
Let $e$ be the number of edges of $\mathcal G$ and set $e_+:=e-e_2$.
If the underlying polytope of $\mathcal{G}$ is a truncation polytope, Marquis proved that $\mathcal C(\mathcal G)$ is a disjoint union of copies of $\mathbb{R}^{e_+-3}$ \cite[Theorem~3.16]{MR2660566}.
Since a truncation $3$-polytope satisfies $e=3(f-2)$, we have $e_+-3=(e-e_2)-3=3f-e_2-9$, in agreement with Theorem~\ref{main theorem}.
Thus, our theorem extends the smoothness conclusion, though not the global cell-decomposition conclusion, to non-simple orderable polytopes satisfying the remaining hypotheses of Theorem~\ref{main theorem}.

Table~\ref{table} compares, for combinatorial $3$-polytopes with four to seven facets, the total number of order-$2$ patterns, the number of such patterns on truncation polytopes, and the number of orderable patterns.
For a fixed combinatorial polytope, the \textit{order-$2$ pattern} of a labeling is the set of edges carrying the label $2$.
Thus, two labelings of the same underlying combinatorial polytope determine the same order-$2$ pattern if and only if precisely the same edges carry the label $2$.
This is the relevant data for the count since orderability depends only on the order-$2$ pattern and not on the particular values of the labels greater than $2$.
The counts were obtained computationally, and the Mathematica code is available at the cited webpage \cite{mathematica}.

\begin{table}[!ht]
	\centering
	\resizebox{\linewidth}{!}{%
		\begin{tabular}{lllll}
			\hline
			\textbf{\# of facets}
			& \textbf{\makecell{All order-$2$ \\ patterns}}
			& \textbf{\makecell{Patterns on \\ truncation polytopes}}
			& \textbf{\makecell{Orderable \\ patterns}}
			& \textbf{\makecell{Ratio of orderable \\ patterns}} \\ \hline
			4 & 64 & 64 & 64 & 1 \\
			5 & 768 & 512 & 654 & 0.851563 \\
			6 & 14848 & 4096 & 7130 & 0.480199 \\
			7 & 421888 & 98304 & 157334 & 0.372928 \\ \hline
		\end{tabular}%
	}
	\caption{Counts of order-$2$ patterns of combinatorial $3$-polytopes with four to seven facets, together with the corresponding truncation and orderable counts.}
	\label{table}
\end{table}

\subsection{Outline of the Paper}

Section~\ref{Section 2} introduces combinatorial and labeled combinatorial polytopes, Coxeter polytopes, and their deformation spaces, while Section~\ref{section 3} recalls Vinberg's criterion for linear reflection groups.
Section~\ref{Restricted deformation spaces} reviews Choi's theorem on restricted deformation spaces, and Section~\ref{Realization spaces} recalls Filpo Molina's smoothness theorem for realization spaces.
Section~\ref{proof of the main theorem} constructs the smooth structure on $\mathcal C(\mathcal G)$ and proves the main theorem.
Section~\ref{section 7} computes two examples, including one showing that the forgetful map $\kappa$ need not be a fiber bundle.

\subsection*{Acknowledgments}

We are grateful to Gye-Seon Lee and Jean-Paul Filpo Molina for valuable information and discussions, and to Ludovic Marquis for helpful information and lectures. 
We also thank the referee for a careful reading of the manuscript and for numerous detailed and constructive comments, which led us to clarify several arguments and improve the exposition.

\section{Deformation spaces of Coxeter polytopes}\label{Section 2}

In Section~\ref{Section 2.1}, we recall the notions of combinatorial polytopes and labeled combinatorial polytopes.
In Section \ref{Section 2.2}, we recall the notion of the \textit{deformation space of Coxeter 3-polytopes} $\mathcal{C} (\mathcal{G})$ realizing a labeled combinatorial polytope $\mathcal{G}$.

\subsection{Combinatorial and labeled combinatorial polytopes}\label{Section 2.1}

We recall the definition of labeled combinatorial polytopes and related notions such as their associated Coxeter groups and Coxeter graphs.

Let $\mathbb{S}^d$ denote the \textit{projective} $d$-\textit{sphere}, defined as the quotient space $(\mathbb{R}^{d+1} \setminus \{0\}) / \sim$, where two nonzero vectors $v$ and $w$ of $\mathbb{R}^{d+1}$ are equivalent if and only if $v = tw$ for some $t > 0$.
For each subset $A \subset \mathbb{R}^{d+1}$, we let $\mathbb{S} (A)$ denote the image of $A \setminus \{0\}$ under the projection $\mathbb{R}^{d+1} \setminus \{0\} \to \mathbb{S}^d$.
An \textit{affine chart} of $\mathbb{S}^d$ is a subset of the form $\{[v] \in \mathbb{S}^d \ | \ \alpha (v) > 0 \}$ for some nonzero linear functional $\alpha$ on $\mathbb{R}^{d+1}$.
Any affine chart $\Omega$ of $\mathbb{S}^d$ can be identified with $\mathbb{R}^d$ in such a way (though not uniquely) that for each nontrivial linear subspace $W \subset \mathbb{R}^{d+1}$, the intersection $\mathbb{S} (W) \cap \Omega$ is either empty or forms an affine subspace of $\Omega = \mathbb{R}^d$.
For points $p=[v]$ and $q=[w]$ in $\mathbb{S}^d$, the \textit{projective line segment} joining $p$ and $q$ is defined by
\[
[p,q]_{\mathbb{S}}
:= \mathbb{S}\bigl(\{sv+tw \mid s,t \geq 0,\ (s,t) \neq (0,0)\}\bigr).
\]
This subset of $\mathbb{S}^d$ is independent of the choice of the representatives $v$ and $w$ of the positive rays $p$ and $q$.
If $p$ and $q$ are contained in a common affine chart $\Omega$, then $[p,q]_{\mathbb{S}}$ is contained in $\Omega$ and agrees with the usual affine line segment under any identification of $\Omega$ with $\mathbb{R}^d$ as above.
A subset $C \subset \mathbb{S}^d$ is said to be \textit{convex} if $[p,q]_{\mathbb{S}} \subset C$ for every pair of points $p,q \in C$.
In particular, if $C$ is contained in an affine chart, this definition agrees with the usual affine notion of convexity in that chart.
Furthermore, the set $C$ is said to be \textit{properly convex} if it is convex and its closure $\overline{C}$ is a convex and bounded subset of some affine chart of $\mathbb{S}^d$.

A convex $d$-polytope is a properly convex subset of $\mathbb{S}^d$ with nonempty interior, given by
\begin{align}\label{convex polytope}
	P = \{ [v] \in \mathbb{S}^d \mid \alpha_1(v) \geq 0, \dots, \alpha_f(v) \geq 0 \}
\end{align}
for some nonzero linear functionals
$\alpha_1,\ldots,\alpha_f \in (\mathbb R^{d+1})^*$.
For a subset presented as in \eqref{convex polytope}, the proper convexity of $P$ is equivalent to
\[
\text{span}\{\alpha_1,\ldots,\alpha_f\}=(\mathbb R^{d+1})^*.
\]
Throughout this paper, we will assume that none of these defining inequalities are redundant, meaning none of the inequalities $\alpha_j \geq 0$ is implied by the others.
The boundary $\partial P$ of $P$ in $\mathbb{S}^d$ is the union of the \textit{faces} of $P$, each of which is a convex $k$-polytope for some $k < d$ in a $k$-dimensional projective subspace of $\mathbb{S}^d$.
We call the 0-dimensional (resp., 1-dimensional) faces of $P$ the \textit{vertices} (resp., \textit{edges}).
We call the codimension 1 (resp., codimension 2) faces of $P$ the \textit{facets} (resp., \textit{ridges}).
Since the inequalities $\alpha_1 \geq 0, \ldots, \alpha_f \geq 0$ are assumed to be non-redundant, the convex $d$-polytope $P$ has exactly $f$ facets.

Let $P \subset \mathbb{S}^d$ be a convex $d$-polytope.
The \textit{face lattice} $\textup{FL}(P)$ of $P$ is the set of faces of $P$ partially ordered by inclusion.
Two convex $d$-polytopes $P$ and $P'$ are \textit{combinatorially equivalent} if there is a bijection $\phi:\textup{FL}(P)\to\textup{FL}(P')$ preserving the inclusion relations: $F_1\subset F_2$ if and only if $\phi(F_1)\subset\phi(F_2)$.
A \textit{combinatorial $d$-polytope} $\mathcal G$ is a combinatorial equivalence class of convex $d$-polytopes.
By a slight abuse of terminology, we also say that a convex $d$-polytope $P$ is combinatorially equivalent to a combinatorial $d$-polytope $\mathcal G$ if the combinatorial equivalence class of $P$ is $\mathcal G$.

We will also use the following equivalent description of a combinatorial polytope.
Every convex $d$-polytope $P$ has a natural structure of a finite regular CW-complex whose open cells are the relative interiors of the faces of $P$.
The same data are encoded by the isomorphism class of this regular CW-complex, where an isomorphism is a homeomorphism that maps each open cell onto an open cell.
Accordingly, we regard a combinatorial $d$-polytope $\mathcal G$ as an isomorphism class of such regular CW-complexes.
The faces of $\mathcal G$ correspond to the closures of the cells in any representative, and the terms \textit{vertices, edges, facets, ridges}, and so on are understood according to their dimensions.
For example, if $P$ is a representative of $\mathcal G$ with facets $P_1,\ldots,P_f$, then these facets represent the facets $\mathcal G_1,\ldots,\mathcal G_f$ of $\mathcal G$, respectively.

Two facets $\mathcal{G}_i$ and $\mathcal{G}_j$ are said to be \textit{adjacent} if their intersection $\mathcal{G}_i \cap \mathcal{G}_j$ is a ridge of $\mathcal{G}$, i.e., a codimension 2 face of $\mathcal{G}$.
For each pair of adjacent facets, $\mathcal{G}_i$ and $\mathcal{G}_j$, we assign an integer $m_{i,j} \geq 2$ to the ridge $\mathcal{G}_i \cap \mathcal{G}_j$.
We set $m_{i,j}:=\infty$ whenever $\mathcal G_i$ and $\mathcal G_j$ are not adjacent.
The combinatorial polytope $\mathcal{G}$, together with the assigned integers, is called a \textit{labeled combinatorial polytope}.
The combinatorial polytope obtained by forgetting the ridge labels is called the \textit{underlying combinatorial polytope} of the labeled combinatorial polytope.

Let $\mathcal G$ be a labeled combinatorial polytope with facets $\mathcal G_1,\ldots,\mathcal G_f$.
The \textit{Coxeter group associated with} $\mathcal G$ is the group with presentation
\begin{equation}\label{labeled polytope Coxeter group}
	\Gamma_{\mathcal G}
	:=
	\left\langle
	\gamma_1,\ldots,\gamma_f
	\ \middle|\
	\gamma_i^2=1,\ 
	(\gamma_i\gamma_j)^{m_{i,j}}=1
	\right\rangle,
\end{equation}
where the relations $\gamma_i^2=1$ are imposed for $i=1,\ldots,f$, and the relations $(\gamma_i\gamma_j)^{m_{i,j}} = 1$ are imposed for $1\leq i<j\leq f$ with $m_{i,j}<\infty$.

\begin{definition}\label{Coxeter graph}
	Let $\Gamma=\Gamma_{\mathcal G}$ be the Coxeter group with presentation \eqref{labeled polytope Coxeter group}.
	The \textit{Coxeter graph} $G_\Gamma$ for this presentation is the weighted graph defined by the following conditions:
	\begin{enumerate}[(i)]
		\item the vertices (or the nodes) are the generators $\gamma_1, \ldots, \gamma_f$;
		\item two vertices $\gamma_i$, $\gamma_j$ are joined by a single edge if and only if $m_{i, j} \in \{3, 4, \ldots, \infty \}$, and in that case the edge is labeled by $m_{i, j}$.
	\end{enumerate}
	The labels $m_{i, j} = 3$ are conventionally omitted in the representations of Coxeter graphs; however, this convention will not be relevant in this paper.
	We say that the Coxeter group $\Gamma$ is \textit{irreducible} if its associated Coxeter graph $G_\Gamma$ is a connected graph.
\end{definition}

\subsection{The deformation space $\mathcal{C} (\mathcal{G})$ of Coxeter polytopes}\label{Section 2.2}

In this section, we define the notion of Coxeter polytopes and their deformation spaces.  
We will also explain the topology assigned to the deformation spaces.

A linear automorphism $r \in \text{SL}_\pm (d+1,\mathbb{R})$ is called a
\textit{linear reflection} if
\[
r^2=\text{Id}
\quad\text{and}\quad
\dim \ker(r+\text{Id})=1 .
\]
Equivalently, there exist a vector $v\in \mathbb{R}^{d+1}$ and a linear
functional $\alpha\in (\mathbb{R}^{d+1})^*$ such that $\alpha(v)=2$ and
\[
r(x)=x-\alpha(x)v, \qquad x\in \mathbb{R}^{d+1}.
\]
Here $v$ spans the $(-1)$-eigenspace of $r$, while the hyperplane
$\ker \alpha$ is fixed pointwise by $r$. 
We use the same symbol $r$ for the induced projective automorphism $r:\mathbb{S}^d\to \mathbb{S}^d$. 
This automorphism sends the class $[v]\in \mathbb{S}^d$ to its antipodal point $[-v]$ and fixes the projective hyperplane
$\mathbb{S} (\ker\alpha)$ pointwise.

\begin{definition}\label{Coxeter polytope}
	Let $\mathcal{G}$ be a labeled combinatorial polytope with $f$ facets $\mathcal{G}_1, \ldots, \mathcal{G}_f$. 
	A \textit{Coxeter $d$-polytope realizing} $\mathcal{G}$ is a convex $d$-polytope $P \subset \mathbb{S}^d$ together with linear reflections $r_1, \ldots, r_f \in \text{SL}_\pm (d+1, \mathbb{R})$ such that:
	\begin{enumerate}[(i)]
		\item the equivalence class of $P$ is $\mathcal{G}$, with facets $P_1, \ldots, P_f \subset P$ corresponding to the facets $\mathcal{G}_1, \ldots, \mathcal{G}_f$;
		\item each $r_j$ is a linear reflection in the hyperplane of $\mathbb{S}^d$ supporting the facet $P_j$;
		\item if the facets $P_i, P_j$ are adjacent and the integer assigned to the ridge $P_i \cap P_j$ is $m_{i,j}$, then the product $r_i r_j$ has order $m_{i,j}$ in the group $\text{SL}_\pm (d+1, \mathbb{R})$; and
		\item the subgroup $\Gamma \subset \text{SL}_\pm (d+1, \mathbb{R})$ generated by the linear reflections $r_1, \ldots, r_f$ satisfies the condition
		\[
		\text{Int} (P) \cap \gamma \cdot \text{Int} (P) = \varnothing \quad \text{for} \ \gamma \in \Gamma \setminus \{1\}.
		\]
	\end{enumerate}
\end{definition}

Condition \textup{(iv)} can be verified using Vinberg's criterion in Theorem \ref{Vinberg 1}.

Note that the group $\text{SL}_\pm (d+1, \mathbb{R})$ acts (on the left) on the Coxeter $d$-polytopes realizing $\mathcal{G}$ by
\begin{align}\label{action on Coxeter polytopes}
	g \cdot (P, r_1, \ldots, r_f) := (g (P), g r_1 g^{-1}, \ldots, g r_f g^{-1}), \quad g \in \text{SL}_\pm (d+1, \mathbb{R}).
\end{align}

\begin{definition}
	Let $\mathcal{G}$ be a labeled combinatorial polytope with $f$ facets $\mathcal{G}_1, \ldots, \mathcal{G}_f$. 
	The \textit{deformation space of Coxeter $d$-polytopes realizing} $\mathcal{G}$ is the space
	\[
	\mathcal{C} (\mathcal{G}) := \{ \text{Coxeter} \ d\text{-polytopes} \ (P, r_1, \ldots, r_f) \ \text{realizing} \ \mathcal{G} \} / \text{SL}_\pm (d+1, \mathbb{R}),\footnote{All group actions considered in this paper act on spaces from the \textit{left}. 
		However, we will write the group on the right-hand side of the space.
	}
	\]
	where the action of $\text{SL}_\pm (d + 1, \mathbb{R})$ is given by \eqref{action on Coxeter polytopes}.
\end{definition}

Let $\mathbb{R}_+$ be the multiplicative group of positive real numbers.
The product Lie group $\text{SL}_\pm (d+1, \mathbb{R}) \times \mathbb{R}_+^f$ acts on the vector space $\left((\mathbb{R}^{d+1})^*\right)^f \times (\mathbb{R}^{d+1})^f$ by
\begin{align}\label{action}
	&(A, c_1, \ldots, c_f) \cdot (\alpha_1, \ldots, \alpha_f, v_1, \ldots, v_f) \notag \\ 
	&\quad := (c_1^{-1} \alpha_1 \circ A^{-1}, \ldots, c_f^{-1} \alpha_f \circ A^{-1}, c_1 A v_1, \ldots, c_f A v_f). 
\end{align}
To topologize $\mathcal{C} (\mathcal{G})$, we identify the space $\mathcal{C} (\mathcal{G})$ as a subset of the quotient space
\begin{align}\label{quotient}
	\frac{\left((\mathbb R^{d+1})^*\right)^f\times(\mathbb R^{d+1})^f}
	{\text{SL}_{\pm}(d+1,\mathbb R)\times\mathbb R_+^f}
\end{align}
in the following manner. 

Let $(P, r_1, \ldots, r_f)$ be a tuple representing an element of $\mathcal{C} (\mathcal{G})$, where $r_1, \ldots, r_f \in \text{SL}_\pm (d+1, \mathbb{R})$ are linear reflections in the facets $P_1, \ldots, P_f$ of the convex $d$-polytope $P$, generating a linear reflection group. 
For each $j \in \{1, \ldots, f\}$, we have an expression $r_j = \text{Id} - \alpha_j \otimes v_j$ for some $\alpha_j \in (\mathbb{R}^{d+1})^*$ and $v_j \in \mathbb{R}^{d+1}$ with $\alpha_j(v_j) = 2$ and $\alpha_j (x) \geq 0$ for $x \in P$.
Then $P$ can be expressed as \eqref{convex polytope}.

For each element $[(P, r_1, \ldots, r_f)] \in \mathcal{C} (\mathcal{G})$, the tuple $(r_1, \ldots, r_f)$ of linear reflections is uniquely determined by the action of the group $\text{SL}_\pm (d+1, \mathbb{R})$ where the action is given by
\[ 
g \cdot (r_1, \ldots, r_f) := (g r_1 g^{-1}, \ldots, g r_f g^{-1}). 
\]

Moreover, for each reflection $r_j = \text{Id} - \alpha_j \otimes v_j$ with $\alpha_j (v_j) = 2$, the tuple $(\alpha_j, v_j)$ is determined by the action of the multiplicative group $\mathbb{R}_+$ given by
\[ 
c \cdot (\alpha_j, v_j) := (c^{-1} \alpha_j, c v_j). 
\]
Hence, the element $[(P, r_1, \ldots, r_f)] \in \mathcal{C} (\mathcal{G})$ determines a unique element
\[ \Phi ([(P, r_1, \ldots, r_f)]) := [(\alpha_1, \ldots, \alpha_f, v_1, \ldots, v_f)] \]
of the orbit space in \eqref{quotient}.

Thus, we have obtained a map 
\[ \Phi: \mathcal{C} (\mathcal{G}) \to \frac{\left((\mathbb R^{d+1})^*\right)^f\times(\mathbb R^{d+1})^f}
{\text{SL}_{\pm}(d+1,\mathbb R)\times\mathbb R_+^f}. \]
By the following lemma, the map $\Phi$ is injective.
Throughout this paper, we give $\mathcal{C} (\mathcal{G})$ the topology induced by the injective map $\Phi$.
This identification will be useful in the proof of our main theorem, particularly in relating the space $\mathcal{C} (\mathcal{G})$ to the \textit{realization space}, which will be discussed in Section~\ref{Realization spaces}.
\begin{lemma}\label{embedding lemma}
	The map $\Phi$ is injective.
	In particular, there is a unique topology on $\mathcal{C} (\mathcal{G})$ such that the map $\Phi$ is a topological embedding.
\end{lemma}

\begin{proof}
	Let $[(P,r_1,\ldots,r_f)]$ and $[(P',r'_1,\ldots,r'_f)]$ be elements of $\mathcal{C}(\mathcal{G})$ having the same image under $\Phi$.
	For each $i=1,\ldots,f$, choose $\alpha_i,\alpha'_i\in(\mathbb{R}^{d+1})^*$ and $v_i,v'_i\in\mathbb{R}^{d+1}$ such that
	\[
	r_i=\text{Id}-\alpha_i\otimes v_i,
	\qquad
	r'_i=\text{Id}-\alpha'_i\otimes v'_i,
	\qquad
	\alpha_i(v_i)=\alpha'_i(v'_i)=2,
	\]
	where $\alpha_i$ and $\alpha'_i$ are nonnegative on the cones over $P$ and $P'$, respectively.
	
	Since the two elements have the same image under $\Phi$, there exist $A\in\text{SL}_{\pm}(d+1,\mathbb{R})$ and $c_1,\ldots,c_f\in\mathbb{R}_+$ such that $\alpha'_i=c_i^{-1}\alpha_i\circ A^{-1}$, $v'_i=c_iAv_i$ for every $i=1,\ldots,f$.
	It follows that
	\begin{align*}
		r'_i
		&=\text{Id}-\alpha'_i\otimes v'_i \\
		&=\text{Id}-(c_i^{-1}\alpha_i\circ A^{-1})\otimes(c_iAv_i) \\
		&=A(\text{Id}-\alpha_i\otimes v_i)A^{-1} \\
		&=Ar_iA^{-1}.
	\end{align*}
	Moreover, since each $c_i$ is positive,
	\begin{align*}
		P'
		&=\left\{[x]\in\mathbb{S}^d\mid \alpha'_i(x)\geq 0\text{ for every }i\right\} \\
		&=\left\{[x]\in\mathbb{S}^d\mid \alpha_i(A^{-1}x)\geq 0\text{ for every }i\right\} \\
		&=A(P).
	\end{align*}
	Thus, we obtain $(P',r'_1,\ldots,r'_f) = A\cdot(P,r_1,\ldots,r_f)$, so the two tuples represent the same element of $\mathcal{C}(\mathcal{G})$.
	Therefore, $\Phi$ is injective.
\end{proof}

\section{The theory of Vinberg}\label{section 3}

Vinberg \cite{MR0302779} characterized the discrete subgroups $\Gamma$ of $\text{SL}_\pm (d+1, \mathbb{R})$ generated by reflections in the facets of a convex polytope in $\mathbb{S}^d$, in terms of the linear functionals and vectors defining these reflections. It is shown that these linear functionals and vectors must satisfy specific real polynomial equalities and inequalities. 

Moreover, Vinberg's results provide a method to construct the universal covering manifold of a given Coxeter $d$-orbifold. 
Using this method, we can identify the deformation space $\mathcal{C} (\mathcal{G})$ as a subspace of the deformation space of real projective structures on the Coxeter orbifold associated with $\mathcal{G}$. 
This process will be briefly mentioned in Theorem \ref{Vinberg 2} and Remark \ref{subspace} but will not be discussed in detail, as it does not appear in the subsequent discussions.

In this section, we present some of Vinberg's results, the first of which will be frequently used in our discussions. 
Additionally, we strengthen Lemma \ref{embedding lemma} by describing the image of the embedding $\Phi$ (Remark \ref{characterization}).

We use the notion of a linear reflection introduced in Section 2.2.
Let $K \subset \mathbb{R}^{d+1}$ be a convex polyhedral cone given by
\[ K := \{ x \in \mathbb{R}^{d+1} \ | \ \alpha_j (x) \geq 0, \ j = 1, \ldots, f \} \]
for some linear functionals $\alpha_1, \ldots, \alpha_f \in (\mathbb{R}^{d+1})^* \setminus \{0\}$.
Let $v_1, \ldots, v_f \in \mathbb{R}^{d + 1}$ be vectors with $\alpha_j (v_j) = 2$ and let $R_j := \text{Id} - \alpha_j \otimes v_j$ be the linear reflection determined by $\alpha_j, v_j$.
We further assume that $K$ has nonempty interior and each inequality $\alpha_j \geq 0$ is not implied by the other $(f-1)$ inequalities.

The subgroup $\Gamma$ of $\text{GL} (d+1, \mathbb{R})$ generated by the $f$ reflections $R_1, \ldots, R_f$ is called a \textit{linear reflection group generated by the reflections} $R_1, \ldots, R_f$ if we further have
\[ \gamma \cdot \text{Int} (K) \cap \text{Int} (K) = \varnothing \quad \text{for} \ \gamma \in \Gamma \setminus \{1\}. \]
Linear reflection groups are also called \textit{linear Coxeter groups or discrete linear groups generated by reflections} \cite{MR0302779}.

In this setting, we have the following theorem.

\begin{theorem} \textup{(\cite{MR0302779}, Theorems~1 and~2(6), and Propositions~6 and~17)} \label{Vinberg 1}
	Let $K$, $\alpha_j$, $v_j$ and $R_j$, $j = 1, \ldots, f$ be as above.
	Then the subgroup of $\textup{GL} (d + 1, \mathbb{R})$ generated by $R_1, \ldots, R_f$ is a linear reflection group if and only if $\alpha_1, \ldots, \alpha_f$, $v_1, \ldots, v_f$ satisfy the following conditions:
	\begin{enumerate}
		\item[\textup{(i)}] $\alpha_i (v_j) \leq 0$ if $i \ne j$;
		\item[\textup{(ii)}] $\alpha_i (v_j) = 0$ if and only if $\alpha_j (v_i) = 0$;
		\item[\textup{(iii)}] $\alpha_i (v_j) \alpha_j (v_i) \geq 4$ or $\alpha_i (v_j) \alpha_j (v_i) = 4 \cos^2 (\frac{\pi}{n_{i,j}})$ for some integer $n_{i,j} \geq 2$.
	\end{enumerate}
	In this case, the polyhedral cone $K$ is a fundamental domain of the discrete subgroup, and the subgroup is isomorphic to the abstract group
\begin{align}\label{Coxeter group}
	\left\langle \gamma_1, \ldots, \gamma_f \mid (\gamma_i \gamma_j)^{n_{i,j}} = 1 \right\rangle,
\end{align}
	via the isomorphism $\gamma_j \mapsto R_j$.
	Here, we have $n_{i, j} = 1$ if $i = j$.
	If $i \ne j$ and $\alpha_i (v_j) \alpha_j (v_i) \geq 4$, then the relation $(\gamma_i \gamma_j)^{n_{i,j}} = 1$ is omitted from the group presentation \eqref{Coxeter group}.
	
	For each $i$, let $K_i$ denote the facet of $K$ supported by $\alpha_i$.
	It suffices to verify conditions \textup{(i)}--\textup{(iii)} only for pairs of adjacent facets: the subgroup generated by $R_1,\ldots,R_f$ is a linear reflection group if and only if these conditions hold for every pair $i\ne j$ such that $K_i$ and $K_j$ are adjacent.
	Moreover, if this subgroup is a linear reflection group and $K_i$ and $K_j$ are not adjacent, then $R_iR_j$ has infinite order.
\end{theorem}

The following theorem and its subsequent Remark \ref{subspace} enable us to relate the deformation space $\mathcal{C} (\mathcal{G})$ of Coxeter polytopes realizing the labeled combinatorial polytope $\mathcal{G}$ with the deformation space of real projective structures on the Coxeter orbifold associated to $\mathcal{G}$ (see Section \ref{introduction}).

\begin{theorem} \textup{(\cite{MR0302779}, Theorem 2)}\label{Vinberg 2}
	Let $\Gamma \subset \textup{SL}_\pm(d+1, \mathbb{R})$ be a linear reflection group generated by reflections $R_1, \ldots, R_f \in  \textup{SL}_\pm(d+1, \mathbb{R})$ in the facets of a convex polyhedral cone $K$.
	For each $x \in K$, let $\Gamma_x$ denote the subgroup of $\Gamma$ generated by reflections in those facets of $K$ which contain $x$.
	Define 
	\[ K^f := \{ x \in K \ | \ \Gamma_x \ \textup{is finite} \}.  \]
	Then the following assertions are true.
	\begin{enumerate}
		\item[\textup{(i)}] $\bigcup_{\gamma \in \Gamma} \gamma K$ is a convex cone.
		\item[\textup{(ii)}] $\Gamma$ acts discretely on the interior $\widetilde{\Omega} := \textup{Int} \left( \bigcup_{\gamma \in \Gamma} \gamma K \right) $ of the cone $\bigcup_{\gamma \in \Gamma} \gamma K$.
		\item[\textup{(iii)}] $\widetilde{\Omega} \cap K = K^f$.
		\item[\textup{(iv)}] The canonical map $K^f \to \widetilde{\Omega} / \Gamma$ is a homeomorphism.
		\item[\textup{(v)}] For every $x \in K$, $\Gamma_x \subset \Gamma$ is the stabilizer of $x$ in $\Gamma$.
		\item[\textup{(vi)}] For every pair $K_i, K_j$ of facets of $K$ supported by the functionals $\alpha_i, \alpha_j$, let $n_{ij}$ denote the order of $R_i R_j$ ($n_{ij}$ may be infinite).
		Then
		\[ R_i^2 = 1, \quad (R_i R_j)^{n_{ij}} = 1 \]
		is a system of defining relations for $\Gamma$.
	\end{enumerate}
\end{theorem}

\begin{remark}\label{subspace}
	Let $K^f, \widetilde{\Omega}$ be as in the theorem.
	Let $\pi: \mathbb{R}^{d + 1} \setminus \{0\} \to \mathbb{S}^d$ be the natural projection and let $\Omega := \pi (\widetilde{\Omega})$, $\widehat{P}:= \pi (K^f)$.
Since $\bigcup_{\gamma \in \Gamma} \gamma K$ is a convex cone by Theorem~\ref{Vinberg 2}, its interior $\widetilde{\Omega}$ is also convex.
Hence its projectivization $\Omega \subset \mathbb{S}^d$ is convex in the sense defined in Section~\ref{Section 2.1}.
	Therefore, the quotient $\Omega / \Gamma$ gives a convex real projective structure on the orbifold $\widehat{P}$.
	In fact, Theorem 2 of \cite{MR2247648} states that every real projective structure on a Coxeter orbifold is convex when the orbifold is 3-dimensional.
\end{remark}

\begin{remark}\label{characterization}
	Recall from Lemma \ref{embedding lemma} that the space $\mathcal{C} (\mathcal{G})$ can be identified as a subspace of the quotient space \eqref{quotient}.
	We can improve the lemma by describing the image of the embedding $\Phi$ as follows.
	It follows from Theorem~\ref{Vinberg 1} and the definition of $\mathcal{C} (\mathcal{G})$ that
\[ [(\alpha,v)]
=
[(\alpha_1,\ldots,\alpha_f,v_1,\ldots,v_f)]
\in
\frac{\left((\mathbb R^{d+1})^*\right)^f\times(\mathbb R^{d+1})^f}
{\text{SL}_{\pm}(d+1,\mathbb R)\times\mathbb R_+^f} \]
	belongs to (the image of) $\mathcal{C} (\mathcal{G})$ if and only if
	\begin{enumerate}[(i)]
		\item $\alpha_i (v_i) = 2$ for $i = 1, \ldots, f$;
		\item $\alpha_i (v_j) \leq 0$ if $i \ne j$;
		\item $\alpha_i (v_j) = 0$ if and only if $\alpha_j (v_i) = 0$;
		\item $\alpha_i (v_j) \alpha_j (v_i) \geq 4$ if $\mathcal{G}_i$ and $\mathcal{G}_j$ are not adjacent;
		\item $\alpha_i (v_j) \alpha_j (v_i) = 4 \cos^2 \left( \frac{\pi}{m_{i,j}} \right) $ if $\mathcal{G}_i$ and $\mathcal{G}_j$ are adjacent, where $m_{i,j}$ is the exponent in \eqref{labeled polytope Coxeter group};
		\item the projectivization of the polyhedral cone $\{x \in \mathbb{R}^{d+1} \ | \ \alpha_j (x) \geq 0 \ \text{for each} \ j\}$ in $\mathbb{S}^d$ is a convex $d$-polytope representing the combinatorial polytope $\mathcal{G}$, so that the facet supported by $\alpha_j$ represents $\mathcal{G}_j$.
	\end{enumerate}
\end{remark}

\section{Restricted deformation spaces}\label{Restricted deformation spaces}

In this section, we recall the precise notion of restricted deformation spaces and the sufficient conditions under which Coxeter 3-orbifolds admit smooth restricted deformation spaces.

Given a discrete subgroup $\Gamma \subset \text{SL}_\pm (d + 1, \mathbb{R})$ preserving a convex open domain $\Omega \subset \mathbb{S}^d$, it admits fundamental domains $D \subset \Omega$ such that $\bigcup_{\gamma \in \Gamma} \gamma \cdot D = \Omega$.
Choi \cite{MR2247648} proved that if a Coxeter 3-orbifold is orderable and satisfies certain generic conditions, then the space consisting of real projective structures whose holonomy representations share a common fixed fundamental domain is a smooth manifold.
Such spaces are called the \textit{restricted deformation spaces}, which are defined precisely as follows.

Let $\mathcal G$ be a labeled combinatorial $3$-polytope with facets $\mathcal G_1,\ldots,\mathcal G_f$.
Let $P,Q\subset\mathbb S^3$ be convex $3$-polytopes realizing $\mathcal G$, and let $P_1,\ldots,P_f$ and $Q_1,\ldots,Q_f$ be their facets corresponding to $\mathcal G_1,\ldots,\mathcal G_f$, respectively.
We say that $P$ and $Q$ are projectively equivalent as realizations of $\mathcal G$ if there exists $A\in\text{SL}_{\pm}(4,\mathbb R)$ such that
\[
A(P)=Q
\qquad\text{and}\qquad
A(P_i)=Q_i
\quad\text{for every }i=1,\ldots,f.
\]
We define an equivalence relation on $\mathcal C(\mathcal G)$ by $[(P,r_1,\ldots,r_f)] \sim [(Q,s_1,\ldots,s_f)]$ if and only if $P$ and $Q$ are projectively equivalent as realizations of $\mathcal G$ in the above sense.
The resulting equivalence class is denoted by $\mathcal C_P(\mathcal G)$ and is called the \textit{restricted deformation space} with fundamental polytope $P$.
Here the subscript $P$ in $\mathcal C_P(\mathcal G)$ represents not only the polytope $P$ itself but also the marking of its facets given by the correspondences $P_i\leftrightarrow\mathcal G_i$ for $i=1,\ldots,f$; for simplicity, we record only $P$ in the notation.

The following two definitions involve our main assumptions on the labeled combinatorial polytopes $\mathcal{G}$.

\begin{definition}\label{orderable}
	A labeled combinatorial $3$-polytope $\mathcal G$ is \textit{orderable} if its facets can be indexed as $\mathcal G_1,\ldots,\mathcal G_f$ so that, for each $i$, at most three edges of $\mathcal G_i$ either have order $2$ or are of the form $\mathcal G_i\cap\mathcal G_j$ for some $j<i$.
\end{definition}

\begin{example}[A labeled combinatorial 3-polytope which is orderable]
	\label{ex:orderable-prism}
	Let $\mathcal{G}$ be the labeled combinatorial 3-polytope shown in
	Figure~\ref{fig:orderable-prism}.
	The underlying polytope is a triangular prism.
	We label the three quadrilateral side facets cyclically by
	$\mathcal{G}_1, \mathcal{G}_2, \mathcal{G}_3$, and the two triangular
	facets by $\mathcal{G}_4$ and $\mathcal{G}_5$.
	The three edges at which two quadrilateral side facets meet, namely $\mathcal{G}_1 \cap \mathcal{G}_2$, $\mathcal{G}_2 \cap \mathcal{G}_3$, $\mathcal{G}_3 \cap \mathcal{G}_1$ are assigned order $2$, and all the remaining edges are assigned
	order $3$.
	
	\begin{center}
		\begin{tikzpicture}[
			scale=0.95,
			every node/.style={inner sep=1pt}
			]
			\coordinate (A)  at (0, 0.8);
			\coordinate (B)  at (0.69282, -0.4);
			\coordinate (C)  at (-0.69282, -0.4);
			\coordinate (Ap) at (0, 2.4);
			\coordinate (Bp) at (2.1650625, -1.25);
			\coordinate (Cp) at (-2.1650625, -1.25);
			
			\draw (A) -- (B) -- (C) -- cycle;
			\draw (Ap) -- (Bp) -- (Cp) -- cycle;
			\draw (A) -- (Ap);
			\draw (B) -- (Bp);
			\draw (C) -- (Cp);
			
			\node at (-1.05, 0.08)
			{\footnotesize $\mathcal{G}_1$};
			\node at (1.05, 0.08)
			{\footnotesize $\mathcal{G}_2$};
			\node at (0, -0.92)
			{\footnotesize $\mathcal{G}_3$};
			\node at (0, 0.04)
			{\footnotesize $\mathcal{G}_4$};
			\node at (2.5, 0.08)
			{\footnotesize $\mathcal{G}_5$};
			
			\node at (-0.2, 1.30) {\tiny 2};
			\node at (1.2, -0.90) {\tiny 2};
			\node at (-1.2, -0.90) {\tiny 2};
			
			\node at (0.5, 0.25) {\tiny 3};
			\node at (0, -0.52) {\tiny 3};
			\node at (-0.5, 0.25) {\tiny 3};
			\node at (1.4, 0.35) {\tiny 3};
			\node at (0, -1.38) {\tiny 3};
			\node at (-1.4, 0.35) {\tiny 3};
		\end{tikzpicture}
		\captionof{figure}{
			An orderable labeled combinatorial 3-polytope.
		}
		\label{fig:orderable-prism}
	\end{center}
	
	This labeled combinatorial polytope is orderable.
	Indeed, we use the given indexing of the facets.
	With respect to this indexing, the numbers of edges of $\mathcal{G}_1,\ldots,\mathcal{G}_5$ counted according to Definition~\ref{orderable} are $2,2,2,3,$ and $3$, respectively.
	In particular, every facet has at most three counted edges, as required by Definition~\ref{orderable}.
\end{example}

\begin{example}[A labeled combinatorial 3-polytope that is not orderable]
	\label{ex:nonorderable-prism}
	If a facet of a labeled combinatorial 3-polytope contains at least
	four edges of order $2$, then the polytope is not orderable, since
	all of these edges are counted in Definition~\ref{orderable},
	independently of the indexing of the facets.
	In the following example, there is no such facet but the labeled combinatorial polytope is still not orderable.
	
	Let $\mathcal{H}$ be the labeled combinatorial 3-polytope shown in
	Figure~\ref{fig:nonorderable-prism}.
	Label the three quadrilateral facets cyclically by
	$\mathcal{H}_1,\mathcal{H}_2,\mathcal{H}_3$, and label the two
	triangular facets by $\mathcal{H}_4$ and $\mathcal{H}_5$.
	Assign order $2$ to the six edges $\mathcal{H}_1\cap\mathcal{H}_3$, $\mathcal{H}_2\cap\mathcal{H}_3$, $\mathcal{H}_1\cap\mathcal{H}_4$, $\mathcal{H}_2\cap\mathcal{H}_4$, $\mathcal{H}_1\cap\mathcal{H}_5$, $\mathcal{H}_2\cap\mathcal{H}_5$ and assign order $3$ to the remaining three edges.
	
	\begin{center}
		\begin{tikzpicture}[
			scale=0.95,
			every node/.style={inner sep=1pt}
			]
			\coordinate (A)  at (0, 0.8);
			\coordinate (B)  at (0.69282, -0.4);
			\coordinate (C)  at (-0.69282, -0.4);
			\coordinate (Ap) at (0, 2.4);
			\coordinate (Bp) at (2.1650625, -1.25);
			\coordinate (Cp) at (-2.1650625, -1.25);
			
			\draw (A) -- (B) -- (C) -- cycle;
			\draw (Ap) -- (Bp) -- (Cp) -- cycle;
			\draw (A) -- (Ap);
			\draw (B) -- (Bp);
			\draw (C) -- (Cp);
			
			\node at (-1.05, 0.08)
			{\footnotesize $\mathcal{H}_1$};
			\node at (1.05, 0.08)
			{\footnotesize $\mathcal{H}_2$};
			\node at (0, -0.92)
			{\footnotesize $\mathcal{H}_3$};
			\node at (0, 0.04)
			{\footnotesize $\mathcal{H}_4$};
			\node at (2.5, 0.08)
			{\footnotesize $\mathcal{H}_5$};
			
			\node at (-0.2, 1.30) {\tiny 3};
			\node at (1.2, -0.90) {\tiny 2};
			\node at (-1.2, -0.90) {\tiny 2};
			
			\node at (0.5, 0.25) {\tiny 2};
			\node at (0, -0.52) {\tiny 3};
			\node at (-0.5, 0.25) {\tiny 2};
			\node at (1.4, 0.35) {\tiny 2};
			\node at (0, -1.38) {\tiny 3};
			\node at (-1.4, 0.35) {\tiny 2};
		\end{tikzpicture}
		
		\captionof{figure}{
			A non-orderable labeled combinatorial 3-polytope with no
			facet containing more than three edges of order $2$.
		}
		\label{fig:nonorderable-prism}
	\end{center}
	
	No facet of $\mathcal{H}$ contains more than three edges of order $2$: each of $\mathcal{H}_1$ and $\mathcal{H}_2$ contains three, whereas each of $\mathcal{H}_3,\mathcal{H}_4$, and $\mathcal{H}_5$ contains two.
	Nevertheless, $\mathcal{H}$ is not orderable.
	Indeed, consider an arbitrary total ordering $<$ on the set $\{\mathcal{H}_1,\ldots,\mathcal{H}_5\}$.
	If $\mathcal{H}_2<\mathcal{H}_1$, the three order-$2$ edges $\mathcal{H}_1\cap\mathcal{H}_3$, $\mathcal{H}_1\cap\mathcal{H}_4$, and $\mathcal{H}_1\cap\mathcal{H}_5$ are counted for $\mathcal{H}_1$, as is the edge $\mathcal{H}_1\cap\mathcal{H}_2$.
	The facet $\mathcal{H}_1$ then has four counted edges.
	Similarly, if $\mathcal{H}_1<\mathcal{H}_2$, the facet $\mathcal{H}_2$ has four counted edges.
	It follows that no total ordering satisfies Definition~\ref{orderable}, so $\mathcal{H}$ is not orderable.
\end{example}

\begin{definition}\label{normal type}
	We say that $\mathcal{G}$ is of \textit{normal type} if $\mathcal{G}$ does \textit{not} satisfy any of the following conditions:
	\begin{enumerate}[(i)]
		\item $\mathcal{G}$ is a cone over a polygon and the integers assigned to the edges of the base polygon are all 2;
		\item $\mathcal{G}$ equals a polygon times a closed interval and the edge orders assigned to the edges of the two base polygons are all 2;
		\item the associated Coxeter group $\Gamma_\mathcal{G}$ (see \eqref{labeled polytope Coxeter group}) is finite.
		\item $\mathcal G$ admits an affine Coxeter group representation,
		i.e. there is a Coxeter polytope $(P,r_1,\ldots,r_f)$ realizing
		$\mathcal G$ such that, writing
		\[
		\Gamma:=\langle r_1,\ldots,r_f\rangle
		\subset \text{SL}_{\pm}(4,\mathbb R),
		\]
		there exists an affine chart $\mathbb A^3\subset\mathbb S^3$
		invariant under $r_1,\ldots,r_f$ such that
		\[
		\operatorname{Int}\left(
		\bigcup_{\gamma\in\Gamma}\gamma\cdot P
		\right)
		\subset \mathbb A^3.
		\]
	\end{enumerate}
\end{definition}

\begin{proposition}[Margulis-Vinberg \cite{MR1748082}]
	Let $\Gamma_\mathcal{G}$ be the Coxeter group associated with $\mathcal{G}$. 
	If $\Gamma_\mathcal{G}$ is irreducible (see Definition \ref{Coxeter graph}), then we have the following trichotomy
	\begin{enumerate}[(i)]
		\item \textit{spherical} (i.e., finite);
		\item \textit{affine} (i.e., infinite and virtually abelian);
		\item \textit{large} (i.e., there exists a finite-index subgroup of $\Gamma_\mathcal{G}$ admitting a surjective homomorphism onto a free group of rank $\geq 2$).
	\end{enumerate}
\end{proposition}

The following proposition provides a broad class of labeled combinatorial $3$-polytopes of normal type.
In particular, finite-volume hyperbolic Coxeter orbifolds form a class of examples, since their associated Coxeter groups are irreducible and large.

\begin{proposition}\label{prop:normal-type-criterion}
	Let $\mathcal{G}$ be a labeled combinatorial 3-polytope and let $\Gamma_\mathcal{G}$ be its associated Coxeter group. 
	If the Coxeter group $\Gamma_\mathcal{G}$ is irreducible and large, then $\mathcal{G}$ is of normal type.
\end{proposition}

\begin{proof}
	The irreducibility of $\Gamma_{\mathcal{G}}$ rules out conditions (i) and (ii) of Definition~\ref{normal type}: either condition would make the Coxeter graph of $\Gamma_{\mathcal{G}}$ disconnected.
	Moreover, the largeness of $\Gamma_{\mathcal{G}}$ implies that it is infinite; hence condition (iii) does not hold.
	
	Let $(P,r_1,\ldots,r_f)$ be a Coxeter $3$-polytope realizing $\mathcal G$, and set
	\[
	\Gamma:=\langle r_1,\ldots,r_f\rangle\subset\text{SL}_{\pm}(4,\mathbb R).
	\]
	It remains to show that the linear reflections $r_1,\ldots,r_f$ do not simultaneously preserve an affine chart of $\mathbb S^3$ containing
	\[
	\text{Int}\!\left(\bigcup_{\gamma\in\Gamma}\gamma\cdot P\right).
	\]
	
	Write $r_i(x)=x-\alpha_i(x)v_i$, where $\alpha_i\in(\mathbb R^4)^*$ and $v_i\in\mathbb R^4$ satisfy $\alpha_i(v_i)=2$, and let $A$ be the $f\times f$ matrix whose $(i,j)$-entry is $\alpha_i(v_j)$.
	Since $\Gamma_{\mathcal G}$ is irreducible and large, the matrix $A$ is of negative type in the sense of \cite[Fact~3.17]{Danciger2024}; in particular, its smallest real eigenvalue is negative.
	
	Suppose that $r_1,\ldots,r_f$ simultaneously preserve an affine chart $\Omega\subset\mathbb{S}^3$ containing
	$\text{Int} \left(\bigcup_{\gamma\in\Gamma}\gamma\cdot P\right)$.
	Let $H\subset\mathbb{R}^4$ be the linear hyperplane such that $\partial\Omega=\mathbb{S}(H)$.
	Since $\Omega$ is invariant under $r_1,\ldots,r_f$, the linear subspace $H$ is invariant under each $r_i$.
	
	For each $i$, we have either $v_i\in H$ or $H\subset\ker\alpha_i$. Indeed, suppose that $H\not\subset\ker\alpha_i$. Then there exists $x\in H$ such that $\alpha_i(x)\neq0$. Since $H$ is invariant under $r_i$, we have $\alpha_i(x)v_i=x-r_i(x)\in H$, and hence $v_i\in H$.
	
	Define
	\[
	I:=\{i\mid v_i\in H\}
	\qquad\text{and}\qquad
	J:=\{j\mid H\subset\ker\alpha_j\}.
	\]
	The sets $I$ and $J$ form a partition of $\{1,\ldots,f\}$. Indeed, the preceding argument shows that their union is $\{1,\ldots,f\}$, while they are disjoint since $\alpha_i(v_i)=2$ for every $i$.
	
	If $i\in I$ and $j\in J$, then $\alpha_j(v_i)=0$.
	By condition~(ii) of Theorem~\ref{Vinberg 1}, we also have $\alpha_i(v_j)=0$.
	Since $(P,r_1,\ldots,r_f)$ realizes $\mathcal{G}$, this implies that $m_{i,j}=2$. Therefore, the vertices corresponding to $i$ and $j$ are not joined by an edge in the Coxeter graph. Since the Coxeter graph of $\Gamma_{\mathcal{G}}$ is connected, it follows that either $I=\varnothing$ or $J=\varnothing$.
	
	Since $P$ is properly convex, the functionals $\alpha_1,\ldots,\alpha_f$ span $(\mathbb R^4)^*$, so
	\[
	\bigcap_{i=1}^f\ker\alpha_i=\{0\}.
	\]
	If $I=\varnothing$, then $H\subset\ker\alpha_i$ for every $i$, which would force $H=\{0\}$ and contradict the fact that $H$ is a hyperplane.
	Thus $J=\varnothing$, and $\text{Span}\{v_1,\ldots,v_f\}\subset H$.
	
	Since $A$ is of negative type, Theorem~3 of \cite{MR0302779} gives a vector $x=(x_1,\ldots,x_f)\in\mathbb R^f$ such that $x_i>0$ for every $i$ and every entry of $Ax$ is negative.
	Define
	\[
	v:=-\sum_{j=1}^f x_jv_j.
	\]
	For each $i$, we have
	\[
	\alpha_i(v)
	=
	-\sum_{j=1}^f x_j\alpha_i(v_j)
	=
	-(Ax)_i
	>
	0.
	\]
	These inequalities place $[v]$ in $\text{Int}(P)\subset\Omega$, while the inclusion $v\in\text{Span}\{v_1,\ldots,v_f\}\subset H$ places the same point in $\mathbb S(H)=\partial\Omega$.
	This contradiction rules out condition~\textup{(iv)} of Definition~\ref{normal type}, and $\mathcal G$ is of normal type.
\end{proof}

We present a rephrased version of Choi's result \cite{MR2247648} on the smoothness of the restricted deformation spaces.

\begin{theorem} \textup{(\cite{MR2247648}, Theorem 4)} \label{restricted}
	Let $\mathcal{G}$ be a labeled combinatorial 3-polytope.
	Let $f$ and $e$ be the numbers of facets and edges of $\mathcal G$, respectively, and let $e_2$ be the number of edges of order $2$.
	Suppose that $\mathcal{G}$ is orderable and is of normal type.
	Let $k (\mathcal{G})$ be the dimension of the stabilizer subgroup of $\text{SL}_\pm (4, \mathbb{R})$ fixing a convex 3-polytope in $\mathbb{S}^3$ combinatorially equivalent to $\mathcal{G}$.
	Then the restricted deformation space $\mathcal{C}_P (\mathcal{G})$ is a smooth manifold of dimension $3f - e - e_2 - k(\mathcal{G})$ if it is not empty.
\end{theorem}

The subgroup of $\text{Stab}(P)$ preserving every corresponding facet $P_i$ has finite index in $\text{Stab}(P)$ and has the same dimension $k(\mathcal G)$.
The following computation is stated in \cite{MR2247648} without a proof; we include a proof for completeness.

\begin{proposition}\label{prop:stabilizer-dimension}
	Let $P \subset \mathbb{S}^3$ be a convex $3$-polytope, and let
	\[
	\textup{Stab}(P)
	:= \left\{ A \in \textup{SL}_\pm(4,\mathbb{R}) \mid A(P)=P \right\}.
	\]
	Then
	\[
	\dim \textup{Stab}(P)
	=
	\begin{cases}
		3, & \text{if $P$ is a tetrahedron},\\
		1, & \text{if $P$ is a cone over a polygon other than a triangle},\\
		0, & \text{otherwise}.
	\end{cases}
	\]
	In particular, the number $k(\mathcal{G})$ in Theorem~\ref{restricted}
	depends only on the underlying polytope of $\mathcal{G}$.
\end{proposition}

\begin{proof}
	We first outline the proof.
	After passing to the identity component of $\text{Stab}(P)$, we choose four vertices of $P$ with linearly independent lifts and thereby represent every element of the identity component by a diagonal matrix.
	We then introduce a graph whose edges record which diagonal entries are forced to agree by the remaining vertices of $P$.
	If this graph has $c$ connected components, the identity component of $\text{Stab}(P)$ has dimension $c-1$.
	Finally, we interpret these connected components as the factors of a projective join decomposition of $P$.
	This shows that $c=4$ when $P$ is a tetrahedron, $c=2$ when $P$ is a cone over a polygon other than a triangle, and $c=1$ in all other cases, yielding the stated dimension formula.
	
	Let $H:=\text{Stab}(P)$. This is a closed Lie subgroup of
	$\text{SL}_\pm(4,\mathbb{R})$. Let $H^{\circ}$ be its identity
	component. Since $H$ acts on the finite set of vertices of $P$, the
	connected group $H^{\circ}$ fixes every vertex of $P$. Moreover,
	$\dim H=\dim H^{\circ}$.
	
	Choose four vertices $p_1,\ldots,p_4$ of $P$ having linearly independent
	lifts $u_1,\ldots,u_4\in\mathbb{R}^4$. Such a choice is possible since
	$P$ is $3$-dimensional. For every $A\in H^{\circ}$, there are positive
	numbers $\lambda_i(A)$ such that
	\[
	A u_i=\lambda_i(A)u_i, \qquad i=1,\ldots,4.
	\]
	Thus, with respect to the basis $(u_1,\ldots,u_4)$, every element of
	$H^{\circ}$ is diagonal.
	
	For each vertex $p$ of $P$, choose a lift of the form
	\[
	u(p)=\sum_{i=1}^4 x_i(p)u_i.
	\]
	Define a graph $\Lambda$ on the vertex set $\{1,2,3,4\}$ by joining $i$
	and $j$ whenever there is a vertex $p$ of $P$ such that
	$x_i(p)x_j(p)\neq 0$. Let $C_1,\ldots,C_c$ be the connected components
	of $\Lambda$.
	
	We now give an explicit description of $H^{\circ}$.
	If $A\in H^{\circ}$ fixes a vertex $p$, then $A u(p)$ is a positive
	scalar multiple of $u(p)$.
	It follows that $\lambda_i(A)=\lambda_j(A)$ whenever
	$x_i(p)x_j(p)\neq 0$.
	Hence the diagonal entries of $A$ are constant on each $C_a$.
	Also, for every vertex $p$ of $P$, all indices $i$ for which
	$x_i(p)\neq 0$ belong to a single component of $\Lambda$.
	Conversely, a positive diagonal matrix whose diagonal entries are
	constant on every $C_a$ fixes every vertex of $P$, and therefore
	preserves $P$.
	After imposing determinant one, these matrices form a connected
	subgroup of $H$ containing the identity.
	It follows that
	\[
	H^{\circ}
	=
	\left\{
	\operatorname{diag}(\lambda_1,\ldots,\lambda_4)
	\ \middle|\
	\begin{array}{l}
		\lambda_i>0,\quad \lambda_1\lambda_2\lambda_3\lambda_4=1,\\
		\lambda_i=\lambda_j\text{ whenever $i,j\in C_a$ for some $a$}
	\end{array}
	\right\}.
	\]
	In particular,
	\begin{equation}\label{eq:stabilizer-components}
		\dim H=c-1.
	\end{equation}
	
	We now interpret the components $C_a$ geometrically. Set
	\[
	U_a:=\text{Span}\{u_i\mid i\in C_a\}.
	\]
	Every vertex of $P$ has all of its nonzero coordinates in a single
	component $C_a$. If $P_a\subset\mathbb{S}(U_a)$ denotes the convex hull
	of the vertices supported on $C_a$, then
	\[
	\mathbb{R}^4=U_1\oplus\cdots\oplus U_c
	\quad\text{and}\quad
	P=P_1*\cdots*P_c,
	\]
	where the right-hand side denotes the projective join, namely the convex
	hull of $P_1\cup\cdots\cup P_c$. The vertices of $P_a$ are precisely
	the vertices of $P$ supported on $C_a$. Moreover,
	\[
	\dim P_a=\dim U_a-1=|C_a|-1,
	\qquad
	\sum_{a=1}^c \bigl(\dim P_a+1\bigr)=4.
	\]
	
	Suppose that $c\geq 2$. If one of the factors $P_a$ is
	$2$-dimensional, then all the remaining factors together have vector-space
	dimension one. Hence there is exactly one remaining factor, which is a
	point, and $P$ is a cone over the polygon $P_a$. If no factor is
	$2$-dimensional, then every factor is either a point or a segment. A point
	has one vertex and a segment has two vertices, so the join has exactly
	\[
	\sum_{a=1}^c \bigl(\dim P_a+1\bigr)=4
	\]
	vertices. Therefore, $P$ is a tetrahedron. Consequently, if $c\geq 2$,
	then $P$ is either a tetrahedron or a cone over a polygon.
	
	If $P$ is a tetrahedron, choose its four vertices as
	$p_1,\ldots,p_4$. Then $\Lambda$ has four connected components, so
	\eqref{eq:stabilizer-components} gives $\dim H=3$.
	
	Suppose next that $P$ is a cone over an $n$-gon with $n\geq 4$. Choose
	three vertices $p_1,p_2,p_3$ of the base polygon and let $p_4$ be the
	apex. Since no three vertices of a convex polygon are collinear, these
	four vertices have linearly independent lifts. Every base vertex has
	zero $u_4$-coordinate. Moreover, a fourth vertex of the base polygon has
	nonzero $u_1$-, $u_2$-, and $u_3$-coordinates; otherwise it would be
	collinear with two of $p_1,p_2,p_3$. Thus $\Lambda$ has exactly the two
	components $\{1,2,3\}$ and $\{4\}$. Equation
	\eqref{eq:stabilizer-components} gives $\dim H=1$.
	
	Finally, if $P$ is neither a tetrahedron nor a cone over a polygon other
	than a triangle, the preceding join argument forces $c=1$. Hence
	\eqref{eq:stabilizer-components} gives $\dim H=0$.
\end{proof}

\begin{remark} \label{rem:geomded}
	This remark and Appendix~\ref{appendix:invariant-subspaces} record a correction to the proof of Proposition~2 of \cite{MR2247648}.
	The discussion below and the Appendix are included for readers interested in this correction and may otherwise be skipped without interrupting the main development of the paper.
	
	We first recall the content of Proposition~2 in the terminology of the present paper.
	Let $\mathcal{G}$ be a labeled combinatorial $3$-polytope of normal type, and let $(P,r_1,\ldots,r_f)$ be a Coxeter $3$-polytope realizing $\mathcal{G}$.
	Consider its holonomy representation
	\[
	\rho:\Gamma_{\mathcal{G}} \to \text{SL}_{\pm}(4,\mathbb{R}), \quad \rho(\gamma_i)=r_i.
	\]
	Proposition~2 of \cite{MR2247648} asserts that there is a neighborhood of $\rho$ in $\text{Hom}(\Gamma_{\mathcal{G}},\text{SL}_{\pm}(4,\mathbb{R}))$ on which the effective conjugation action of $\text{PGL}(4,\mathbb{R})$ is proper and free.
	
	The part of the proof of \cite[Proposition~2]{MR2247648} intended to rule out holonomy-invariant projective subspaces contains errors and gaps.
	In particular, the argument does not exclude a holonomy-invariant projective plane disjoint from the Vinberg domain.
	Such a plane occurs for an affine Coxeter representation as the projective plane at infinity of its invariant affine chart; this case is excluded by Definition~\ref{normal type}\textup{(iv)}.
	In Appendix~\ref{appendix:invariant-subspaces}, we give a self-contained replacement for the part of the proof that rules out an invariant projective line, an invariant disjoint union of two projective lines, and an invariant projective plane.
	The remaining parts of the proof of \cite[Proposition~2]{MR2247648} are unchanged.
	
	Proposition~2 is used in \cite{MR2247648} to prove its Theorem~4, which is recalled in the present paper as Theorem~\ref{restricted}.
	Apart from this indirect use through Theorem~\ref{restricted}, no result in the present paper uses Proposition~2 directly.
\end{remark}

For two non-projectively-equivalent convex 3-polytopes $P, P' \subset \mathbb{S}^3$ representing a common combinatorial polytope $\mathcal{G}$, the associated restricted deformation spaces $\mathcal{C}_P (\mathcal{G})$ and $\mathcal{C}_{P'} (\mathcal{G})$ may not be homeomorphic to each other, even if $\mathcal{G}$ satisfies the hypothesis of Theorem \ref{restricted}.
It may happen that one restricted deformation space $\mathcal{C}_P (\mathcal{G})$ is empty while another restricted deformation space $\mathcal{C}_{P'} (\mathcal{G})$ is not.
Moreover, even if both $\mathcal{C}_P (\mathcal{G})$ and $\mathcal{C}_{P'} (\mathcal{G})$ are nonempty, those spaces may not be homeomorphic to each other in general.
Both phenomena occur simultaneously in the example discussed in Section~\ref{counterexample}.

Let $\mathcal G$ be a labeled combinatorial polytope.
Each element $[(P,r_1,\ldots,r_f)]\in\mathcal C(\mathcal G)$ determines the projective equivalence class of the realization $(P;P_1,\ldots,P_f)$, where $P_i$ corresponds to $\mathcal G_i$.
To keep track of these underlying realizations as the Coxeter structure varies, we recall the projective realization space in the next section.
In Section~\ref{Section 6.2}, we define the corresponding forgetful map $\mathcal C(\mathcal G)\to\mathcal{RS}(\mathcal G)$ and verify that its nonempty fibers are precisely the restricted deformation spaces.

\section{Realization spaces}\label{Realization spaces}

In this section, we recall the notion of realization spaces for convex $3$-polytopes.
A realization of a fixed combinatorial $3$-polytope $\mathcal G$ is represented by a tuple indexed by the vertices of $\mathcal G$.
This indexing determines a correspondence between the faces of $\mathcal G$ and those of the realized polytope, and this correspondence is part of the realization data.
The projective realization space is obtained by quotienting these indexed realizations by ambient projective transformations.
Steinitz \cite{steinitz1916} studied a related realization space consisting of 3-dimensional polytopes in the affine 3-space $\mathbb{A}^3$ modulo affine equivalence, proving that the realization space of each affine 3-polytope is a cell.
See also \cite{MR1482230}, which introduces the realization spaces of affine 3-polytopes and includes a proof of Steinitz's result.

For the proof of Theorem \ref{main theorem}, we require an analogous result for the realization spaces of \textit{projective} 3-polytopes. 
At the end of this section, we will recall a theorem regarding the smoothness of the realization spaces of projective 3-polytopes.

For each subset $A \subset \mathbb{S}^3$, let $\text{conv} (A) \subset \mathbb{S}^3$ denote the convex hull of $A$ in $\mathbb{S}^3$.

\begin{definition}\label{realization}
	Let $\mathcal G$ be a combinatorial 3-polytope as defined in Section~\ref{Section 2.1}.
	Let $\mathcal{V}$ be the set of vertices of $\mathcal{G}$.
	A \textit{realization of} $\mathcal{G}$ is a tuple $(p_v)_{v \in \mathcal{V}} \in (\mathbb{S}^3)^\mathcal{V}$ such that
	\begin{enumerate}[(i)]
		\item $P := \text{conv} \{ p_v \ | \ v \in \mathcal{V} \} \subset \mathbb{S}^3$ is a convex 3-polytope combinatorially equivalent to $\mathcal{G}$;
		\item for each subset $\mathcal F\subset\mathcal V$, the subset $\text{conv}\{p_v\mid v\in\mathcal F\}\subset P$ is a face of $P$ if and only if $\mathcal F$ is the set of vertices of a face of $\mathcal G$ of the same dimension.
	\end{enumerate}
\end{definition}

\begin{definition}\label{pre-realization space}
	The \textit{pre-realization space} $\overline{\mathcal{RS}} (\mathcal{G})$ is the set of all the realizations of $\mathcal{G}$.
\end{definition}

Note that the pre-realization space $\overline{\mathcal{RS}} (\mathcal{G})$ is a subset of the product $(\mathbb{S}^3)^\mathcal{V}$.
We endow $\overline{\mathcal{RS}}(\mathcal{G})$ with the subspace topology.

Note that the group $\text{SL}_\pm (4, \mathbb{R})$ of projective automorphisms of $\mathbb{S}^3$ acts on the space $\overline{\mathcal{RS}} (\mathcal{G})$ by 
\begin{align}\label{pre-realization action}
	A \cdot (p_v)_{v \in \mathcal{V}} := (A \cdot p_v)_{v \in \mathcal{V}},
\end{align}
where we identify $A \in \text{SL}_\pm (4, \mathbb{R})$ as a projective automorphism $A: \mathbb{S}^3 \to \mathbb{S}^3$.

\begin{definition}\label{projective realization}
	The \textit{projective realization space} of $\mathcal G$ is the quotient space
	\[
	\mathcal{RS}(\mathcal G):=
	\overline{\mathcal{RS}}(\mathcal G)/\text{SL}_{\pm}(4,\mathbb R).
	\]
\end{definition}

\begin{remark}\label{face correspondence of realizations}
	Although a realization is written as a tuple of vertices, it contains more information than the underlying convex polytope alone: each element also records a \textit{marking} of the faces of its underlying polytope.
	More precisely, let
	\[
	\mathbf{p} = (p_v)_{v \in \mathcal{V}}
	\in \overline{\mathcal{RS}}(\mathcal{G})
	\]
	be a realization of $\mathcal{G}$, and let
	\[
	P(\mathbf{p})
	:=
	\text{conv}\{p_v \mid v \in \mathcal{V}\}.
	\]
	The index $v$ identifies the vertex $v$ of $\mathcal{G}$ with the vertex $p_v$ of $P(\mathbf{p})$.
	By condition \textup{(ii)} of Definition~\ref{realization}, this vertex correspondence also determines a correspondence between all faces of $\mathcal{G}$ and all faces of $P(\mathbf{p})$.
	
	For notational convenience, suppose that $\mathcal{G}$ has $f$ facets, and index them as $\mathcal{G}_1, \ldots,  \mathcal{G}_f$.
	In the orderable case, we will choose this indexing to satisfy the orderability condition in Definition~\ref{orderable}.
	For each $i$, let $\mathcal{V}_i \subset \mathcal{V}$ denote the set of vertices of the facet $\mathcal{G}_i$.
	The facet of $P(\mathbf{p})$ corresponding to $\mathcal{G}_i$ is then
	\[
	P_i(\mathbf{p})
	:=
	\text{conv}\{p_v \mid v \in \mathcal{V}_i\}.
	\]
	Thus, once the facets of $\mathcal{G}$ have been indexed, the facets of every realization are uniquely indexed as $P_1(\mathbf{p}), \ldots, P_f(\mathbf{p})$.
	
	The action in \eqref{pre-realization action} preserves these
	correspondences. 
	More precisely, we have
	\[
	A\bigl(P_i(\mathbf{p})\bigr)
	=
	P_i\bigl(A \cdot \mathbf{p}\bigr)
	\]
	for every $A \in \textup{SL}_{\pm}(4,\mathbb{R})$ and every
	$i \in \{1,\ldots,f\}$.
	Passing to the quotient $\mathcal{RS}(\mathcal G)$ therefore retains the correspondence between the faces of $\mathcal G$ and those of the realized polytope.
	In particular, we do not take any further quotient by combinatorial automorphisms of $\mathcal G$.
	If this correspondence is forgotten, one obtains a different space by taking the further finite quotient by $\text{Aut}(\mathcal G)$ acting by reindexing the vertices.
	
	This convention agrees with that of Filpo Molina
	\cite{FilpoMolina}, where a polytope is treated as a point
	configuration with indexed vertices and ambient transformations in $\text{SL}_\pm (4, \mathbb{R})$ preserve the vertex indices.
\end{remark}

Fix an affine chart $\mathbb A^3\subset\mathbb S^3$ and identify it with $\mathbb R^3$ by affine coordinates.
Let $\mathcal G$ be a combinatorial $3$-polytope with vertex set $\mathcal V$.
An \textit{affine realization of} $\mathcal G$ is a tuple $(p_v)_{v\in\mathcal V}\in(\mathbb A^3)^{\mathcal V}$ such that its affine convex hull $P:=\text{conv}_{\mathbb A^3}\{p_v\mid v\in\mathcal V\}$ is a convex $3$-polytope and the assignment $v\mapsto p_v$ is a bijection onto the vertex set of $P$ that extends to an isomorphism between the face lattices of $\mathcal G$ and $P$.
The \textit{affine pre-realization space} $\overline{\mathcal{RS}}_{\mathrm{aff}}(\mathcal G)$ is the set of all such tuples, endowed with the subspace topology inherited from $(\mathbb A^3)^{\mathcal V}$.

The affine group $\text{Aff}(3,\mathbb R)$ consists of the transformations $x\mapsto Mx+b$, where $M\in\text{GL}(3,\mathbb R)$ and $b\in\mathbb R^3$.
It acts on $\overline{\mathcal{RS}}_{\mathrm{aff}}(\mathcal G)$ by
\[
(M,b)\cdot(p_v)_{v\in\mathcal V}
:=
(Mp_v+b)_{v\in\mathcal V},
\quad
M\in\text{GL}(3,\mathbb R),\quad b\in\mathbb R^3.
\]
The \textit{affine realization space} of $\mathcal G$ is the quotient
\[
\mathcal{RS}_{\mathrm{aff}}(\mathcal G)
:=
\overline{\mathcal{RS}}_{\mathrm{aff}}(\mathcal G)/
\text{Aff}(3,\mathbb R),
\]
endowed with the quotient topology.

A classical result of Steinitz \cite{steinitz1916} states that $\mathcal{RS}_{\mathrm{aff}}(\mathcal G)$ is homeomorphic to $\mathbb R^{e-6}$, where $e$ is the number of edges of $\mathcal G$.
Filpo Molina \cite{FilpoMolina} extended the result to the projective geometry.
We state the result in our terminology as follows.

\begin{theorem} \textup{(\cite{FilpoMolina}, Theorem B)}\label{realization spaces}
	Let $\mathcal{G}$ be a combinatorial $3$-polytope with $e$ edges that is not a cone over a polygonal base.
	Then the projective realization space $\mathcal{RS} (\mathcal{G})$ is a smooth manifold of dimension $e - 9$.
\end{theorem}

Note that if $P \subset \mathbb{S}^3$ is a convex 3-polytope that is not combinatorially equivalent to a cone over a polygonal base, then $e \geq 9$.
This can be verified as follows.
Let $e_j$ denote the number of edges of the $j$-th facet of $P$.
Since $e_j \geq 3$ for each $j$, we have
\[
e = \frac{1}{2} \sum_{j = 1}^{f} e_j
\geq \frac{3f}{2}.
\]
This implies that $e \geq 9$ if $f \geq 6$.
If $f = 5$, then $P$ is either a triangular prism or a cone over a quadrilateral.
Each triangular prism has exactly 9 edges, and cones over quadrilaterals are excluded by our assumption.
If $f = 4$, then $P$ is a tetrahedron, which is always a cone over a triangle.

\begin{remark}\label{realization arbitrary}
	For $d>3$, the projective realization space of a combinatorial $d$-polytope is defined as in Definitions~\ref{realization}--\ref{projective realization}, with the ambient dimension $3$ replaced by $d$ and the acting group replaced by $\text{SL}_{\pm}(d+1,\mathbb R)$.
	However, these spaces need not be topological manifolds; see \cite{MR4220994}.
	This paper concerns $3$-dimensional polytopes, and Theorem~\ref{realization spaces} shows that the projective realization spaces considered here are smooth manifolds.
\end{remark}

\section{The proof of the main theorem}\label{proof of the main theorem}

In this section, we prove Theorem \ref{main theorem}.
The starting point is to consider a map $\mathcal{C} (\mathcal{G}) \to \mathcal{RS} (\mathcal{G})$ that ``forgets'' the reflection data (explained in Section \ref{Section 6.2}).
The smoothness of $\mathcal{C} (\mathcal{G})$ will essentially follow from the smoothness of the realization space $\mathcal{RS} (\mathcal{G})$ and the nonempty fibers of the map $\mathcal{C} (\mathcal{G}) \to \mathcal{RS} (\mathcal{G})$.

The realization space $\mathcal{RS} (\mathcal{G})$, as described in Section \ref{Realization spaces}, is expressed in terms of the vertices of the polytopes.
It will be more convenient, however, to describe the realization space in terms of the facets of the polytopes.
We achieve this in Section \ref{reparametrization} by embedding $\mathcal{RS} (\mathcal{G})$ in a quotient space of $(V^*)^f$.

In Section \ref{Section 6.2}, we consider an open map $\mathcal{C} (\mathcal{G}) \to \mathcal{RS} (\mathcal{G})$ whose nonempty fibers correspond to the restricted deformation spaces discussed in Section \ref{Restricted deformation spaces}.
Section \ref{Section 6.3} uses this map to construct a smooth manifold $\widetilde{\mathcal{D}} (\mathcal{G})$.
Finally, in Section \ref{Section 6.4}, we show that the quotient of $\widetilde{\mathcal{D}} (\mathcal{G})$ by the Lie group $\text{SL}_\pm (4, \mathbb{R}) \times \mathbb{R}_+^f$ is homeomorphic to $\mathcal{C} (\mathcal{G})$, thereby inducing a smooth structure on $\mathcal{C} (\mathcal{G})$.

Throughout this section, we adopt the following settings.
Let $\mathcal G$ be an orderable labeled combinatorial $3$-polytope satisfying the conditions of Theorem~\ref{main theorem}.
In particular, the underlying combinatorial polytope of $\mathcal G$ is not a cone over a polygon.
Let $\mathcal V$ be the set of vertices of $\mathcal G$.
We use the facet indexing $\mathcal G_1,\ldots,\mathcal G_f$ introduced in Remark~\ref{face correspondence of realizations}.
Since $\mathcal G$ is orderable, this indexing is chosen to satisfy the condition in Definition~\ref{orderable}.
We denote by $m_{i,j}$ the integer $\geq 2$ assigned to the edge (= ridge) $\mathcal{G}_i \cap \mathcal{G}_j$, whenever it is defined.
Let $V = \mathbb{R}^4$ be a 4-dimensional real vector space, and let $\mathbb{S}^3 := \mathbb{S}(V)$ denote the real projective 3-sphere obtained by projectivizing $V$. 
The group of linear automorphisms of $V$ with determinant $\pm 1$ is denoted by $\text{SL}_\pm(V) = \text{SL}_\pm(4, \mathbb{R})$, which can also be viewed as the group of projective automorphisms of $\mathbb{S}^3$.
Lastly, let $\mathbb{R}_+$ be the one-dimensional group of positive real numbers, and define $\mathbb{G} := \text{SL}_\pm (V) \times \mathbb{R}_+^f$, which is a product Lie group of dimension $15 + f$.

\subsection{Reparametrizing the realization space}\label{reparametrization}

In Section~\ref{Realization spaces}, a realization of
$\mathcal G$ was described in terms of its indexed vertices.
For the arguments below, it will be more convenient to describe
the same realization in terms of its corresponding facets.
Let
\[
\mathbf p=(p_v)_{v\in\mathcal V}
\in\overline{\mathcal{RS}}(\mathcal G),
\qquad
P(\mathbf p):=
\text{conv}\{p_v\mid v\in\mathcal V\}.
\]
For each $i\in\{1,\ldots,f\}$, let $P_i(\mathbf p)$ be the facet
of $P(\mathbf p)$ corresponding to $\mathcal G_i$, as defined in
Remark~\ref{face correspondence of realizations}.
Choose a supporting functional
$\alpha_i\in V^*$ for the facet $P_i(\mathbf p)$ such that $\alpha_i$ is nonnegative on the cone over $P(\mathbf p)$.
Then the positive ray $\mathbb R_+\alpha_i$ is uniquely determined by $\mathbf p$.
Hence the realization $\mathbf p$, together with its face correspondence, determines a unique point
\[
[(\alpha_1,\ldots,\alpha_f)]
\in (V^*)^f/\mathbb R_+^f.
\]
Each element $A\in\text{SL}_{\pm}(V)$ acts on $V^*$ by
\[
A\cdot\alpha:=\alpha\circ A^{-1}.
\]
It follows that the projective class of $\mathbf p$ determines a
unique point of
\[
(V^*)^f/
\bigl(\mathbb R_+^f\times\text{SL}_{\pm}(V)\bigr)
=
(V^*)^f/\mathbb G.
\]

We now make this facet-coordinate description precise.
Let $\overline{\mathcal{RS}}(\mathcal{G})$ be the pre-realization space of $\mathcal{G}$ (see Definition \ref{pre-realization space}).
Then, $\mathcal{RS}(\mathcal{G})$ is the quotient space $\overline{\mathcal{RS}}(\mathcal{G}) / \text{SL}_\pm(V)$. 
We first construct an embedding $\overline{\mathcal{RS}}(\mathcal{G}) \to (V^*)^f / \mathbb{R}_+^f$, where $(V^*)^f / \mathbb{R}_+^f$ is the quotient space obtained by the action of the group $\mathbb{R}_+^f$ on $(V^*)^f$ given by
\[
(c_1, \ldots, c_f) \cdot (\alpha_1, \ldots, \alpha_f) := (c_1^{-1} \alpha_1, \ldots, c_f^{-1} \alpha_f).
\]

Let
\[
\mathbf p=(p_v)_{v\in\mathcal V}
\in\overline{\mathcal{RS}}(\mathcal G),
\]
and let $P:=P(\mathbf p)$. For each $i$, set
\[
P_i:=P_i(\mathbf p),
\]
where $P_i(\mathbf p)$ is the facet corresponding to
$\mathcal G_i$ as in
Remark~\ref{face correspondence of realizations}.

We choose elements $\alpha_1, \ldots, \alpha_f \in V^*$ such that for each $i$ the hyperplane $\ker \alpha_i \subset V$ is projected onto the hypersphere of $\mathbb{S}^3$ containing $P_i$.
Additionally, we impose the condition $\alpha_i(x) \geq 0$ for each $x \in V$ lifting a point of $P \subset \mathbb{S}^3$. 
The tuple $(\alpha_1, \ldots, \alpha_f) \in (V^*)^f$ is uniquely determined up to multiplication by positive real numbers. In this way, we obtain a map $\overline{\iota}: \overline{\mathcal{RS}}(\mathcal{G}) \to (V^*)^f / \mathbb{R}_+^f$.

\begin{lemma}\label{pre-realization embedding}
	The map $\overline{\iota}: \overline{\mathcal{RS}}(\mathcal{G}) \to (V^*)^f / \mathbb{R}_+^f$ is a topological embedding.
\end{lemma}

\begin{proof}
	Suppose that $\overline{\iota}(p)=\overline{\iota}(p')$ for two realizations $p=(p_v)_{v\in\mathcal V}$,
	$p'=(p'_v)_{v\in\mathcal V}$.
	Let $P:=\text{conv}\{p_v\mid v\in\mathcal V\}$, $P':=\text{conv}\{p'_v\mid v\in\mathcal V\}$.
	For each $i$, choose inward-pointing supporting functionals $\alpha_i,\alpha'_i\in V^*$ for the facets of $P$ and $P'$ corresponding to $\mathcal G_i$, respectively, such that $\alpha_i$ and $\alpha'_i$ are nonnegative on the cones over $P$ and $P'$.
	Equality of their classes in $(V^*)^f/\mathbb R_+^f$ implies that there are $c_i>0$ such that $\alpha'_i=c_i\alpha_i$.
	Replacing $\alpha'_i$ by $c_i^{-1}\alpha'_i$, we may assume that $\alpha'_i=\alpha_i$ for every $i$. 
	Hence
	\[
	P=P'
	=\{[x]\in\mathbb S^3\mid \alpha_i(x)\ge 0
	\text{ for every }i\}.
	\]
	For $v\in\mathcal V$, set $I(v):=\{i\mid v\in\mathcal G_i\}$.
	By the face correspondence in the definition of a realization, $p_v$ is the unique vertex of $P$ contained in all facets corresponding to $\{\mathcal G_i\mid i\in I(v)\}$.
	The same characterization holds for $p'_v$.
	Therefore $p_v=p'_v$ for every $v\in\mathcal V$, and $\overline{\iota}$ is injective.
	Locally, a supporting hyperplane is obtained from a fixed set of affinely independent vertices of the corresponding facet by solving a nonsingular system of linear equations. 
	Conversely, the coordinates of a vertex are obtained from a fixed set of incident supporting hyperplanes in the same way. Hence both $\overline{\iota}$ and its inverse on its image are continuous.
\end{proof}

The group $\text{SL}_\pm (V)$ acts on the space $\overline{\mathcal{RS}}(\mathcal{G})$ by the rule \eqref{pre-realization action}, and also acts on $(V^*)^f / \mathbb{R}_+^f$ by the rule
\[
A \cdot [(\alpha_1, \ldots, \alpha_f)] := [(\alpha_1 \circ A^{-1}, \ldots, \alpha_f \circ A^{-1})].
\]
It can be easily checked that the embedding $\overline{\iota}: \overline{\mathcal{RS}}(\mathcal{G}) \to (V^*)^f / \mathbb{R}_+^f$ in Lemma \ref{pre-realization embedding} is equivariant with respect to these actions. 
Therefore, we obtain an embedding
\[
\iota: \mathcal{RS}(\mathcal{G}) \to \left((V^*)^f / \mathbb{R}_+^f\right) / \text{SL}_\pm (V) \cong (V^*)^f / \mathbb{G}.
\]
(Recall $\mathbb{G} = \text{SL}_\pm (V) \times \mathbb{R}_+^f$.)

In the proof of Theorem \ref{main theorem}, we use the smoothness of the realization space $\mathcal{RS}(\mathcal{G})$.
It turns out to be convenient to use the images of the maps $\overline{\iota}$ and $\iota$ instead of their domains $\overline{\mathcal{RS}}(\mathcal{G})$ and $\mathcal{RS}(\mathcal{G})$.
Hence, we introduce the following notations. We define
\[
\overline{\mathcal{E}} (\mathcal{G}) := \overline{\iota}(\overline{\mathcal{RS}}(\mathcal{G})) \subset (V^*)^f / \mathbb{R}_+^f,
\quad \mathcal{E} (\mathcal{G}) := \iota (\mathcal{RS} (\mathcal{G})) \subset (V^*)^f / \mathbb{G}.
\]
Lastly, we define $\widetilde{\mathcal{E}} (\mathcal{G})$ to be the preimage of $\mathcal{E} (\mathcal{G})$ under the projection $(V^*)^f \to (V^*)^f / \mathbb{G}$.
We have the following commuting diagram, where the vertical maps are the quotient maps induced by the corresponding group actions.
\begin{center}
	\begin{tikzcd}
		&  & \widetilde{\mathcal{E}} (\mathcal{G}) \subset (V^*)^f \arrow[d, two heads]                 \\
		\overline{\mathcal{RS}}(\mathcal{G}) \arrow[d, two heads] \arrow[rr, "\cong"'] \arrow[rr, "\overline{\iota}"] &  & \overline{\mathcal{E}} (\mathcal{G}) \subset (V^*)^f / \mathbb{R}_+^f \arrow[d, two heads] \\
		\mathcal{RS} (\mathcal{G}) = \overline{\mathcal{RS}}(\mathcal{G}) / \text{SL}_\pm (V) \arrow[rr, "\cong"'] \arrow[rr, "\iota"] &  & \mathcal{E} (\mathcal{G}) \subset (V^*)^f / \mathbb{G}                                             
	\end{tikzcd}
\end{center}

We define the \textit{face lattice} of $\mathcal{G}$ to be the set of all faces of $\mathcal{G}$, partially ordered by inclusion, and denote it by $\text{FL}(\mathcal{G})$.
By the construction of the map $\overline{\iota}$, we have the following description of $\widetilde{\mathcal{E}} (\mathcal{G}) \subset (V^*)^f$.

\begin{lemma}\label{realization space membership}
	Let $\alpha_1, \dots, \alpha_f \in V^*$. Then $(\alpha_1, \dots, \alpha_f) \in \widetilde{\mathcal{E}}(\mathcal{G})$ if and only if the set
	\[
	Q := \{ [x] \in \mathbb{S}^3 \mid \alpha_i(x) \geq 0 \ \text{for all} \ i \}
	\]
	is a convex 3-polytope combinatorially equivalent to $\mathcal{G}$, and there exists a face-lattice isomorphism $\phi: \textup{FL}(\mathcal{G}) \to \textup{FL}(Q)$ sending each facet $\mathcal{G}_i$ to the facet of $Q$ supported by $\alpha_i$.
\end{lemma}

\begin{proof}
	The ``only if'' part follows directly from the construction of the map $\overline{\iota}$.
	
	Conversely, the face-lattice isomorphism $\phi:\text{FL}(\mathcal G) \to \text{FL}(Q)$ induces an indexing of the vertices of $Q$ by $q_v:=\phi(v)$, $v\in\mathcal V$. 
	The realization $(q_v)_{v \in \mathcal{V}}$ is then mapped to $[(\alpha_1, \dots, \alpha_f)]$ via the map $\overline{\iota}$.
\end{proof}

The next two lemmas equip $\widetilde{\mathcal E}(\mathcal G)$ with the smooth structure needed for the quotient construction.

\begin{lemma}\label{freeness lemma}
	Let $\mathcal G$ be a combinatorial 3-polytope as defined in Section~\ref{Section 2.1}, with facets $\mathcal G_1,\ldots,\mathcal G_f$.
	Suppose $\mathcal{G}$ is not a cone over a polygon.
	Let $\widetilde{\mathcal{E}} (\mathcal{G}) \subset (V^*)^f$ be the subspace constructed in the above process. 
	Then the action of $\mathbb{G}$ on $\widetilde{\mathcal{E}}(\mathcal{G})$ is free.
\end{lemma}

\begin{proof}
	Let $\alpha = (\alpha_1, \ldots, \alpha_f) \in \widetilde{\mathcal{E}}(\mathcal{G})$, and let 
	\[ g = (A, c_1, \ldots, c_f) \in \mathbb{G} = \text{SL}_\pm (V) \times \mathbb{R}_+^f. \]
	Suppose that $g \cdot \alpha = \alpha$. 
	Consider the polytope $Q = \{ [x] \in \mathbb{S}^3 \mid \alpha_i (x) \geq 0 \ \text{for all} \ i \}$, which is combinatorially equivalent to $\mathcal{G}$.
	
	Since $g \cdot \alpha = (c_1^{-1} \alpha_1 \circ A^{-1}, \ldots, c_f^{-1} \alpha_f \circ A^{-1})$, the assumption $g \cdot \alpha = \alpha$ is equivalent to
	\begin{align}\label{freeness}
		c_i^{-1} \alpha_i \circ A^{-1} = \alpha_i \quad \text{for} \ i = 1, \ldots, f.
	\end{align}
	
	First, we note that the forms $\alpha_1, \ldots, \alpha_f$ span $V^*$. 
	If they did not span $V^*$, then there would be a point $[x] \in \mathbb{S}^3$ such that $\alpha_i (x) = 0$ for all $i$, which would imply that both $[x]$ and $[-x]$ belong to $Q$, contradicting the fact that $Q$ is properly convex.
	
	We can assume, by reindexing the facets of $\mathcal{G}$ if necessary, that $\alpha_1, \ldots, \alpha_4$ are linearly independent.
	Then for each $i \in \{5, \ldots, f\}$, there are unique coefficients $d_{i,j} \in \mathbb{R}$, $j = 1,2,3,4$, such that
	\begin{align}\label{linear combination}
		\alpha_i = \sum_{j = 1}^4 d_{i,j} \alpha_j.
	\end{align}
	
	We make two observations:
	
	\begin{enumerate}[1.]
		\item For each $i\in\{5,\ldots,f\}$, at most one of the coefficients $d_{i,j}$ is zero.
		Indeed, suppose that at least two of the coefficients vanish.
		For simplicity, we assume that $d_{i,3}=d_{i,4}=0$.
		Choose a lift $x\in V$ of a point in the interior of $Q$.
		Then $\alpha_l(x)>0$ for every $l\in\{1,\ldots,f\}$.
		Therefore, $d_{i,1}$ and $d_{i,2}$ cannot both be nonpositive, since otherwise $\alpha_i(x)=d_{i,1} \alpha_1(x)+ d_{i,2} \alpha_2(x)\leq 0$.
		Thus at least one of $d_{i,1}$ and $d_{i,2}$ is positive.
		If $d_{i,1},d_{i,2}\geq 0$, then the inequalities $\alpha_1\geq 0$ and $\alpha_2\geq 0$ imply $\alpha_i\geq 0$.
		If $d_{i,1}>0$ and $d_{i,2}<0$, then $\alpha_1 = \frac{1}{d_{i,1}}\alpha_i - \frac{d_{i,2}}{d_{i,1}}\alpha_2$,
		so the inequalities $\alpha_i\geq 0$ and $\alpha_2\geq 0$ imply $\alpha_1\geq 0$.
		The case where $d_{i,1} < 0$ and $d_{i,2} > 0$ is treated similarly.
		In every case, one of the three inequalities $\alpha_i\geq 0$, $\alpha_1\geq 0$, $\alpha_2\geq 0$ is implied by the other two. 
		This contradicts the nonredundancy of the facet-defining inequalities of $Q$.
		\item For each $ j \in \{1,2,3,4\} $, there is at least one $ i \in \{5, \ldots, f\} $ such that $ d_{i,j} \neq 0 $.
		Indeed, suppose that for some $j \in \{1,2,3,4\}$, we had $d_{i,j} = 0$ for all $i \geq 5$.
		Without loss of generality, assume $ j = 1 $. Since $ \alpha_1, \ldots, \alpha_4 $ are linearly independent, there is a unique point $[x] \in \mathbb{S}^3$ such that $ \alpha_2(x) = \alpha_3(x) = \alpha_4(x) = 0 $ and $ \alpha_1(x) > 0 $. If $ d_{i,1} = 0 $ for all $ i \in \{5, \ldots, f\} $, then $ \alpha_i(x) = 0 $ for all such $ i $, implying that $ Q $ is the cone over the facet $ Q_1 $ with apex $ [x] $, contradicting our assumption.
	\end{enumerate}
	
	Now, the equations $ c_j^{-1} \alpha_j \circ A^{-1} = \alpha_j $ for $ j = 1,2,3,4 $ imply that the matrix representation of $ A $ with respect to the dual basis of $ \{\alpha_1, \ldots, \alpha_4\} $ is given by $ \text{diag}(c_1^{-1}, c_2^{-1}, c_3^{-1}, c_4^{-1}) $. Since $ A \in \text{SL}_\pm (V) $, we have $ \pm 1 = \det(A) = c_1^{-1} c_2^{-1} c_3^{-1} c_4^{-1} $.
	Since $c_j > 0$ for all $j$, we obtain $c_1 c_2 c_3 c_4 = 1$.
	
	Finally, to conclude the proof, it suffices to show that $ c_1 = c_2 = c_3 = c_4 $, since this will imply that $ A = \text{Id} $, and by equation \eqref{freeness}, $ c_j = 1 $ for all $j$.
	
	Let $i \in \{5, \ldots, f\}$.
	Applying \eqref{freeness} and \eqref{linear combination} gives
	\[ \alpha_i \circ A^{-1}
	= \sum_{j = 1}^{4} d_{i, j} \alpha_j \circ A^{-1}
	= \sum_{j = 1}^{4} d_{i, j} c_j \alpha_j. \]
	On the other hand, we have
	\[ \alpha_i \circ A^{-1}
	= c_i \alpha_i
	= c_i \sum_{j = 1}^{4} d_{i, j} \alpha_j. \]
	Comparing the coefficients, we obtain $(c_i - c_j) d_{i, j} = 0$ for all $i \in \{5, \ldots, f\}$ and $j \in \{1,2,3,4\}$.
	
	We consider the index $i = 5$.
	If $d_{5,1}, d_{5,2}, d_{5,3}, d_{5,4}$ are all nonzero, then it follows that $c_5 = c_1 = c_2 = c_3 = c_4$, so we obtain the desired result.
	Suppose that one of $d_{5,1}, d_{5,2}, d_{5,3}, d_{5,4}$ is zero.
	For simplicity, we assume that $d_{5,1} = 0$.
	Then Observation 1 implies that $d_{5,2}, d_{5,3}, d_{5,4}$ are nonzero, and this implies that $c_5 = c_2 = c_3 = c_4$.
	By Observation 2, there is another index $i \in \{6, \ldots, f\}$ such that $d_{i, 1} \ne 0$.
	Then Observation 1 again implies that there are at least two indices $u, v \in \{2,3,4\}$ such that $d_{i, u}, d_{i, v} \ne 0$, so $c_i = c_1 = c_u = c_v$.
	Hence, we obtain $c_1 = c_2 = c_3 = c_4$.
	
	We conclude that $A = \text{diag}(c_1^{-1}, c_2^{-1}, c_3^{-1}, c_4^{-1})= \text{Id}$ and $g = (\text{Id}, 1, \ldots, 1)$.
	Therefore, the action of $\mathbb{G}$ on $\widetilde{\mathcal{E}}(\mathcal{G})$ is free.
\end{proof}

\begin{lemma}\label{principal bundle}
	Let $\mathcal{G}$ be a combinatorial 3-polytope.
	Let $f$ and $e$ be the numbers of facets and edges of $\mathcal{G}$, respectively. 
	Suppose that $\mathcal{G}$ is not a cone over a polygon.
	Let $\pi: \widetilde{\mathcal{E}}(\mathcal{G}) \to \mathcal{E}(\mathcal{G})$ be the quotient map, defined as the restriction of the quotient map $(V^*)^f \to (V^*)^f / \mathbb{G}$.
	Then $\pi$ is a principal $\mathbb{G}$-bundle.
	
	In particular, $\widetilde{\mathcal{E}}(\mathcal{G})$ is a topological manifold of dimension 
	\[
	\dim \mathcal{E}(\mathcal{G}) + \dim \mathbb{G} = (e - 9) + (15 + f) = e + f + 6,
	\]
	and admits a smooth structure such that $\pi$ is a smooth submersion.
\end{lemma}

\begin{proof}
	The action of $\mathbb{G}$ on $\widetilde{\mathcal{E}}(\mathcal{G})$ is free by Lemma~\ref{freeness lemma}.
	Since $\pi$ is the orbit map, its fibers are precisely the $\mathbb{G}$-orbits. Thus $\mathbb{G}$ acts simply transitively on every fiber of $\pi$.
	
	We now outline the remaining argument.
	Around a fixed point $\alpha\in\widetilde{\mathcal E}(\mathcal G)$, we choose four linearly independent defining functionals and express each of the remaining functionals in terms of them.
	The pattern of nonzero coefficients at $\alpha$ determines a connected bipartite graph.
	A spanning tree of this graph is then used, together with the normalization of the first four functionals and the determinant condition, to select a unique representative from each $\mathbb G$-orbit in a suitable saturated neighborhood of $\alpha$.
	We prove that the normalized representative depends continuously on the original point and is constant along each orbit, so that it descends to a continuous local section of $\pi$ and yields an explicit principal-bundle trivialization.
	Finally, these local trivializations give the asserted dimension of $\widetilde{\mathcal E}(\mathcal G)$, and \cite[Proposition~II.9]{MR2574141} is used to equip it with a compatible smooth structure for which $\pi$ is a smooth principal-bundle projection.
	
	Fix $\alpha=(\alpha_1,\ldots,\alpha_f) \in \widetilde{\mathcal{E}}(\mathcal{G})$.
	The linear forms $\alpha_1,\ldots,\alpha_f$ span $V^*$, since the polytope determined by $\alpha$ is properly convex. 
	After reindexing the facets, we may therefore assume that $\alpha_1,\ldots,\alpha_4$ are linearly independent.
	
	Let $\widetilde W$ be the open subset of $\widetilde{\mathcal{E}}(\mathcal{G})$ consisting of those $\beta=(\beta_1,\ldots,\beta_f)$ for which $\beta_1,\ldots,\beta_4$ are linearly independent.
	The set $\widetilde W$ is $\mathbb{G}$-invariant. 
	For every $i\geq 5$, there are unique continuous functions $t_{i,j}: \widetilde W\to\mathbb{R}$, $j=1,2,3,4$, defined by
	\[
	\beta_i=\sum_{j=1}^{4}t_{i,j}(\beta)\beta_j.
	\]
	
	If $g=(A,c_1,\ldots,c_f)\in\mathbb{G}$ and $\gamma=g\cdot\beta$, then
	\begin{equation}\label{eq:t-transformation}
		t_{i,j}(\gamma)
		=
		\frac{c_j}{c_i}t_{i,j}(\beta).
	\end{equation}
	Indeed, equation \eqref{eq:t-transformation} follows from
	\[
	\begin{aligned}
		\gamma_i
		&=
		c_i^{-1}\beta_i\circ A^{-1}  \\
		&=
		c_i^{-1}
		\sum_{j=1}^{4}t_{i,j}(\beta)\beta_j\circ A^{-1} \\
		&=
		\sum_{j=1}^{4}
		\frac{c_j}{c_i}t_{i,j}(\beta)\gamma_j.
	\end{aligned}
	\]
	In particular, since every $c_i$ is positive, both the vanishing and the sign of each $t_{i,j}$ are invariant under the $\mathbb{G}$-action.
	
	Define a bipartite graph $\mathcal H_\alpha$ with vertex set
	\[
	\{1,2,3,4\}\sqcup\{5,\ldots,f\},
	\]
	where $j\in\{1,2,3,4\}$ and $i\in\{5,\ldots,f\}$ are joined by an
	edge if and only if $t_{i,j}(\alpha)\neq 0$.
	By Observation~1 in the proof of Lemma~\ref{freeness lemma}, every vertex $i\geq 5$ is adjacent to at
	least three vertices in $\{1,2,3,4\}$. 
	By Observation~2, every vertex in $\{1,2,3,4\}$ is adjacent to at least one vertex in $\{5,\ldots,f\}$.
	
	These two observations imply that $\mathcal H_\alpha$ is connected.
	Indeed, every connected component contains a vertex $i\geq 5$ and
	hence contains at least three vertices from $\{1,2,3,4\}$. Thus two
	distinct connected components would contain at least six distinct
	vertices from the four-element set $\{1,2,3,4\}$, which is
	impossible. Notice that this argument also covers the case $f=5$.
	
	Choose a spanning tree $T$ of $\mathcal H_\alpha$ and set
	\[
	\sigma_{i,j}
	:=
	\text{sgn}\bigl(t_{i,j}(\alpha)\bigr),
	\qquad
	(i,j)\in E(T).
	\]
	Define
	\[
	\widetilde U
	:=
	\left\{
	\beta\in\widetilde W
	\ \middle|\
	\sigma_{i,j}t_{i,j}(\beta)>0
	\text{ for every }(i,j)\in E(T)
	\right\}.
	\]
	Then $\widetilde U$ is an open neighborhood of $\alpha$.
	Equation~\eqref{eq:t-transformation}, together with the positivity
	of the numbers $c_i$, shows that $\widetilde U$ is
	$\mathbb{G}$-invariant. Hence $\widetilde U$ is saturated with respect to the quotient map $\pi$.
	
	Put $U:=\pi(\widetilde U)$.
	Since $\pi^{-1}(U)=\widetilde U$ and $\widetilde U$ is open, $U$ is an open neighborhood of
	$[\alpha]$ in $\mathcal{E}(\mathcal{G})$.
	
	We now construct a distinguished representative of every orbit in
	$\widetilde U$. Fix $\beta\in\widetilde U$. Define positive numbers
	\[
	q_1(\beta),\ldots,q_f(\beta)
	\]
	by setting $q_1(\beta)=1$ and imposing
	\begin{equation}\label{eq:q-tree-relations}
		q_i(\beta)
		=
		|t_{i,j}(\beta)|q_j(\beta)
		\qquad
		\text{for every }(i,j)\in E(T),
	\end{equation}
	where $i\geq 5$ and $j\leq 4$.
	
	Since $T$ is a tree, the equations \eqref{eq:q-tree-relations} determine all the numbers $q_r(\beta)$ uniquely. 
	More explicitly, orient $T$ away from the root vertex $1$. Along the unique path from $1$ to a vertex $r$, the number $q_r(\beta)$ is obtained recursively by multiplying or dividing by one of the nonzero functions $|t_{i,j}(\beta)|$. 
	Thus every $q_r$ is a continuous positive function on $\widetilde U$.
	
	Let $D(\beta):=\det(\beta_1,\beta_2,\beta_3,\beta_4)$ and define
	\[
	\lambda(\beta)
	:=
	\left(
	\frac{|D(\beta)|}
	{q_1(\beta)q_2(\beta)q_3(\beta)q_4(\beta)}
	\right)^{1/4},
	\qquad
	c_r(\beta):=\lambda(\beta)q_r(\beta).
	\]
	Then we have $c_1(\beta)c_2(\beta)c_3(\beta)c_4(\beta) = |D(\beta)|$.
	
	Define a linear map $A_\beta\colon V \to V$ by
	\[
	A_\beta(x)
	=
	\bigl(
	c_1(\beta)^{-1}\beta_1(x),
	c_2(\beta)^{-1}\beta_2(x),
	c_3(\beta)^{-1}\beta_3(x),
	c_4(\beta)^{-1}\beta_4(x)
	\bigr).
	\]
	We have
	\[
	\begin{aligned}
		\det A_\beta
		&=
		\frac{D(\beta)}
		{c_1(\beta)c_2(\beta)c_3(\beta)c_4(\beta)}
		=
		\frac{D(\beta)}{|D(\beta)|}
		\in\{\pm1\}.
	\end{aligned}
	\]
	Therefore we obtain $A_\beta\in\text{SL}_\pm(V)$.
	Set $g_\beta := \bigl(A_\beta,c_1(\beta),\ldots,c_f(\beta)\bigr)\in\mathbb{G}$ and $\beta':=g_\beta\cdot\beta$.
	Define $\widetilde{s}:\widetilde{U}\to\widetilde{\mathcal E}(\mathcal G)$ by
	\[
	\widetilde{s}(\beta):=\beta'=g_\beta\cdot\beta.
	\]
	
	For $j=1,2,3,4$, the definition of $A_\beta$ gives
	\[
	\beta'_j
	=
	c_j(\beta)^{-1}\beta_j\circ A_\beta^{-1}
	=
	e_j^*,
	\]
	where $e_1^*,\ldots,e_4^*$ is the standard dual basis of $V^*$.
	Moreover, by \eqref{eq:t-transformation}, we obtain
	\begin{align}\label{transition t}
		t_{i,j}(\beta')
		=
		\frac{c_j(\beta)}{c_i(\beta)}t_{i,j}(\beta).
	\end{align}
	
	For every edge $(i,j)\in E(T)$, we therefore obtain
	\begin{equation}\label{eq:tree-normalization}
		\begin{aligned}
			t_{i,j}(\beta')
			&=
			\frac{q_j(\beta)}{q_i(\beta)}t_{i,j}(\beta)
			=
			\frac{t_{i,j}(\beta)}
			{|t_{i,j}(\beta)|}
			=
			\sigma_{i,j}.
		\end{aligned}
	\end{equation}
	
	Although the normalization in \eqref{eq:tree-normalization} is imposed only on the coefficients associated with the edges of $T$, it determines the entire point $\beta'$. 
	Indeed, the tree equations determine all ratios among $c_1(\beta),\ldots,c_f(\beta)$, while the determinant condition $c_1(\beta)c_2(\beta)c_3(\beta)c_4(\beta)=|D(\beta)|$ determines their common scale. 
	Once these numbers are determined, the formula \ref{transition t} determines every remaining coefficient. Since $\beta'_1,\ldots,\beta'_4$ are fixed, these coefficients determine every $\beta'_i$ with $i\geq 5$.
	
	We next prove the uniqueness of the normalized representative.
	Suppose that $\eta$ and $\eta'$ are two points in the same
	$\mathbb{G}$-orbit satisfying
	\[
	\eta_j=\eta'_j=e_j^*
	\qquad
	(j=1,2,3,4)
	\]
	and
	\[
	t_{i,j}(\eta)
	=
	t_{i,j}(\eta')
	=
	\sigma_{i,j}
	\qquad
	((i,j)\in E(T)).
	\]
	Write
	\[
	\eta'=h\cdot\eta,
	\qquad
	h=(B,d_1,\ldots,d_f)\in\mathbb{G}.
	\]
	The first four normalization conditions (i.e., $d_j^{-1} e_j^* \circ B^{-1} = e_j^*$ for $j = 1,2,3,4$) imply $B=\text{diag}(d_1^{-1},d_2^{-1},d_3^{-1},d_4^{-1})$.
	For every $(i,j)\in E(T)$,
	equation~\eqref{eq:t-transformation} gives
	\[
	\sigma_{i,j}
	=
	t_{i,j}(\eta')
	=
	\frac{d_j}{d_i}t_{i,j}(\eta)
	=
	\frac{d_j}{d_i}\sigma_{i,j}.
	\]
	Equation~\eqref{eq:t-transformation} gives $d_i=d_j$ along every edge of $T$.
	Since $T$ is connected, all the numbers $d_r$ are equal to a common positive number $d$.
	The condition $B\in\text{SL}_{\pm}(V)$ gives $d^{-4}=|\det B|=1$, so $d=1$ and $h$ is the identity.
	Every orbit in $\widetilde U$ therefore contains a unique normalized point.
	
	The maps
	\[
	\beta\longmapsto q_r(\beta),\quad
	\lambda(\beta),\quad
	c_r(\beta),\quad
	A_\beta,\quad
	g_\beta
	\]
	are continuous, so $\widetilde s:\widetilde U\to\widetilde{\mathcal E}(\mathcal G)$ is continuous.
	If $\beta,\gamma\in\widetilde U$ lie in the same $\mathbb G$-orbit, the uniqueness of the normalized representative gives $\widetilde s(\beta)=\widetilde s(\gamma)$.
	The map $\widetilde s$ is constant on the fibers of the quotient map $\pi|_{\widetilde U}$ and descends to a continuous map $s:U\to\widetilde{\mathcal E}(\mathcal G)$.
	The identity $\pi\circ s=\text{id}_U$ shows that $s$ is a continuous local section of $\pi$.
	
	In fact, this section gives an explicit principal-bundle
	trivialization. Define
	\[
	\Phi:
	U\times\mathbb{G}
	\to
	\pi^{-1}(U),
	\qquad
	\Phi(u,h):=h\cdot s(u).
	\]
	The simple transitivity of the action on each fiber shows that
	$\Phi$ is bijective. 
	Its inverse is
	\[
	\Phi^{-1}(\beta)
	=
	\bigl(\pi(\beta),g_\beta^{-1}\bigr),
	\qquad
	\beta\in\pi^{-1}(U)=\widetilde U,
	\]
	since
	\[
	g_\beta\cdot\beta=s(\pi(\beta)).
	\]
	Both $\Phi$ and $\Phi^{-1}$ are continuous. Hence $\Phi$ is a
	$\mathbb{G}$-equivariant homeomorphism, where $\mathbb{G}$ acts on
	the second factor of $U\times\mathbb{G}$ by left multiplication.
	Since $\alpha$ was arbitrary, $\pi$ is a locally trivial principal
	$\mathbb{G}$-bundle.
	
	Finally, $\mathcal E(\mathcal G)$ is a smooth manifold of dimension $e-9$ by Theorem~\ref{realization spaces}, while $\mathbb G$ is a Lie group of dimension $15+f$.
	The local trivializations constructed above therefore show that $\widetilde{\mathcal E}(\mathcal G)$ is a topological manifold of dimension $(e-9)+(15+f)=e+f+6$.
	By \cite[Proposition~II.9]{MR2574141}, the continuous principal $\mathbb G$-bundle $\pi:\widetilde{\mathcal E}(\mathcal G)\to\mathcal E(\mathcal G)$ is continuously equivalent to a smooth principal $\mathbb G$-bundle.
	Transporting the smooth structure through this bundle equivalence, we obtain a smooth structure on $\widetilde{\mathcal E}(\mathcal G)$, compatible with its given topology, for which $\pi$ is a smooth principal-bundle projection.
\end{proof}

\subsection{The projection $\mathcal{C} (\mathcal{G}) \to \mathcal{RS} (\mathcal{G})$ and the restricted deformation spaces}\label{Section 6.2}

In this section, we consider a natural map $\mathcal{C} (\mathcal{G}) \to \mathcal{RS} (\mathcal{G})$ and observe that each nonempty fiber of the map can be identified with a restricted deformation space $\mathcal{C}_Q (\mathcal{G})$.
In view of Theorem \ref{restricted} and Theorem \ref{realization spaces}, the space $\mathcal{RS} (\mathcal{G})$ and the nonempty fibers of the map $\mathcal{C} (\mathcal{G}) \to \mathcal{RS} (\mathcal{G})$ thus admit smooth structures.

For computational purposes, especially in constructing the smooth structure on $\mathcal{C} (\mathcal{G})$, we use the embedding of $\mathcal{C} (\mathcal{G})$ onto a subspace $\mathcal{D} (\mathcal{G}) \subset ((V^*)^f \times V^f) / \mathbb{G}$ as described in Lemma \ref{embedding lemma}.
This embedding allows us to transform the map $\mathcal{C} (\mathcal{G}) \to \mathcal{RS} (\mathcal{G})$ into a topologically equivalent map $\mathcal{D} (\mathcal{G}) \to \mathcal{E} (\mathcal{G})$ given by restricting the projection $((V^*)^f \times V^f) / \mathbb{G} \to (V^*)^f / \mathbb{G}$.

In constructing the smooth structure on $\mathcal{C} (\mathcal{G})$, it will also be convenient to work with the lift $\widetilde{\mathcal{D}} (\mathcal{G}) \subset (V^*)^f \times V^f$ of $\mathcal{D} (\mathcal{G})$ under the projection $(V^*)^f \times V^f \to ((V^*)^f \times V^f) / \mathbb{G}$, in a manner analogous to our definition of $\widetilde{\mathcal{E}} (\mathcal{G})$.

We will examine the map $\widetilde{\mathcal{D}} (\mathcal{G}) \to \widetilde{\mathcal{E}} (\mathcal{G})$ and observe that the nonempty fibers of this map can also be identified with the restricted deformation spaces.

Let $\mathcal V$ denote the vertex set of $\mathcal G$. 
For a Coxeter polytope $(Q,r_1,\ldots,r_f)$ realizing $\mathcal G$, index the vertices of $Q$ as $(q_v)_{v\in\mathcal V}$ according to the given face correspondence; namely, $q_v\in Q_j$ if and only if $v\in\mathcal G_j$.
We define
\[
\kappa: \mathcal C(\mathcal G)\to
\mathcal{RS}(\mathcal G),\qquad
[(Q,r_1,\ldots,r_f)]
\mapsto
[(q_v)_{v\in\mathcal V}].
\]
To see that this map is well-defined, suppose that two representatives are related by
\[
(Q',r_1',\ldots,r_f')
=
\bigl(A(Q),Ar_1A^{-1},\ldots,Ar_fA^{-1}\bigr)
\]
for some $A\in\text{SL}_{\pm}(V)$. 
Then $A(Q_j)=Q_j'$ for every $j$. 
Consequently, for each $v\in\mathcal V$, the point $A(q_v)$ belongs precisely to the facets $Q_j'$ for which $v\in\mathcal G_j$, and hence $A(q_v)=q_v'$. 
Thus $(q_v)_{v\in\mathcal V}$ and $(q_v')_{v\in\mathcal V}$ represent the same element of $\mathcal{RS}(\mathcal G)$.

Recall that by Lemma \ref{embedding lemma}, we have an embedding $\Phi: \mathcal{C} (\mathcal{G}) \to ((V^*)^f \times V^f) / \mathbb{G}$. 
Let $\mathcal{D} (\mathcal{G}) := \Phi (\mathcal{C} (\mathcal{G}))$ denote the image of this embedding (see the diagram in Figure \ref{commuting diagram}). 
Let 
\[ \eta: (V^*)^f \times V^f \to ((V^*)^f \times V^f) / \mathbb{G} \]
be the projection induced by the action of $\mathbb{G}$ on $(V^*)^f \times V^f$, and define $\widetilde{\mathcal{D}} (\mathcal{G}) := \eta^{-1} (\mathcal{D} (\mathcal{G}))$. 
An element $(\alpha_1, \ldots, \alpha_f, v_1, \ldots, v_f) \in (V^*)^f \times V^f$ belongs to $\widetilde{\mathcal{D}} (\mathcal{G})$ if and only if it satisfies the conditions of Remark \ref{characterization}.

Consider the natural projection
\[
\rho : (V^*)^f \times V^f \to (V^*)^f,
\qquad
(\alpha_1,\ldots,\alpha_f,v_1,\ldots,v_f)
\mapsto
(\alpha_1,\ldots,\alpha_f).
\]
The map $\rho$ is $\mathbb{G}$-equivariant and maps
$\widetilde{\mathcal{D}}(\mathcal{G})$ into
$\widetilde{\mathcal{E}}(\mathcal{G})$.
We use the same symbol $\rho:\widetilde{\mathcal D}(\mathcal G) \to \widetilde{\mathcal E}(\mathcal G)$ for this restriction.
By equivariance, it induces a map
\[
\overline{\rho}: \mathcal D(\mathcal G) \to \mathcal E(\mathcal G).
\]
We verify that, under the embeddings $\Phi$ and $\iota$, this
induced map agrees with the map $\mathcal{C}(\mathcal{G}) \to \mathcal{RS}(\mathcal{G})$ defined above.
Let
\[
x=[(Q,r_1,\ldots,r_f)] \in \mathcal{C}(\mathcal{G})
\]
and write
\[
\Phi(x)
=
[(\alpha_1,\ldots,\alpha_f,v_1,\ldots,v_f)],
\qquad
r_i=\text{Id}-\alpha_i\otimes v_i,
\]
where each $\alpha_i$ is chosen to be nonnegative on the cone
over $Q$.
Then
\[
Q=
\{[y]\in\mathbb{S}^3
\mid \alpha_i(y)\geq 0
\text{ for every } i\},
\]
and the facet $Q_i$ corresponding to $\mathcal{G}_i$ is supported by the projective hyperplane $\mathbb{S}(\ker\alpha_i)$.
Consequently, if $(q_v)_{v\in\mathcal{V}}$ is the realization associated with $Q$, then, by the definition of
$\iota$,
\[
\iota([(q_v)_{v\in\mathcal{V}}])
=
[(\alpha_1,\ldots,\alpha_f)]
=
\overline{\rho}(\Phi(x)).
\]
Thus $\overline{\rho}$ forgets precisely the reflection vectors $v_1,\ldots,v_f$ while retaining the projective class of the underlying realization. 
Hence the diagram in Figure \ref{commuting diagram} commutes.

\begin{figure}[h]
	\centering
	\begin{tikzcd}
		\widetilde{\mathcal{D}}(\mathcal{G}) \subset (V^*)^f \times V^f \arrow[rr, "\rho"] \arrow[d, "\eta"] &  & \widetilde{\mathcal{E}} (\mathcal{G}) \subset (V^*)^f \arrow[d, "\pi"] \\
		\mathcal{D} (\mathcal{G}) \subset ((V^*)^f \times V^f) / \mathbb{G} \arrow[rr, "\overline{\rho}"]                                             &  & \mathcal{E} (\mathcal{G}) \subset (V^*)^f / \mathbb{G}                          \\
		\mathcal{C}(\mathcal{G})
		\arrow[u, "\cong"]
		\arrow[u, "\Phi"']
		\arrow[rr, "\kappa"]
		&
		&
		\mathcal{RS}(\mathcal{G})
		\arrow[u, "\cong"]
		\arrow[u, "\iota"']                  
	\end{tikzcd}
	\caption{The natural forgetful maps and their lifted and quotient versions.}
	\label{commuting diagram}
\end{figure}

Recall that the nonempty fibers of
\[
\kappa:
\mathcal C(\mathcal G)
\to
\mathcal{RS}(\mathcal G)
\]
are precisely the restricted deformation spaces (see Section \ref{Restricted deformation spaces}).

\begin{lemma}\label{fiber identification}
	
	Let $\alpha=(\alpha_1,\ldots,\alpha_f)\in\widetilde{\mathcal E}(\mathcal G)$ and put
	\[
	Q_\alpha := \bigl\{[x]\in\mathbb S^3 \mid \alpha_i(x)\geq 0 \text{ for every }i \bigr\}.
	\]
	If the fiber $\rho^{-1} (\alpha)$ is nonempty, then the maps $\Phi$ and $\eta$ restrict to homeomorphisms
	\[
	\mathcal C_{Q_\alpha}(\mathcal G) \xrightarrow[\cong]{\ \Phi\ } \overline{\rho}^{-1}(\pi(\alpha)) \xleftarrow[\cong]{\ \eta_\alpha\ } \rho^{-1}(\alpha),
	\]
	where $\eta_\alpha := \left.\eta \right|_{\rho^{-1}(\alpha)}$.
\end{lemma}

\begin{proof}
	Let
	\[
	x=[(Q,r_1,\ldots,r_f)]
	\quad\text{and}\quad
	x'=[(Q',r'_1,\ldots,r'_f)]
	\]
	be elements of $\mathcal C(\mathcal G)$, and let $(q_v)_{v\in\mathcal V}$ and $(q'_v)_{v\in\mathcal V}$ be their corresponding realizations.
	
	By the definition of $\kappa$, we have $\kappa(x)=\kappa(x')$ if and only if there exists $A\in\text{SL}_{\pm}(V)$ such that $A(q_v)=q'_v$ for every $v\in\mathcal V$.
	
	For each $i$, let $\mathcal V_i$ be the vertex set of the
	facet $\mathcal G_i$. By Remark~\ref{face correspondence of realizations},
	\[
	Q_i
	=
	\text{conv}
	\{q_v\mid v\in\mathcal V_i\},
	\quad
	Q'_i
	=
	\text{conv}
	\{q'_v\mid v\in\mathcal V_i\}.
	\]
	It follows that
	\[
	A(q_v)=q'_v
	\quad\text{for every }v
	\]
	implies
	\[
	A(Q)=Q'
	\quad\text{and}\quad
	A(Q_i)=Q'_i
	\quad\text{for every }i.
	\]
	
	Conversely, suppose that $A(Q)=Q'$ and $A(Q_i)=Q'_i$ for every $i$.
	For each vertex $v\in\mathcal V$, set $I(v):=\{i\mid v\in\mathcal G_i\}$.
	The face correspondence gives
	\[
	\{q_v\}
	=
	\bigcap_{i\in I(v)}Q_i,
	\qquad
	\{q'_v\}
	=
	\bigcap_{i\in I(v)}Q'_i.
	\]
	Hence $A(q_v)=q'_v$ for every $v\in\mathcal V$.
	
	Therefore, $\kappa(x)=\kappa(x')$ if and only if the polytopes $Q$ and $Q'$ with labeling of their facets $Q_i, Q_i'$ are projectively equivalent as realizations of $\mathcal G$.
	By the definition in
	Section~\ref{Restricted deformation spaces}, this is
	precisely the condition that $x$ and $x'$ belong to the
	same restricted deformation space. Thus the nonempty
	fibers of $\kappa$ are exactly the restricted deformation
	spaces.
	
	Now fix
	$\alpha\in\widetilde{\mathcal E}(\mathcal G)$.
	By Lemma~\ref{realization space membership}, the point $\pi(\alpha)\in\mathcal E(\mathcal G)$ corresponds under $\iota$ to the projective class of the realization $Q_\alpha$. 
	The equality $\iota\circ\kappa=\overline{\rho}\circ\Phi$ and the fact that $\Phi$ and $\iota$ are homeomorphisms onto their images imply that $\Phi$ restricts to a homeomorphism
	\[
	\mathcal C_{Q_\alpha}(\mathcal G)
	\to
	\overline{\rho}^{-1}(\pi(\alpha)).
	\]
	
	It remains to consider $\eta_\alpha$. Set
	\[
	O_\alpha
	:=
	\pi^{-1}(\pi(\alpha))
	=
	\mathbb G\cdot\alpha \subset \widetilde{\mathcal{E}} (\mathcal{G}).
	\]
	Since $\pi:\widetilde{\mathcal E}(\mathcal G)\to\mathcal E(\mathcal G)$ is a principal $\mathbb G$-bundle by
	Lemma~\ref{principal bundle}, the orbit map
	\[
	o_\alpha:
	\mathbb G
	\to
	O_\alpha,
	\qquad
	g\mapsto g\cdot\alpha,
	\]
	is a homeomorphism.
	
	Consider
	\[
	\Theta:
	\mathbb G\times\rho^{-1}(\alpha)
	\to
	\rho^{-1}(O_\alpha),
	\qquad
	(g,z)\mapsto g\cdot z.
	\]
	This map is a homeomorphism. Indeed, for
	$z\in\rho^{-1}(O_\alpha)$, define
	\[
	g_z
	:=
	o_\alpha^{-1}(\rho(z)).
	\]
	Then
	\[
	\Theta^{-1}(z)
	=
	\bigl(
	g_z,
	g_z^{-1}\cdot z
	\bigr),
	\]
	which depends continuously on $z$.
	
	Under $\Theta$, the action of $\mathbb G$ on
	$\rho^{-1}(O_\alpha)$ corresponds to left translation on
	the first factor:
	\[
	h\cdot(g,z)
	=
	(hg,z).
	\]
	Consequently,
	\[
	\rho^{-1}(O_\alpha)/\mathbb G
	\cong
	\rho^{-1}(\alpha).
	\]
	
	On the other hand, the commutativity of the diagram gives
	\[
	\eta^{-1}
	\bigl(
	\overline{\rho}^{-1}(\pi(\alpha))
	\bigr)
	=
	\rho^{-1}(O_\alpha).
	\]
	The set on the right is $\mathbb G$-invariant, and the
	restriction of the open quotient map $\eta$ to this
	saturated subset is again a quotient map. Hence
	\[
	\rho^{-1}(O_\alpha)/\mathbb G
	\cong
	\overline{\rho}^{-1}(\pi(\alpha)).
	\]
	Under these identifications, the resulting homeomorphism
	is precisely $\eta_\alpha$.
\end{proof}

In the remainder of this section, we use the homeomorphisms in Lemma~\ref{fiber identification} to identify these three spaces whenever the corresponding fibers are nonempty.

\begin{proposition}\label{prop:embedded-fibers}
	Let $\alpha=(\alpha_1,\ldots,\alpha_f)\in\widetilde{\mathcal{E}}(\mathcal{G})$, and suppose that the fiber $\rho^{-1}(\alpha)$ is nonempty.
	Identify $\rho^{-1}(\alpha)\subset\{\alpha\}\times V^f$ with a subset of $V^f$ via the natural homeomorphism $\{\alpha\}\times V^f\to V^f$.
	Then $\rho^{-1}(\alpha)$ is an embedded submanifold of $V^f$, and its embedded-submanifold smooth structure agrees with the smooth structure on the fiber obtained in the proof of Theorem~\ref{restricted}.
	
	Moreover, let $\beta\in\widetilde{\mathcal{E}}(\mathcal{G})$ be such that $\rho^{-1}(\beta)$ is nonempty.
	Suppose that $U,U'\subset V^f$ are open, that $F:U\to U'$ is smooth, and that
	\[
	\rho^{-1}(\alpha)\subset U,
	\qquad
	\rho^{-1}(\beta)\subset U',
	\qquad
	F\bigl(\rho^{-1}(\alpha)\bigr)\subset\rho^{-1}(\beta).
	\]
	Then the restriction
	\[
	F|_{\rho^{-1}(\alpha)}:
	\rho^{-1}(\alpha)\to\rho^{-1}(\beta)
	\]
	is smooth.
\end{proposition}

\begin{proof}
	Since the convex polytope determined by $\alpha$ is properly convex, the linear functionals $\alpha_1,\ldots,\alpha_f$ span $V^*$.
	Hence the linear map
	\[
	D_\alpha:V\to\mathbb{R}^f,
	\qquad
	D_\alpha(x)
	=
	\bigl(\alpha_1(x),\ldots,\alpha_f(x)\bigr),
	\]
	is injective.
	
	For each $j\in\{1,\ldots,f\}$, set
	\[
	A_j:=\{x\in V\mid\alpha_j(x)=2\},
	\]
	and let
	\[
	A:=A_1\times\cdots\times A_f\subset V^f.
	\]
	Then $A$ is an affine subspace of $V^f$, and condition $\alpha_j(v_j)=2$ gives $\rho^{-1}(\alpha)\subset A$.
	
	Let $p_j:\mathbb{R}^f\to\mathbb{R}^{f-1}$ be the coordinate projection that deletes the $j$-th coordinate, and define
	\[
	\begin{aligned}
		\Psi_\alpha:A
		&\to
		(\mathbb{R}^{f-1})^f
		\cong\mathbb{R}^{f(f-1)},\\
		(v_1,\ldots,v_f)
		&\mapsto
		\bigl(
		p_1(D_\alpha(v_1)),\ldots,
		p_f(D_\alpha(v_f))
		\bigr).
	\end{aligned}
	\]
	This map is an affine embedding.
	Indeed, the coordinate deleted by $p_j$ is identically equal to $2$ on $D_\alpha(A_j)$, and $D_\alpha$ is injective.
	It follows that
	\[
	B_\alpha:=\Psi_\alpha(A)
	\]
	is a $3f$-dimensional affine subspace of $\mathbb{R}^{f(f-1)}$.
	
	The coordinates of $\Psi_\alpha(v_1,\ldots,v_f)$ are precisely the values $\alpha_i(v_j)$ indexed by the ordered pairs $(i,j)$ with $i\neq j$.
	Thus
	\[
	S_\alpha:=
	\Psi_\alpha\bigl(\rho^{-1}(\alpha)\bigr)
	\subset B_\alpha
	\]
	is the solution set of the conditions in Remark~\ref{characterization}.
	Viewed inside $\mathbb{R}^{f(f-1)}$, it is described by the Vinberg polynomial equations together with the affine equations defining $B_\alpha$, subject to the Vinberg inequalities.
	The full-rank computation in the proof of \cite[Theorem~4]{MR2247648}, when expressed in these coordinates, shows that $S_\alpha$ is an embedded submanifold of $\mathbb{R}^{f(f-1)}$.
	
	Since $S_\alpha\subset B_\alpha$, it is also an embedded submanifold of the affine subspace $B_\alpha$.
	Since $\Psi_\alpha:A\to B_\alpha$ is an affine diffeomorphism, it follows that $\rho^{-1}(\alpha)$ is an embedded submanifold of $A$, and hence of $V^f$.
	The same argument applied to $\beta$ shows that $\rho^{-1}(\beta)$ is an embedded submanifold of $V^f$.
	The resulting embedded-submanifold structures are precisely the smooth structures on the fibers obtained in the proof of Theorem~\ref{restricted}.
	
	Finally, the composite
	\[
	\rho^{-1}(\alpha)
	\hookrightarrow U
	\xrightarrow{F}
	U'
	\hookrightarrow V^f
	\]
	is smooth and has image contained in the embedded submanifold $\rho^{-1}(\beta)$.
	By the defining property of the smooth structure on an embedded submanifold, it follows that
	\[
	F|_{\rho^{-1}(\alpha)}:
	\rho^{-1}(\alpha)\to\rho^{-1}(\beta)
	\]
	is smooth.
\end{proof}

\subsection{Construction of a smooth structure on $\widetilde{\mathcal{D}} (\mathcal{G})$}\label{Section 6.3}

In this section, we construct a smooth atlas on $\widetilde{\mathcal{D}}(\mathcal{G})$.
We divide the argument into five parts.
In Section~\ref{subsubsec:defining-data}, we record the conditions defining $\widetilde{\mathcal{D}}(\mathcal{G})$ and choose the indices supplied by the orderability condition.
In Section~\ref{subsubsec:standard-chart-construction}, we use these data to transport reflection vectors from a fixed fiber of $\rho$ to nearby fibers and construct the standard charts.
In Section~\ref{subsubsec:standard-chart-topology}, we prove that each standard chart is a topological embedding onto an open subset of $\widetilde{\mathcal{D}}(\mathcal{G})$.
In Section~\ref{subsubsec:standard-chart-smoothness}, we prove that the transition maps are smooth and compute the dimension of $\widetilde{\mathcal{D}}(\mathcal{G})$ when it is nonempty.
Finally, in Section~\ref{subsubsec:forgetful-openness}, we use the local product form of the standard charts to prove that the forgetful maps in Figure~\ref{commuting diagram} are open.
The resulting smooth structure, together with Lemma~\ref{proper and free}, will be used in Section~\ref{Section 6.4} to define a smooth structure on $\mathcal{C}(\mathcal{G})$.

\subsubsection{The defining conditions and the choice of indices}\label{subsubsec:defining-data}

We first record the conditions that characterize $\widetilde{\mathcal{D}}(\mathcal{G})$ as a subspace of $(V^*)^f \times V^f$.
We then use the orderability condition to choose, for each facet $\mathcal{G}_i$, three adjacent facets so that the prescribed Vinberg relations can be completed to a system of four linear equations for the corresponding reflection vector.
The key point is that the four coefficient functionals in each such system form a basis of $V^*$.

Recall that the subspace $\widetilde{\mathcal{D}} (\mathcal{G}) \subset (V^*)^f \times V^f$ is defined as the preimage $\eta^{-1} (\mathcal{D} (\mathcal{G}))$ under the natural projection 
\[ \eta: (V^*)^f \times V^f \to ((V^*)^f \times V^f) / \mathbb{G}. \]
Moreover, the space $\mathcal{D} (\mathcal{G}) \subset ((V^*)^f \times V^f) / \mathbb{G}$ is defined as the image of the embedding $\Phi: \mathcal{C} (\mathcal{G}) \to ((V^*)^f \times V^f) / \mathbb{G}$ from Lemma \ref{embedding lemma}.
According to Theorem \ref{Vinberg 1} and Lemma \ref{realization space membership}, an element 
\[ (\alpha, v) = (\alpha_1, \ldots, \alpha_f, v_1, \ldots, v_f) \in (V^*)^f \times V^f \]
belongs to $\widetilde{\mathcal{D}} (\mathcal{G})$ if and only if it satisfies the following conditions:
\begin{enumerate}[(V1)]
	\item \label{V1} $\alpha_i (v_i) = 2$ for $i = 1, \ldots, f$;
	\item \label{V2} $\alpha_i (v_j) \leq 0$ if $i \ne j$;
	\item \label{V3} $\alpha_i (v_j) = 0$ if and only if $\alpha_j (v_i) = 0$;
	\item \label{V4} if $\mathcal{G}_i$ and $\mathcal{G}_j$ are adjacent and $m_{i, j}$ is the integer assigned to the edge $\mathcal{G}_i\cap \mathcal{G}_j$, then $\alpha_i (v_j) \alpha_j (v_i) = 4 \cos^2 \left( \frac{\pi}{m_{i, j}} \right)$;
	\item \label{V5} if $\mathcal{G}_i$ and $\mathcal{G}_j$ are not adjacent, then $\alpha_i (v_j) \alpha_j (v_i) \geq 4$;
	\item \label{V6} the set
	\[
	Q := \{ [x] \in \mathbb{S}^3 \mid \alpha_i(x) \geq 0 \ \text{for all} \ i \}
	\]
	is a convex 3-polytope combinatorially equivalent to $\mathcal{G}$, and there exists a face-lattice isomorphism $\textup{FL}(\mathcal{G}) \to \textup{FL}(Q)$ mapping each facet $\mathcal{G}_i$ to the facet of $Q$ supported by $\alpha_i$.
\end{enumerate}

We now choose the indices used in the recursive chart construction by applying the orderability condition.
Since $\mathcal{G}$ is orderable, we reindex the facets $\mathcal{G}_1, \ldots, \mathcal{G}_f$ if necessary, so that for each facet $\mathcal{G}_i$, there are at most three facets $\mathcal{G}_j$ adjacent to $\mathcal{G}_i$ such that either:
\begin{enumerate}[-]
	\item $j < i$, or
	\item $m_{i, j} = 2$.
\end{enumerate}

Then, for each $i \in \{1, \ldots, f\}$, we choose three distinct indices $k_{i, 1}, k_{i, 2}, k_{i, 3} \in \{1, \ldots, f\}$ such that:
\begin{enumerate}[(i)]
	\item For $j = k_{i, 1}, k_{i, 2}, k_{i, 3}$, the facets $\mathcal{G}_i$ and $\mathcal{G}_j$ are adjacent.
	\item The first $a_i$ indices $j \in \{k_{i, 1}, \ldots, k_{i, a_i}\}$ satisfy $m_{i, j} = 2$.
	\item The next $b_i$ indices $j \in \{k_{i, a_i + 1}, \ldots, k_{i, a_i + b_i}\}$ satisfy $m_{i, j} \geq 3$ and $j < i$.
	\item The last $3 - a_i - b_i$ indices $j \in \{k_{i, a_i + b_i + 1}, \ldots, k_{i, 3}\}$ satisfy $m_{i, j} \geq 3$ and $j > i$.
	\item If $\mathcal{G}_i$ and $\mathcal{G}_j$ are adjacent and satisfy either $j < i$ or $m_{i, j} = 2$, then $j$ is included in $\{k_{i, 1}, k_{i, 2}, k_{i, 3}\}$.
\end{enumerate}

\begin{remark}\label{rem:selected-basis}
	
	For every $\alpha=(\alpha_1,\ldots,\alpha_f) \in \widetilde{\mathcal{E}}(\mathcal{G})$ and every \(i\in\{1,\ldots,f\}\), the four linear functionals $\alpha_i, \alpha_{k_{i,1}},\alpha_{k_{i,2}},\alpha_{k_{i,3}}$ form a basis of $V^*=(\mathbb{R}^4)^*$.
	Indeed, let
	\[
	P := \bigl\{[x]\in\mathbb{S}^3 \mid \alpha_l(x)\geq 0 \text{ for every } l \bigr\},
	\]
	and, for each $l$, let $P_l$ be the facet of $P$ supported by $H_l:=\mathbb{S}(\ker\alpha_l)$.
	By condition \textup{(i)} above, each of the facets $P_{k_{i,1}},P_{k_{i,2}},P_{k_{i,3}}$ is adjacent to
	$P_i$. 
	Hence $E_r:=P_i\cap P_{k_{i,r}}$, $r=1,2,3$ are three distinct edges of the convex polygon $P_i$.
	The corresponding supporting projective lines of $P_i$ are $L_r:=H_i\cap H_{k_{i,r}}$, $r=1,2,3$.

	We claim that $L_1\cap L_2\cap L_3=\varnothing$.
	Suppose that $q$ belongs to this intersection.
	If $q\in P_i$, then $q$ belongs to the three distinct edges $E_1,E_2,E_3$, which is impossible for a convex polygon.
	If $q\notin P_i$, choose a projective line $L$ in $H_i$ which passes through $q$ and is disjoint from $P_i$ as the line at infinity.
	In the resulting affine chart of $\mathbb{S} (\ker \alpha_i) \setminus L$ which contains $P_i$, the subset $P_i$ is a bounded convex
	polygon and $L_1,L_2,L_3$ become three distinct parallel supporting lines, since their common intersection point lies at infinity. 
	This is again impossible, since a bounded convex polygon has at most two parallel supporting lines.

	Consequently, we have $H_i \cap H_{k_{i,1}} \cap H_{k_{i,2}} \cap H_{k_{i,3}} = \varnothing$.
	In the vector space $V$ this means $\ker\alpha_i \cap\ker\alpha_{k_{i,1}} \cap\ker\alpha_{k_{i,2}} \cap\ker\alpha_{k_{i,3}} = \{0\}$.
	Since $\dim V^* = 4$, the four functionals $\alpha_i,\alpha_{k_{i,1}},\alpha_{k_{i,2}}, \alpha_{k_{i,3}}$ form a basis of $V^*$.

\end{remark}

Since $\mathcal{G}$ is orderable, for each $i$, there are at most three facets $\mathcal{G}_j$ adjacent to $\mathcal{G}_i$ such that either $m_{i, j} = 2$ or $j < i$.
Therefore, the indices $k_{i, 1}, \ldots, k_{i, a_i + b_i}$ with $a_i + b_i \leq 3$ in the above requirement do exist.
Note that the requirement (v) indicates that $k_{i, 1}, \ldots, k_{i, a_i + b_i}$ include \textit{all} such indices.
For the remaining facets $\mathcal{G}_j$ adjacent to $\mathcal{G}_i$, we must have $j > i$.
In (iv), we arbitrarily select $3 - a_i - b_i$ indices $k_{i, a_i + b_i + 1}, \ldots, k_{i, 3}$ with this property.

We make such a choice due to a technical reason, and the choice of indices $k_{i, 1}, k_{i, 2}, k_{i, 3}$ achieves the following.

Given $\beta=(\beta_1,\ldots,\beta_f) \in\widetilde{\mathcal{E}}(\mathcal{G})$, we will construct $w_1',\ldots,w_f'$ inductively so that $(\beta,w')$ satisfies conditions \hyperref[V1]{(V1)}, \hyperref[V3]{(V3)}, and \hyperref[V4]{(V4)}.
Suppose that $w_1',\ldots,w_{i-1}'$ have already been constructed. 
These conditions prescribe the values $\beta_j(w_i')$ for $j\in \{i,k_{i,1},\ldots,k_{i,a_i+b_i}\}$.
For example, if $j\in \{k_{i,a_i+1},\ldots,k_{i,a_i+b_i}\}$, then condition \hyperref[V4]{(V4)} gives
\[
\beta_j(w_i')
=
\frac{4\cos^2\left(\frac{\pi}{m_{i,j}}\right)}
{\beta_i(w_j')},
\]
where $w_j'$ has already been determined since $j<i$.
Since $a_i+b_i\leq 3$, these conditions provide at most four linear equations for $w_i'$.
We supplement them by prescribing the values $\beta_j(w_i')$, $j\in \{k_{i,a_i+b_i+1},\ldots,k_{i,3}\}$.
By Remark~\ref{rem:selected-basis}, the four coefficient functionals $\beta_i$, $\beta_{k_{i,1}}$, $\beta_{k_{i,2}}$, $\beta_{k_{i,3}}$ form a basis of $V^*$. 
Therefore, the resulting system has a unique solution $w_i'$.

\subsubsection{Construction of the standard charts}\label{subsubsec:standard-chart-construction}

Fix $(\alpha,v)=(\alpha_1,\ldots,\alpha_f,v_1,\ldots,v_f) \in \widetilde{\mathcal{D}}(\mathcal{G})$.
Here we construct a local product map around $(\alpha,v)$ by transporting each $w \in \rho^{-1}(\alpha)$ to a tuple $w'(\beta,w) \in \rho^{-1}(\beta)$ for every nearby $\beta \in \widetilde{\mathcal{E}}(\mathcal{G})$.
The transport is defined recursively using the indices chosen in Section~\ref{subsubsec:defining-data}, and Lemma~\ref{standard chart lemma} establishes its existence, uniqueness, and continuity.

The desired transport map is given by the following lemma.

\begin{lemma}\label{standard chart lemma}
	Let $(\alpha, v) \in \widetilde{\mathcal{D}} (\mathcal{G})$. 
	Then there exist sufficiently small neighborhoods $U_\alpha \subset \widetilde{\mathcal{E}} (\mathcal{G})$ and $W_v \subset \rho^{-1} (\alpha)$ of $\alpha$ and $v$, respectively, such that there is a unique continuous map $w' : U_\alpha \times W_v \to V^f$ satisfying the following properties for all $\beta \in U_\alpha$ and $w \in W_v$:
	\begin{enumerate}
		\item[(i)] $(\beta, w') \in \widetilde{\mathcal{D}} (\mathcal{G})$,
		\item[(ii)] $w'(\alpha, w) = w$,
		\item[(iii)] $\beta_j (w_i') = \alpha_j (w_i)$ for $i \in \{1, \ldots, f\}$, $j \in \{k_{i, a_i + b_i + 1}, \ldots, k_{i, 3}\}$.
	\end{enumerate}
\end{lemma}

\begin{proof}
	For $(\beta, w) \in U_\alpha \times W_v$, write
	\[
	w'(\beta,w)
	=
	\bigl(
	w_1'(\beta,w),\ldots,w_f'(\beta,w)
	\bigr).
	\]
	We first outline the proof.
	The argument consists of three steps.
	First, we construct $w_1',\ldots,w_f'$ inductively.
	At the $i$-th step, condition~\textup{(iii)} together with conditions \hyperref[V1]{(V1)}, \hyperref[V3]{(V3)}, and \hyperref[V4]{(V4)} determines a system of four linear equations for $w_i'$.
	By Remark~\ref{rem:selected-basis}, the coefficient functionals of this system form a basis of $V^*$, so the system has a unique solution depending continuously on $(\beta,w)$.
	The same induction gives $w'(\alpha,w)=w$ and proves the uniqueness of $w'$.
	Second, we shrink $U_\alpha$ and $W_v$ so that the values $\beta_p(w_q')$ remain strictly negative whenever $p\neq q$ and $\mathcal G_p\cap\mathcal G_q$ is not an edge of order $2$.
	This establishes condition \hyperref[V2]{(V2)} and completes the verification of condition \hyperref[V3]{(V3)}, while conditions \hyperref[V1]{(V1)}, \hyperref[V4]{(V4)}, and \hyperref[V6]{(V6)} follow directly from the construction and the fact that $\beta\in\widetilde{\mathcal E}(\mathcal G)$.
	Finally, Vinberg's criterion shows that the constructed reflections generate a linear reflection group, and the infinite order of the products corresponding to nonadjacent facets then gives condition \hyperref[V5]{(V5)}.
	Together, these steps show that $(\beta,w'(\beta,w))$ belongs to $\widetilde{\mathcal D}(\mathcal G)$ and establish all the assertions of the lemma.
	
	We proceed inductively on $i$.
	Suppose that the continuous functions $w_1',\ldots,w_{i-1}'$ have already been defined and satisfy
	\[
	w_j'(\alpha,w)=w_j \qquad \text{for every $w\in W_v$ and $j<i$}.
	\]
	For each $j\in \{k_{i,a_i+1},\ldots,k_{i,a_i+b_i}\}$, we have $j<i$ and $m_{i,j}\geq 3$.
	Since $(\alpha,w)\in\widetilde{\mathcal D}(\mathcal G)$ for every $w\in W_v$, conditions \hyperref[V2]{(V2)} and \hyperref[V4]{(V4)} give $\alpha_i(w_j)<0$.
	In particular, $\alpha_i(v_j)<0$.
	By the induction hypothesis and continuity, after shrinking
	$U_\alpha$ and $W_v$ if necessary, we may assume that $\beta_i\bigl(w_j'(\beta,w)\bigr)<0$ for every such $j$ and every $(\beta,w)\in U_\alpha\times W_v$.
	
	We define $w_i'(\beta,w)\in V$ to be the unique solution of the system
	\begin{equation}\label{system}
		\begin{aligned}
			\beta_i(w_i')
			&=2,\\
			\beta_j(w_i')
			&=0
			&&
			\text{for }
			j\in\{k_{i,1},\ldots,k_{i,a_i}\},\\
			\beta_j(w_i')
			&=
			\frac{
				4\cos^2\left(\frac{\pi}{m_{i,j}}\right)
			}{
				\beta_i(w_j')
			}
			&&
			\text{for }
			j\in
			\{k_{i,a_i+1},\ldots,k_{i,a_i+b_i}\},\\
			\beta_j(w_i')
			&=\alpha_j(w_i)
			&&
			\text{for }
			j\in
			\{k_{i,a_i+b_i+1},\ldots,k_{i,3}\}.
		\end{aligned}
	\end{equation}
	By the remark preceding the lemma, the four coefficient
	functionals
	\[
	\beta_i,\quad
	\beta_{k_{i,1}},\quad
	\beta_{k_{i,2}},\quad
	\beta_{k_{i,3}}
	\]
	form a basis of $V^*$.
	Therefore, the system \eqref{system} has a unique solution.
	Its right-hand side depends continuously on $(\beta,w)$, so
	Cramer's rule shows that $w_i'$ is continuous.
	
	When $\beta=\alpha$, the vector $w_i$ is a solution for the system \eqref{system}.
	Indeed, the first equation follows from \hyperref[V1]{(V1)}.
	The second group follows from \hyperref[V3]{(V3)} and \hyperref[V4]{(V4)}, since the corresponding edges have order $2$.
	The third group follows from \hyperref[V4]{(V4)} and the induction hypothesis $w_j'(\alpha,w)=w_j$.
	The equations in the last group are satisfied since, after setting $\beta=\alpha$ and $w_i'=w_i$, each of them reduces to the identity $\alpha_j(w_i)=\alpha_j(w_i)$.
	Hence
	\[
	w_i'(\alpha,w)=w_i
	\qquad
	\text{for every $w\in W_v$}.
	\]
	This completes the induction and constructs a continuous map
	\[
	w'\colon U_\alpha\times W_v\to V^f
	\]
	satisfying condition \textup{(iii)} and, in fact,
	\[
	w'(\alpha,w)=w
	\qquad
	\text{for every $w\in W_v$}.
	\]
	In particular, $w'(\alpha,v)=v$.
	
	The same induction proves uniqueness.
	Indeed, any map satisfying condition \textup{(iii)} and conditions \hyperref[V1]{(V1)}, \hyperref[V3]{(V3)}, and \hyperref[V4]{(V4)} must satisfy the system \eqref{system} for each $i$.
	For $j\in\{k_{i,1},\ldots,k_{i,a_i}\}$, the edge order is $2$; hence \hyperref[V4]{(V4)} and \hyperref[V3]{(V3)} force both corresponding values to vanish.
	For $j\in\{k_{i,a_i+1},\ldots,k_{i,a_i+b_i}\}$, condition \hyperref[V4]{(V4)} gives the third group of equations, while condition \textup{(iii)} gives the last group.
	Since the coefficient functionals form a basis, each $w_i'$ is uniquely determined once $w_1',\ldots,w_{i-1}'$ are known.
	
	We now verify conditions \hyperref[V1]{(V1)} through
	\hyperref[V6]{(V6)} for the pairs $(\beta, w') \in U_\alpha \times W_v$, shrinking $U_\alpha$ and $W_v$ once
	more when necessary.
	
	Condition \hyperref[V1]{(V1)} is immediate from the
	construction, since
	\[
	\beta_i(w_i')=2
	\qquad
	\text{for every }i\in\{1,\ldots,f\}.
	\]
	
	Condition \hyperref[V4]{(V4)} is also immediate from the
	construction.
	Indeed, if the edge $\mathcal G_i\cap\mathcal G_j$ has order $2$, then the second group of equations in \eqref{system} gives
	\[
	\beta_i(w_j')=\beta_j(w_i')=0.
	\]
	(Note that this does not complete the proof of \hyperref[V3]{(V3)}; it remains to show that if $\mathcal{G}_i \cap \mathcal{G}_j$ is not an edge of order 2, the values $\beta_i (w_j')$ and $\beta_j (w_i')$ are nonzero.)
	If its order $m_{i,j}$ is at least $3$, assume, after
	interchanging $i$ and $j$ if necessary, that $j<i$.
	Then the defining equation for $w_i'$ gives
	\[
	\beta_j(w_i')\beta_i(w_j')
	=
	4\cos^2\left(\frac{\pi}{m_{i,j}}\right).
	\]
	
	We next verify condition \hyperref[V2]{(V2)} and complete
	the verification of condition \hyperref[V3]{(V3)}.
	As noted above, if $\mathcal G_i\cap\mathcal G_j$ is an edge of order $2$, then $\beta_i(w_j')=\beta_j(w_i')=0$ by construction.
	Consider now an ordered pair $(p,q)$ with $p\neq q$ such that
	either $\mathcal G_p$ and $\mathcal G_q$ are not adjacent or
	their common edge has order at least $3$.
	Since $(\alpha,v)\in\widetilde{\mathcal D}(\mathcal G)$, we
	have
	\[
	\alpha_p(v_q)<0.
	\]
	Indeed, for an adjacent pair this follows from
	\hyperref[V2]{(V2)} and \hyperref[V4]{(V4)} for
	$(\alpha,v)$, whereas for a nonadjacent pair it follows from
	\hyperref[V2]{(V2)} and \hyperref[V5]{(V5)} for
	$(\alpha,v)$.
	
	The function
	\[
	(\beta,w)
	\mapsto
	\beta_p\bigl(w_q'(\beta,w)\bigr)
	\]
	is continuous and has the strictly negative value
	$\alpha_p(v_q)$ at $(\alpha,v)$.
	Since there are only finitely many such ordered pairs, we may
	shrink $U_\alpha$ and $W_v$ simultaneously so that
	\[
	\beta_p(w_q')<0
	\]
	for every such pair and every
	$(\beta,w)\in U_\alpha\times W_v$.
	Together with the vanishing relations for the edges of order
	$2$, this proves condition \hyperref[V2]{(V2)}.
	Moreover, for every pair of distinct indices $i,j$, the values $\beta_i(w_j')$ and $\beta_j(w_i')$ both vanish if $\mathcal G_i\cap\mathcal G_j$ is an edge of order $2$, and are both strictly negative otherwise.
	Hence condition \hyperref[V3]{(V3)} also holds.
	
	Condition \hyperref[V6]{(V6)} is immediate from
	\[
	\beta\in U_\alpha
	\subset
	\widetilde{\mathcal E}(\mathcal G)
	\]
	and Lemma~\ref{realization space membership}.
	
	It remains to verify condition \hyperref[V5]{(V5)}.
	For each $i$, let
	\[
	R_i' := \text{Id} - \beta_i \otimes w_i'.
	\]
	For every adjacent pair $\mathcal G_i,\mathcal G_j$, conditions \hyperref[V1]{(V1)}--\hyperref[V4]{(V4)} satisfy the local hypotheses of Vinberg's criterion stated in Theorem~\ref{Vinberg 1}.
	Therefore, the reflections
	\[
	R_i':=\text{Id}-\beta_i\otimes w_i',
	\qquad
	i=1,\ldots,f,
	\]
	generate a linear reflection group with fundamental cone determined by $\beta$.
	If $\mathcal G_i$ and $\mathcal G_j$ are not adjacent, \cite[Theorem~2(6)]{MR0302779} shows that $R_i'R_j'$ has infinite order.
	Proposition~17 of \cite{MR0302779} then rules out the finite-order alternative and gives
	\[
	\beta_i(w_j')\beta_j(w_i')\geq4.
	\]
	This proves condition \hyperref[V5]{(V5)}.
	
	We have now verified conditions \hyperref[V1]{(V1)} through
	\hyperref[V6]{(V6)}.
	Therefore,
	\[
	\bigl(\beta,w'(\beta,w)\bigr)
	\in
	\widetilde{\mathcal D}(\mathcal G)
	\]
	for every $(\beta,w)\in U_\alpha\times W_v$.
	This proves condition \textup{(i)}.
	The continuity and uniqueness of $w'$, together with
	conditions \textup{(ii)} and \textup{(iii)}, were established
	in the construction above.
	This completes the proof.
\end{proof}

For each $(\alpha, v) \in \widetilde{\mathcal{D}} (\mathcal{G})$, let $w': U_\alpha \times W_v \to V^f$ be the map given by Lemma \ref{standard chart lemma}.
Then we define a map $\Phi_{\alpha, v}: U_\alpha \times W_v \to \widetilde{\mathcal{D}} (\mathcal{G})$ by $\Phi_{\alpha, v} (\beta, w) := (\beta, w')$.
We refer to such maps as \textit{standard charts around the point} $(\alpha, v)$.

\subsubsection{Topological properties of the standard charts}\label{subsubsec:standard-chart-topology}

The maps $\Phi_{\alpha,v}$ constructed above are continuous and preserve the base coordinate.
To show that they provide local product charts, it remains to prove that each map is injective and has open image.
We do this by recursively recovering the original fiber coordinate $w$ from a nearby point $(\beta,w')$ in the image, thereby constructing a continuous local inverse.

\begin{lemma}\label{standard-chart-embedding}
	Let $(\alpha, v) \in \widetilde{\mathcal{D}} (\mathcal{G})$ and let $U_\alpha \subset \widetilde{\mathcal{E}} (\mathcal{G}), W_v \subset \rho^{-1} (\alpha)$ be open neighborhoods of $\alpha, v$ as in Lemma \ref{standard chart lemma}.
	Then the standard chart $\Phi_{\alpha, v}: U_\alpha \times W_v \to \widetilde{\mathcal{D}} (\mathcal{G})$ is a topological embedding onto an open subset of $\widetilde{\mathcal{D}} (\mathcal{G})$.
\end{lemma}

\begin{proof}
	The map $\Phi_{\alpha,v}$ is continuous by Lemma~\ref{standard chart lemma}.
	We construct a continuous local inverse of $\Phi_{\alpha,v}$ near each point of its image.
	
	For each $i\in\{1,\ldots,f\}$, set
	\[
	A_i := \{k_{i,1},\ldots,k_{i,a_i}\},\quad B_i := \{k_{i,a_i+1},\ldots,k_{i,a_i+b_i}\}, \quad \text{and} \quad C_i := \{k_{i,a_i+b_i+1},\ldots,k_{i,3}\}.
	\]
	Recall that every $j\in B_i$ satisfies $j<i$.
	
	Fix $(\beta,w)\in U_\alpha\times W_v$ and write
	\[
	(\beta,w') := \Phi_{\alpha,v}(\beta,w).
	\]
	For $(\gamma,y)\in (V^*)^f\times V^f$ sufficiently close to $(\beta,w')$, we define
	\[
	X(\gamma,y) = \bigl( X_1(\gamma,y),\ldots,X_f(\gamma,y) \bigr) \in V^f
	\]
	inductively as follows.
	For each $i$, let $X_i(\gamma,y)$ be the unique solution of the system
	\begin{equation}\label{system 2}
		\begin{aligned}
			\alpha_i(X_i(\gamma,y)) &=2,\\ 
			\alpha_j(X_i(\gamma,y)) &=0 &&\text{for }j\in A_i,\\
			\alpha_j(X_i(\gamma,y)) &= \frac{ 4\cos^2\left(\frac{\pi}{m_{i,j}}\right) }{ \alpha_i(X_j(\gamma,y)) } &&\text{for }j\in B_i,\\
			\alpha_j(X_i(\gamma,y)) &= \gamma_j(y_i) &&\text{for }j\in C_i.
		\end{aligned}
	\end{equation}
	Since $j<i$ for every $j\in B_i$, the vectors appearing on the right-hand side have already been constructed.
	Moreover, the four linear functionals $\alpha_i, \alpha_{k_{i,1}}, \alpha_{k_{i,2}}, \alpha_{k_{i,3}}$ form a basis of $V^*$.
	Hence each system \ref{system 2} has a unique solution.
	
	We claim that $X(\beta,w')=w$.
	Indeed, this follows inductively from the uniqueness of the solutions.
	For $j\in C_i$, condition~\textup{(iii)} of Lemma~\ref{standard chart lemma} gives $\beta_j(w_i')=\alpha_j(w_i)$, while the remaining equations follow from the fact that $(\alpha,w)\in\widetilde{\mathcal D}(\mathcal G)$.
	Thus $w_i$ satisfies the defining system for $X_i(\beta,w')$ for every $i$, and consequently $X_i(\beta,w')=w_i$.
	
	For $j\in B_i$, we have $\alpha_i(w_j)<0$.
	Therefore, after restricting to sufficiently small neighborhoods of $\beta$ and $w'$, all the denominators in the systems \ref{system 2} are nonzero.
	The recursive construction then shows that $X(\gamma,y)$ depends continuously on $(\gamma,y)$.
	
	Choose open neighborhoods $S_\beta\subset(V^*)^f$ and $T_{w'}\subset V^f$ of $\beta$ and $w'$, respectively, such that $
	S_\beta\cap\widetilde{\mathcal E}(\mathcal G)\subset U_\alpha$.
	Since $W_v$ is open in the subspace $\rho^{-1}(\alpha)\subset V^f$, choose an open set $\widehat W_v\subset V^f$ such that $W_v=\widehat W_v\cap\rho^{-1}(\alpha)$.
	Because $X(\beta,w')=w\in W_v$ and $X$ is continuous, we may shrink $S_\beta$ and $T_{w'}$ so that $X(\gamma,y)\in\widehat W_v$ for every $(\gamma,y)\in\widetilde{\mathcal D}(\mathcal G)\cap(S_\beta\times T_{w'})$.
	
	We now verify that $(\alpha,X(\gamma,y))\in\widetilde{\mathcal D}(\mathcal G)$.
	The defining systems for the $X_i$ give conditions \hyperref[V1]{(V1)}, \hyperref[V3]{(V3)}, and \hyperref[V4]{(V4)}.
	For every pair $i\neq j$ with $m_{i,j}=2$, they give $\alpha_i(X_j(\gamma,y))=0$.
	After shrinking the neighborhoods once more, continuity and the equality $X(\beta,w')=w$ give $\alpha_i(X_j(\gamma,y))<0$ for every remaining pair $i\neq j$.
	This proves condition \hyperref[V2]{(V2)}.
	Since $\alpha\in\widetilde{\mathcal E}(\mathcal G)$, Lemma~\ref{realization space membership} gives condition \hyperref[V6]{(V6)}.
	
	For each $i$, set
	\[
	R_i^X:=\text{Id}-\alpha_i\otimes X_i(\gamma,y).
	\]
	By Vinberg's criterion stated in Theorem~\ref{Vinberg 1}, the reflections $R_1^X,\ldots,R_f^X$ generate a linear reflection group with fundamental cone determined by $\alpha$.
	For nonadjacent facets $\mathcal G_i$ and $\mathcal G_j$, \cite[Theorem~2(6) and Proposition~17]{MR0302779} gives
	\[
	\alpha_i(X_j(\gamma,y))\alpha_j(X_i(\gamma,y))\geq4.
	\]
	This proves condition \hyperref[V5]{(V5)}.
	We have shown that $(\alpha,X(\gamma,y))\in\widetilde{\mathcal D}(\mathcal G)$, and therefore
	\[
	X(\gamma,y)\in\rho^{-1}(\alpha)\cap\widehat W_v=W_v.
	\]
	
	Set
	\[
	\mathcal N
	:=
	\widetilde{\mathcal D}(\mathcal G)
	\cap
	(S_\beta\times T_{w'})
	\]
	and define $F:\mathcal N\to U_\alpha\times W_v$, $F(\gamma,y)
	:=
	\bigl(\gamma,X(\gamma,y)\bigr)$.
	The map $F$ is continuous.
	We show that it is an inverse of $\Phi_{\alpha,v}$ on $\mathcal N$.
	
	First, let $(\gamma,y)\in\mathcal N$.
	For each $i$, the vector $y_i$ satisfies
	\begin{align*}
		\gamma_i(y_i)
		&=2,\\
		\gamma_j(y_i)
		&=0
		&&\text{for }j\in A_i,\\
		\gamma_j(y_i)
		&=
		\frac{
			4\cos^2\left(\frac{\pi}{m_{i,j}}\right)
		}{
			\gamma_i(y_j)
		}
		&&\text{for }j\in B_i,
	\end{align*}
	since $(\gamma,y)\in\widetilde{\mathcal D}(\mathcal G)$.
	For $j\in C_i$, the definition of $X_i(\gamma,y)$ gives
	\[
	\gamma_j(y_i)=\alpha_j(X_i(\gamma,y)).
	\]
	It follows that $y_i$ is the unique solution of the system defining
	the $i$-th component of $\Phi_{\alpha,v} \bigl(\gamma,X(\gamma,y)\bigr)$.
	Proceeding inductively over $i$, we obtain $\Phi_{\alpha,v}\circ F = \text{id}_{\mathcal N}$.
	
	Conversely, let $(\delta,z)\in U_\alpha\times W_v$ satisfy $ \Phi_{\alpha,v}(\delta,z)\in\mathcal N$, 
	and write $\Phi_{\alpha,v}(\delta,z) = (\delta,z')$.
	Since $(\alpha,z)\in\widetilde{\mathcal D}(\mathcal G)$, the vectors $z_i$ satisfy the first three groups of equations defining $X_i(\delta,z')$.
	For $j\in C_i$, condition~\textup{(iii)} of Lemma~\ref{standard chart lemma} gives
	\[
	\delta_j(z_i')
	=
	\alpha_j(z_i).
	\]
	Thus $z_i$ satisfies the entire system defining
	$X_i(\delta,z')$.
	By induction and uniqueness, we obtain
	\[
	X(\delta,z')=z.
	\]
	Hence $F\circ\Phi_{\alpha,v} = \text{id}_{\Phi_{\alpha,v}^{-1}(\mathcal N)}$.
	
	It follows that the restriction $\Phi_{\alpha,v}: \Phi_{\alpha,v}^{-1}(\mathcal N) \to \mathcal N $ is a homeomorphism.
	In particular, every point in the image of $\Phi_{\alpha,v}$ has an
	open neighborhood in
	$\widetilde{\mathcal D}(\mathcal G)$ contained in that image.
	Therefore, the image of $\Phi_{\alpha,v}$ is open in $\widetilde{\mathcal{D}} (\mathcal{G})$.
	
	Finally, if two points of $U_\alpha\times W_v$ have the same image,
	choose the neighborhood $\mathcal N$ constructed above around their
	common image.
	Both points belong to $\Phi_{\alpha,v}^{-1}(\mathcal N)$, on which
	$\Phi_{\alpha,v}$ is injective.
	Thus the two points are equal.
	Therefore, $\Phi_{\alpha,v}$ is injective and is a homeomorphism onto
	an open subset of $\widetilde{\mathcal D}(\mathcal G)$.
\end{proof}

\subsubsection{Smooth compatibility and dimension}\label{subsubsec:standard-chart-smoothness}

The images of the standard charts $\Phi_{\alpha,v}$ cover $\widetilde{\mathcal{D}}(\mathcal{G})$.
To show that they form a smooth atlas, it remains to prove that the transition maps are smooth.
Because every standard chart preserves the base coordinate, a transition map fixes its first coordinate, while its second coordinate is recovered recursively from systems of four linear equations.
After establishing smooth compatibility, we compute the dimension of $\widetilde{\mathcal{D}}(\mathcal{G})$.

\begin{lemma}\label{transition}
	Let $\Phi^{(l)}: U_\alpha^{(l)} \times W_v^{(l)} \to \widetilde{\mathcal{D}} (\mathcal{G})$, $l = 1,2$, be two standard charts around the points $(\alpha^{(l)}, v^{(l)}) \in \widetilde{\mathcal{D}} (\mathcal{G})$, for $l = 1,2$.
	Suppose that the intersection 
	\[
	C := \Phi^{(1)} \left(U_\alpha^{(1)} \times W_v^{(1)}\right) \cap \Phi^{(2)} \left(U_\alpha^{(2)} \times W_v^{(2)}\right)
	\]
	is nonempty.
	Then the map
	\[
	(\Phi^{(2)})^{-1} \circ \Phi^{(1)}: (\Phi^{(1)})^{-1} (C) \to (\Phi^{(2)})^{-1} (C)
	\]
	is smooth.
	
	In particular, the standard charts of $\widetilde{\mathcal{D}} (\mathcal{G})$ define a smooth structure on $\widetilde{\mathcal{D}} (\mathcal{G})$.
	If $\widetilde{\mathcal D}(\mathcal G)$ is nonempty, then
	\[
	\dim\widetilde{\mathcal D}(\mathcal G)=4f-e_2+6.
	\]
\end{lemma}

\begin{proof}
	Let $(\beta^{(1)}, w^{(1)}) \in (\Phi^{(1)})^{-1} (C)$ and let $(\beta^{(2)}, w^{(2)}) := (\Phi^{(2)})^{-1} \circ \Phi^{(1)} (\beta^{(1)}, w^{(1)})$.
	Define $(\beta, w') := \Phi^{(1)} (\beta^{(1)}, w^{(1)}) = \Phi^{(2)} (\beta^{(2)}, w^{(2)})$.
	
	Each standard chart preserves the first coordinate, so the equality defining $(\beta,w')$ implies $\beta^{(1)}=\beta=\beta^{(2)}$.
	Furthermore, the following relations hold due to the fact that $(\alpha^{(2)}, w^{(2)}) \in \widetilde{\mathcal{D}} (\mathcal{G})$ and from the construction of the standard charts:
	\begin{align*}
		\alpha_i^{(2)} (w_i^{(2)}) &= 2, \\
		\alpha_j^{(2)} (w_i^{(2)}) &= 0 \quad &\text{for} \quad &j \in \{k_{i, 1}, \ldots, k_{i, a_i}\}, \\
		\alpha_j^{(2)} (w_i^{(2)}) &= \frac{4 \cos^2 \left( \frac{\pi}{m_{i, j}} \right) }{\alpha_i^{(2)} (w_j^{(2)})} \quad &\text{for} \quad &j \in \{k_{i, a_i + 1}, \ldots, k_{i, a_i + b_i}\}, \\
		\alpha^{(2)}_j (w^{(2)}_i) &= \beta^{(2)}_j (w'_i) = \beta^{(1)}_j (w'_i) = \alpha^{(1)}_j (w^{(1)}_i) \quad &\text{for} \quad &j \in \{k_{i,a_i+b_i+1}, \ldots, k_{i,3}\}.\notag
	\end{align*}
	The last equalities follow from the fact that $\beta^{(1)} = \beta^{(2)}$ and the condition (iii) of Lemma \ref{standard chart lemma}.
	
	Since $\beta^{(1)}=\beta^{(2)}=\beta$, the transition map leaves the first coordinate unchanged.
	Moreover, the displayed systems recover $w^{(2)}$ recursively from $w^{(1)}$.
	Indeed, for the indices $j\in\{k_{i,a_i+b_i+1},\ldots,k_{i,3}\}$, the last line gives $ \alpha_j^{(2)}(w_i^{(2)}) = \alpha_j^{(1)}(w_i^{(1)})$.
	The other values of $\alpha_j^{(2)}(w_i^{(2)})$ are determined recursively by conditions \hyperref[V1]{(V1)}, \hyperref[V3]{(V3)}, and \hyperref[V4]{(V4)}.
	Thus all the coefficient functionals in the systems determining $w^{(2)}$ are the fixed functionals $\alpha_i^{(2)},\alpha_{k_{i,1}}^{(2)},\alpha_{k_{i,2}}^{(2)},\alpha_{k_{i,3}}^{(2)}$, while their right-hand sides depend only on $w^{(1)}$ and not on $\beta$.
	
	Fix a point in the domain of the transition map.
	At this point, all the denominators occurring in the recursive construction are nonzero, and hence they remain nonzero on a sufficiently small neighborhood.
	By Cramer's rule, the recursive construction therefore defines a smooth map $F:\mathcal U\longrightarrow V^f$
	on an open neighborhood $\mathcal U\subset V^f$ of $w^{(1)}$ such that $w^{(2)}=F(w^{(1)})$.
	Since the fibers $\rho^{-1}(\alpha^{(1)})$ and $\rho^{-1}(\alpha^{(2)})$ are embedded submanifolds of $V^f$ by Proposition~\ref{prop:embedded-fibers}, the restriction of $F$ to the relevant open subset of $\rho^{-1}(\alpha^{(1)})$ is a smooth map into $\rho^{-1}(\alpha^{(2)})$.
	It follows that the transition map is locally of the form $(\beta,w)\mapsto(\beta,F(w))$, and is therefore smooth.
	
	Since the underlying combinatorial polytope of $\mathcal G$ is not a cone over a polygon, Proposition~\ref{prop:stabilizer-dimension} gives $k(\mathcal G)=0$.
	Thus Theorem~\ref{restricted} gives $\dim\mathcal C_P(\mathcal G)=3f-e-e_2$.
	
	Combining this with Lemma~\ref{principal bundle}, we obtain
	\begin{align*}
		\dim \widetilde{\mathcal{D}}(\mathcal G)
		&= \dim \widetilde{\mathcal{E}}(\mathcal G)+\dim\rho^{-1}(\alpha)\\
		&= \dim \widetilde{\mathcal{E}}(\mathcal G)+\dim\mathcal C_P(\mathcal G)\\
		&= (e+f+6)+(3f-e-e_2)\\
		&= 4f-e_2+6,
	\end{align*}
	where $\alpha\in\widetilde{\mathcal E}(\mathcal G)$ is any element such that $\rho^{-1}(\alpha)$ is nonempty, and $P\subset\mathbb S^3$ is the convex $3$-polytope determined by $\alpha$.
\end{proof}

\subsubsection{Openness of the forgetful maps}\label{subsubsec:forgetful-openness}

The local product form of the standard charts has one further consequence.
In each standard chart, the map $\rho$ is represented by the projection onto the first factor.
We first use this observation to prove that $\rho$ is open and then descend the conclusion through the quotient diagram in Figure~\ref{commuting diagram} to $\overline{\rho}$ and $\kappa$.

\begin{corollary}\label{open-forgetful-map}
	Suppose that the underlying combinatorial polytope of
	$\mathcal G$ is not a cone over a polygon.
	Then the maps $\rho$, $\overline{\rho}$, and $\kappa$ in
	Figure~\ref{commuting diagram} are open.
	In particular, $\kappa(\mathcal C(\mathcal G))$ is an open
	submanifold of $\mathcal{RS}(\mathcal G)$.
\end{corollary}

\begin{proof}
	We first prove that $\rho$ is open.
	For every standard chart $\Phi_{\alpha,v}$, its definition gives $\rho\circ\Phi_{\alpha,v}=\text{pr}_1$, where $\text{pr}_1:U_\alpha\times W_v\to U_\alpha$ is the projection onto the first factor.
	Thus $\rho$ is locally represented by a projection and, in particular, is a smooth submersion.
	
	Let $\mathcal O\subset\widetilde{\mathcal D}(\mathcal G)$ be open.
	For each standard chart $\Phi_{\alpha,v}$, set
	\[
	\mathcal O_{\alpha,v} := \Phi_{\alpha,v}^{-1} \bigl( \mathcal O \cap \Phi_{\alpha,v}(U_\alpha\times W_v) \bigr).
	\]
	This is an open subset of $U_\alpha\times W_v$.
	Since $\text{pr}_1$ is open, $\text{pr}_1(\mathcal O_{\alpha,v})$ is open in $U_\alpha$, and hence in $\widetilde{\mathcal E}(\mathcal G)$.
	By the identity $\rho\circ\Phi_{\alpha,v}=\text{pr}_1$, this set is precisely the image under $\rho$ of the part of $\mathcal O$ contained in the image of $\Phi_{\alpha,v}$.
	Since the images of the standard charts cover $\widetilde{\mathcal D}(\mathcal G)$, the set $\rho(\mathcal O)$ is open.
	Therefore, $\rho$ is an open map.
	
	We next prove that $\overline{\rho}$ is open.
	Let $\mathcal U\subset\mathcal D(\mathcal G)$ be open.
	Then $\eta^{-1}(\mathcal U)$ is open in $\widetilde{\mathcal D}(\mathcal G)$, and hence $\rho(\eta^{-1}(\mathcal U))$ is open in $\widetilde{\mathcal E}(\mathcal G)$.
	The principal-bundle projection $\pi$ is open.
	By the commutativity of Figure~\ref{commuting diagram}, we have
	\[
	\overline{\rho}(\mathcal U) = \pi\bigl(\rho(\eta^{-1}(\mathcal U))\bigr).
	\]
	Therefore, $\overline{\rho}(\mathcal U)$ is open in $\mathcal E(\mathcal G)$.
	
	Finally, $\Phi:\mathcal C(\mathcal G)\to\mathcal D(\mathcal G)$ and $\iota:\mathcal{RS}(\mathcal G)\to\mathcal E(\mathcal G)$ are homeomorphisms onto their images.
	The commutativity of Figure~\ref{commuting diagram} gives $\iota\circ\kappa=\overline{\rho}\circ\Phi$.
	Hence $\kappa=\iota^{-1}\circ\overline{\rho}\circ\Phi$ is open.
\end{proof}

\subsection{Proof of the main theorem}\label{Section 6.4}

In this section, we finalize the proof of Theorem \ref{main theorem}.

The formula
\begin{align*}
	&(A,c_1,\ldots,c_f) \cdot (\alpha_1,\ldots,\alpha_f,v_1,\ldots,v_f) \\
	&\qquad := \bigl(c_1^{-1}\alpha_1\circ A^{-1}, \ldots, c_f^{-1}\alpha_f\circ A^{-1}, c_1Av_1, \ldots, c_fAv_f \bigr)
\end{align*}
restricts to a continuous action of $\mathbb G = \text{SL}_\pm(V)\times\mathbb R_+^f$ on $\widetilde{\mathcal D}(\mathcal G)$.
We first verify that this action is smooth with respect to the smooth structure constructed in Section~\ref{Section 6.3}.

\begin{lemma}\label{proper and free}
	Suppose that the underlying combinatorial polytope of
	$\mathcal G$ is not a cone over a polygon.
	Then the action of $\mathbb G$ on
	$\widetilde{\mathcal D}(\mathcal G)$ is smooth, proper, and free.
\end{lemma}

\begin{proof}
	We first prove smoothness.
	For each $i\in\{1,\ldots,f\}$, set
	\[ A_i := \{k_{i,1},\ldots,k_{i,a_i}\}, \quad B_i := \{k_{i,a_i+1},\ldots,k_{i,a_i+b_i}\}, \]
	and
	\[
	C_i := \{k_{i,a_i+b_i+1},\ldots,k_{i,3}\}.
	\]
	Recall that every $j\in B_i$ satisfies $j<i$.
	
	Let
	\[
	\mu: \mathbb G\times\widetilde{\mathcal D}(\mathcal G) \to \widetilde{\mathcal D}(\mathcal G), \quad \mu(g,z):=g\cdot z,
	\]
	denote the action map.
	Fix $\bigl(g_0,(\alpha,v)\bigr) \in \mathbb G\times\widetilde{\mathcal D}(\mathcal G)$, and put $(\bar\alpha,\bar v) := g_0\cdot(\alpha,v)$.
	
	Choose standard charts $\Phi_{\alpha,v}: U_\alpha\times W_v \to \widetilde{\mathcal D}(\mathcal G)$ and $\Phi_{\bar\alpha,\bar v} : U_{\bar\alpha}\times W_{\bar v} \to \widetilde{\mathcal D}(\mathcal G)$.
	By the construction of the local inverse in the proof of Lemma~\ref{standard-chart-embedding}, there is an open neighborhood $\mathcal N$ of $(\bar\alpha,\bar v)$ in $\widetilde{\mathcal D}(\mathcal G)$ such that $\Phi_{\bar\alpha,\bar v}^{-1}(\gamma,z) = \bigl(\gamma,X(\gamma,z)\bigr)$ for $(\gamma,z)\in\mathcal N$, where
	\[
	X(\gamma,z)
	=
	\bigl(X_1(\gamma,z),\ldots,X_f(\gamma,z)\bigr)
	\]
	is defined inductively by
	\begin{align}
		\bar\alpha_i\bigl(X_i(\gamma,z)\bigr)
		&=2,
		\label{action-system-1}
		\\
		\bar\alpha_j\bigl(X_i(\gamma,z)\bigr)
		&=0
		&&
		\text{for }j\in A_i,
		\label{action-system-2}
		\\
		\bar\alpha_j\bigl(X_i(\gamma,z)\bigr)
		&=
		\frac{
			4\cos^2\left(\frac{\pi}{m_{i,j}}\right)
		}{
			\bar\alpha_i\bigl(X_j(\gamma,z)\bigr)
		}
		&&
		\text{for }j\in B_i,
		\label{action-system-3}
		\\
		\bar\alpha_j\bigl(X_i(\gamma,z)\bigr)
		&=
		\gamma_j(z_i)
		&&
		\text{for }j\in C_i.
		\label{action-system-4}
	\end{align}
	
	Since the action is continuous with respect to the original topology, after shrinking the two standard charts and choosing a sufficiently small open neighborhood $\mathcal O$ of $g_0$ in $\mathbb G$, we may assume that $\mu\left(\mathcal O \times \Phi_{\alpha,v}(U_\alpha\times W_v) \right) \subset \mathcal N$.
	
	Let $g=(T,c_1,\ldots,c_f)\in\mathcal O$, $(\beta,w)\in U_\alpha\times W_v$, and write $\Phi_{\alpha,v}(\beta,w)=(\beta,y)$.
	Thus $y=\bigl(y_1,\ldots,y_f\bigr)$ is the tuple constructed in Lemma~\ref{standard chart lemma}.
	Write $g\cdot(\beta,y) = (\gamma,z)$.
	Then we have $\gamma_j = c_j^{-1}\beta_j\circ T^{-1}$, $z_i = c_iTy_i$.
	
	For every $j\in C_i$, we have
	\begin{align*}
		\gamma_j(z_i) = \bigl(c_j^{-1}\beta_j\circ T^{-1}\bigr) (c_iTy_i) = \frac{c_i}{c_j}\beta_j(y_i) = \frac{c_i}{c_j}\alpha_j(w_i).
	\end{align*}
	The last equality follows from condition~\textup{(iii)} of
	Lemma~\ref{standard chart lemma}.
	
	It follows that the second component of $\Phi_{\bar\alpha,\bar v}^{-1} \bigl(g\cdot\Phi_{\alpha,v}(\beta,w) \bigr)$ is the tuple
	\[
	\Xi(g,w) = \bigl( \Xi_1(g,w),\ldots,\Xi_f(g,w) \bigr)
	\]
	defined inductively by
	\begin{align}
		\bar\alpha_i\bigl(\Xi_i(g,w)\bigr)
		&=2,
		\label{transformed-action-system-1}
		\\
		\bar\alpha_j\bigl(\Xi_i(g,w)\bigr)
		&=0
		&&
		\text{for }j\in A_i,
		\label{transformed-action-system-2}
		\\
		\bar\alpha_j\bigl(\Xi_i(g,w)\bigr)
		&=
		\frac{
			4\cos^2\left(\frac{\pi}{m_{i,j}}\right)
		}{
			\bar\alpha_i\bigl(\Xi_j(g,w)\bigr)
		}
		&&
		\text{for }j\in B_i,
		\label{transformed-action-system-3}
		\\
		\bar\alpha_j\bigl(\Xi_i(g,w)\bigr)
		&=
		\frac{c_i}{c_j}\alpha_j(w_i)
		&&
		\text{for }j\in C_i.
		\label{transformed-action-system-4}
	\end{align}
	
	For every $i$, the four coefficient functionals $\bar\alpha_i$, $\bar\alpha_{k_{i,1}}$, $\bar\alpha_{k_{i,2}}$, $\bar\alpha_{k_{i,3}}$ form a basis of $V^*$ by Remark~\ref{rem:selected-basis}.
	Moreover, if $j\in B_i$, then $\bar\alpha_i(\bar v_j)<0$.
	Hence, after shrinking $\mathcal O$ and $W_v$ if necessary, all the denominators in \eqref{transformed-action-system-3} remain nonzero.
	Since every $j\in B_i$ satisfies $j<i$, induction on $i$ and Cramer's rule show that $\Xi: \mathcal O\times W_v \to V^f$ is smooth.
	
	After a further shrinking, its image lies in $W_{\bar v} \subset \rho^{-1}(\bar\alpha)$.
	Since $\rho^{-1}(\bar\alpha)$ is an embedded submanifold of $V^f$ by Proposition~\ref{prop:embedded-fibers}, the map $\Xi: \mathcal O\times W_v \to W_{\bar v}$ is smooth as a map between smooth manifolds.
	
	We have therefore obtained the coordinate expression
	\[
	\Phi_{\bar\alpha,\bar v}^{-1}
	\circ
	\mu
	\circ
	\bigl(
	\text{id}_{\mathbb G}
	\times
	\Phi_{\alpha,v}
	\bigr)
	(g,\beta,w)
	=
	\bigl(
	g\cdot\beta,
	\Xi(g,w)
	\bigr).
	\]
	By Lemma~\ref{principal bundle}, $\widetilde{\mathcal E}(\mathcal G)$ has been endowed with a smooth structure for which $
	\widetilde{\mathcal E}(\mathcal G) \to \mathcal E(\mathcal G)$ is a smooth principal $\mathbb G$-bundle.
	In particular, the action map $\mathbb G\times\widetilde{\mathcal E}(\mathcal G) \to \widetilde{\mathcal E}(\mathcal G)$ is smooth.
	Thus both components of the preceding coordinate expression are smooth.
	Since $\bigl(g_0,(\alpha,v)\bigr)$ was arbitrary, the action of $\mathbb G$ on $\widetilde{\mathcal D}(\mathcal G)$ is smooth.
	
	We next prove freeness and properness.
	Since the action of $\mathbb G$ on $\widetilde{\mathcal E}(\mathcal G)$ is free, its action on $\widetilde{\mathcal D}(\mathcal G)$ is also free.
	
	By Lemma~\ref{principal bundle}, the action of $\mathbb G$ on $\widetilde{\mathcal E}(\mathcal G)$ is proper.
	Let $K \subset \widetilde{\mathcal D}(\mathcal G) \times \widetilde{\mathcal D}(\mathcal G)$ be compact, and suppose that
	\[
	\bigl( g_n\cdot(\alpha^{(n)},v^{(n)}), (\alpha^{(n)},v^{(n)}) \bigr) \in K
	\]
	for sequences $g_n\in\mathbb G$, $(\alpha^{(n)},v^{(n)}) \in \widetilde{\mathcal D}(\mathcal G)$.
	Set
	\[
	K_{\mathcal E}
	:=
	(\rho\times\rho)(K)
	\subset
	\widetilde{\mathcal E}(\mathcal G)
	\times
	\widetilde{\mathcal E}(\mathcal G).
	\]
	The set $K_{\mathcal E}$ is compact, and $\bigl(g_n\cdot\alpha^{(n)}, \alpha^{(n)} \bigr) \in K_{\mathcal E}$.
	The properness of the action on $\widetilde{\mathcal E}(\mathcal G)$ therefore implies that $\{g_n\}$ has a convergent subsequence.
	After passing to this subsequence, compactness of $K$ gives a further
	subsequence for which
	\[
	\bigl(g_n\cdot(\alpha^{(n)},v^{(n)}),
	(\alpha^{(n)},v^{(n)})\bigr)
	\]
	converges in
	$\widetilde{\mathcal D}(\mathcal G)\times
	\widetilde{\mathcal D}(\mathcal G)$.
	In particular, $(\alpha^{(n)},v^{(n)})$ converges.
	Hence the action of $\mathbb G$ on $\widetilde{\mathcal D}(\mathcal G)$ is proper.
\end{proof}

We conclude the proof of our main theorem.

\begin{proof}[Proof of Theorem~\ref{main theorem}]
	By the embedding $\Phi$ in Lemma~\ref{embedding lemma} and the construction in Section~\ref{Section 6.2}, the space $\mathcal C(\mathcal G)$ is homeomorphic to $\mathcal D(\mathcal G)=\widetilde{\mathcal D}(\mathcal G)/\mathbb G$.
	Since the underlying combinatorial polytope of $\mathcal G$ is not a cone over a polygon, Lemma~\ref{proper and free} shows that the action of $\mathbb G$ on $\widetilde{\mathcal D}(\mathcal G)$ is smooth, proper, and free.
	Therefore, the quotient admits a smooth structure for which the quotient map is a smooth submersion.
	
	Assume now that $\mathcal C(\mathcal G)$ is nonempty.
	By Lemma~\ref{transition}, we have $\dim\widetilde{\mathcal D}(\mathcal G)=4f-e_2+6$.
	Since $\dim\mathbb G=15+f$, we obtain
	\[
	\begin{aligned}
		\dim\mathcal C(\mathcal G)
		&=
		\dim\widetilde{\mathcal D}(\mathcal G)-\dim\mathbb G\\
		&=
		(4f-e_2+6)-(15+f)\\
		&=
		3f-e_2-9.
	\end{aligned}
	\]
\end{proof}

\begin{remark}
	In dimension four, Choi--Lee--Marquis constructed a perfect labeled polytope whose deformation space is homeomorphic to a circle; see \cite{MR4093439} for the definition of \textit{perfect} labeled polytope and \cite[Theorem~10.7]{MR4093439} for the example.
	We do not know whether there exists an orderable labeled combinatorial $3$-polytope of normal type for which a connected
	component of $\mathcal C(\mathcal G)$ is not homeomorphic to an open cell.
	By \cite[Remark~3.20]{MR2660566}, every simple orderable combinatorial $3$-polytope is a truncation polytope.
	Therefore, any such example would necessarily have to be non-simple.
\end{remark}

\section{Examples}\label{section 7}

The deformation space $\mathcal{C}(\mathcal{G})$ is naturally identified with a subspace of $((V^*)^f \times V^f) / \mathbb{G}$.
In the preceding section, we introduced the map $\kappa: \mathcal{C}(\mathcal{G}) \to \mathcal{RS} (\mathcal{G})$, which associates each Coxeter 3-polytope $[(\alpha, v)] = [(\alpha_1, \dots, \alpha_f, v_1, \dots, v_f)]$ with the projective equivalence class $[\alpha] = [(\alpha_1, \dots, \alpha_f)]$ of polytopes.
Each nonempty fiber of $\kappa$ corresponds to a restricted deformation space, and it is a smooth manifold in the settings of Theorem \ref{main theorem}.
By Corollary~\ref{open-forgetful-map}, the map $\kappa$ is open onto the open submanifold $\kappa(\mathcal C(\mathcal G))\subset\mathcal{RS}(\mathcal G)$.

In this section, we study two labeled combinatorial $3$-polytopes by explicitly computing their realization spaces and the nonempty fibers of the forgetful map $\kappa$.
In Section~\ref{Example 7.1}, we obtain an example for which $\mathcal{C}(\mathcal{G})$ is homeomorphic to an open $2$-cell; the map $\kappa$ is not surjective, but it is a fiber bundle over its image.
In Section~\ref{counterexample}, we consider an example for which $\mathcal{C}(\mathcal{G})$ is homeomorphic to $\mathbb{R}^4 \sqcup \mathbb{R}^4$ and show that $\kappa$ is neither surjective nor a fiber bundle over its image.
Together, these examples illustrate both the computation of the global deformation space and the possible variation of the restricted deformation spaces with the underlying realization.

\subsection{A deformation space homeomorphic to an open $2$-cell}\label{Example 7.1}

Let $\mathcal{G}$ be the labeled combinatorial $3$-polytope shown in Figure~\ref{example 1}.
This example serves as a model computation in which both the restricted deformation spaces and the global deformation space are determined explicitly.
In Section~\ref{example-1-hypotheses}, we verify the hypotheses of Theorem~\ref{main theorem} and record the dimensions of $\mathcal{C}(\mathcal{G})$, $\mathcal{E}(\mathcal{G})$, and the nonempty fibers of $\kappa$.
In Section~\ref{example-1-realization}, we normalize the supporting functionals and identify $\mathcal{E}(\mathcal{G})$ with $(1,\infty)$.
In Section~\ref{example-1-fibers}, we compute the fibers over these normalized realizations, show that they are nonempty exactly for $1<d<\frac{4}{3}$, and conclude that $\mathcal{C}(\mathcal{G})$ is homeomorphic to an open $2$-cell.

\begin{figure}[t]
	\centering
	\begin{tikzpicture}
		\draw (-0.67, -0.38) -- (0.67, -0.38);
		\draw (-0.67, -0.38) -- (0, 0.77);
		\draw (0.67, -0.38) -- (0, 0.77);
		
		\draw (-2, -1.15) -- (2, -1.15);
		\draw (-2, -1.15) -- (0, 2.31);
		\draw (2, -1.15) -- (0, 2.31);
		
		\draw (0, 0.77) -- (0, 2.31);
		\draw (-0.67, -0.38) -- (-2, -1.15);
		\draw (0.67, -0.38) -- (2, -1.15);
		
		\draw (0, 2.31) -- (0.67, -0.38);
		
		\node at (0, 0) {\footnotesize $\mathcal{G}_1$};
		\node at (0.22, 0.9) {\footnotesize $\mathcal{G}_2$};
		\node at (-0.67, 0.47) {\footnotesize $\mathcal{G}_3$};
		\node at (0, -0.77) {\footnotesize $\mathcal{G}_4$};
		\node at (2, 0.9) {\footnotesize $\mathcal{G}_5$};
		\node at (0.79, 0.44) {\footnotesize $\mathcal{G}_6$};
		
		\node at (-0.43, 0.29) {\tiny 2};
		\node at (0.23, 0.09) {\tiny 2};
		\node at (0, -0.48) {\tiny 2};
		\node at (-0.1, 1.54) {\tiny 2};
		\node at (0.43, 1.06) {\tiny 2};
		\node at (1.23, -0.87) {\tiny 2};
		\node at (1.1, 0.68) {\tiny 2};
		
		\node at (-1.23, -0.87) {\tiny 3};
		\node at (0, -1.25) {\tiny 3};
		\node at (-1.1, 0.68) {\tiny 3};
	\end{tikzpicture}
	\caption{A labeled combinatorial $3$-polytope with a $2$-dimensional deformation space; see Section~\ref{Example 7.1}.}
	\label{example 1}
\end{figure}

\subsubsection{Verification of the hypotheses and dimension computations}\label{example-1-hypotheses}

	The displayed indexing
	$\mathcal{G}_1,\ldots,\mathcal{G}_6$
	satisfies the orderability condition.
	Indeed, for $i=1,\ldots,6$, the numbers of edges of
	$\mathcal{G}_i$ which either have order $2$ or are of the form
	$\mathcal{G}_i\cap\mathcal{G}_j$ for some $j<i$ are 3,3,2,3,3,3, respectively.
	Thus $\mathcal{G}$ is orderable.
	We next verify that $\mathcal{G}$ is of normal type.
	Let $\Gamma_{\mathcal{G}} = \langle \gamma_1,\ldots,\gamma_6\rangle$ be its associated Coxeter group.
	Its Coxeter graph is connected: the vertices $\gamma_3,\gamma_4,\gamma_5$ form a triangle whose edges have label $3$, while
	$\gamma_1,\gamma_2,\gamma_6$ are joined to this triangle by the edges corresponding to $m_{1,5}=m_{2,4}=m_{3,6}=\infty$, respectively.
	Hence $\Gamma_{\mathcal{G}}$ is irreducible.
	Moreover, consider the $(2,3,\infty)$-triangle Coxeter group
	\[
	\Delta := \left\langle a,b,c \ \middle|\ a^2=b^2=c^2=(ab)^2=(bc)^3=1 \right\rangle.
	\]
	The assignment
	\[
	\gamma_1,\gamma_6\mapsto 1, \quad \gamma_2\mapsto a, \quad \gamma_3\mapsto b, \quad \gamma_4,\gamma_5\mapsto c
	\]
	respects all the defining relations of $\Gamma_{\mathcal{G}}$ and defines a surjective homomorphism $\Gamma_{\mathcal{G}}\to\Delta$.
	Let $\Delta^{+} := \ker\!\left( \Delta \to \mathbb{Z}/2\mathbb{Z} \right)$, where each of $a,b,c$ is mapped to the nontrivial element.
	Then $\Delta^{+}$ has index $2$ in $\Delta$.
	It is generated by $x:=ab$, $y:=bc$, and has the presentation
	\[
	\Delta^{+}
	\cong
	\left\langle
	x,y
	\ \middle|\
	x^2=y^3=1
	\right\rangle
	\cong
	\text{PSL}_2(\mathbb{Z}).
	\]
	By \cite[Lemma~1]{MR140583}, the commutator subgroup $[\Delta^{+},\Delta^{+}]$ is a free group of rank $2$ and has
	index $6$ in $\Delta^{+}$.
	Hence $\Delta$ is virtually nonabelian free and, in particular, is large.
	Since $\Gamma_{\mathcal{G}}$ surjects onto $\Delta$, the inverse
	image of $[\Delta^{+},\Delta^{+}]$ is a finite-index subgroup of
	$\Gamma_{\mathcal{G}}$ admitting a surjective homomorphism onto
	a free group of rank $2$.
	Thus $\Gamma_{\mathcal{G}}$ is large.
	Therefore, $\mathcal{G}$ is of normal type by Proposition~\ref{prop:normal-type-criterion}.
	In the notation of Theorem~\ref{main theorem}, we have $f=6$, $e=10$, $e_2=7$.
	The underlying combinatorial polytope is not a cone over a polygon, and hence $k(\mathcal G)=0$ by Proposition~\ref{prop:stabilizer-dimension}.
	It follows from Theorem~\ref{main theorem} that $\mathcal{C}(\mathcal{G})$ is a smooth manifold of dimension $3f-e_2-9=2$.
	Moreover, Theorem~\ref{realization spaces} shows that $\mathcal{RS}(\mathcal{G})\cong\mathcal{E}(\mathcal{G})$ is a smooth manifold of dimension $e-9=1$.
	Finally, the nonempty fibers of $\kappa:\mathcal C(\mathcal G)\to\mathcal{RS}(\mathcal G)$ are precisely the nonempty restricted deformation spaces.
	By Theorem~\ref{restricted}, each such fiber is a smooth manifold of dimension $3f-e-e_2-k(\mathcal G)=1$.
	
	\subsubsection{The realization space}\label{example-1-realization}
	
	We describe $\mathcal{E}(\mathcal{G})$ by choosing a normalized representative in each $\mathbb{G}$-orbit.
	We first normalize the five supporting functionals defining a triangular prism and then determine the sixth functional from the incidence relations of $\mathcal{G}$.
	
	Let $\alpha = (\alpha_1, \ldots, \alpha_6) \in \widetilde{\mathcal{E}} (\mathcal{G})$ be an arbitrary element.
	We note that the first five linear inequalities $\alpha_1 \geq 0, \ldots, \alpha_5 \geq 0$ determine a convex 3-polytope, which is a triangular prism.
	The projective realization space of a triangular prism consists of a single point.
	Indeed, the polar dual of a triangular prism is a triangular bipyramid.
	No four vertices of a triangular bipyramid lie in a projective plane: the only nontrivial case consists of the two apices and two vertices of the equatorial triangle, and a plane containing them would force the segment joining the apices to meet the equatorial plane on the edge joining those two vertices, whereas it meets the interior of the equatorial triangle.
	Thus, the five vertices form a projective frame in $\mathbb S^3$.
	Any two ordered projective frames are related by an element of $\text{SL}_{\pm}(4,\mathbb R)$.
	Passing to polar duals, any two triangular prisms are therefore projectively equivalent as realizations with indexed facets.
	Hence, we may normalize the first five supporting functionals so that the triangular prism is defined by the inequalities $\beta_1 \geq 0, \ldots, \beta_5 \geq 0$, where
	\begin{align*}
		\beta_i &= e_i^* \quad \text{for} \quad i = 1,2,3,4, \\
		\beta_5 &= - e_1^* + e_2^* + e_3^* + e_4^*,
	\end{align*}
	and $\{e_1^*, \ldots, e_4^*\}$ is the standard basis of the dual vector space $V^* = (\mathbb{R}^4)^*$.
	Here we have chosen the indices $1, \ldots ,5$ so that the facets supported by $\beta_1, \ldots, \beta_5$ match the adjacency relations shown in Figure \ref{example 1}.
	
	We next find $\beta_6 \in V^*$ such that the six inequalities $\beta_i \geq 0$, $i = 1, \ldots, 6$, determine a convex polytope combinatorially equivalent to $\mathcal{G}$.
	Let $T \subset \mathbb{S}^3$ be the triangular prism determined by the first five inequalities $\beta_1 \geq 0, \ldots, \beta_5 \geq 0$.
	For each $i \in \{1, \ldots, 5\}$, let $T_i$ denote the facet of $T$ supported by $\beta_i$.
	Whenever $1 \leq i < j < k \leq 5$ and $T_i \cap T_j \cap T_k$ is a vertex of $T$, we denote this vertex by $v_{ijk}$.
	
	The sixth supporting hyperplane $\mathbb{S} \{\beta_6 = 0\} \subset \mathbb{S}^3$ must pass through the vertices $v_{124}$ and $v_{235}$.
	Using the explicit forms of the first five functionals $\beta_1, \ldots, \beta_5$, we obtain
	\[
	v_{124} = [(0,0,1,0)], \quad v_{235} = [(1,0,0,1)].
	\]
	Substituting these vertices into $\beta_6$, we find that $\beta_6$ has the form 
	\begin{align}\label{beta_6}
		\beta_6 = d_1 e_1^* + d_2 e_2^* - d_1 e_4^*
	\end{align}
	for some $d_1, d_2 \in \mathbb{R}$.
	
	The facets $T_1$, $T_2$, and $T_3$ meet at the vertex $v_{123}=[(0,0,0,1)]$, which must lie in the half-space $\mathbb S\{\beta_6>0\}$.
	Evaluating \eqref{beta_6} at $v_{123}$ gives $-d_1>0$, and hence $d_1<0$.
	Since multiplying $\beta_6$ by a positive scalar leaves the hyperplane $\mathbb{S} \{\beta_6 = 0\}$ invariant, we may assume
	\[
	\beta_6 = - e_1^* + d e_2^* + e_4^*
	\]
	for some $d \in \mathbb{R}$.
	
	We find the range of $d$ such that the inequalities $\beta_1 \geq 0, \ldots, \beta_6 \geq 0$ determine a convex 3-polytope combinatorially equivalent to $\mathcal{G}$.
	Let $E_{45}$ denote the edge $T_4 \cap T_5$ of the triangular prism $T$.
	The inequalities $\beta_1 \geq 0, \ldots, \beta_6 \geq 0$ determine a convex 3-polytope combinatorially equivalent to $\mathcal{G}$ if and only if the hyperplane $\mathbb{S} \{\beta_6 = 0\}$ intersects the interior of the edge $E_{45}$ transversely.
	This condition holds if and only if the three vertices $v_{134}$, $v_{123}$, and $v_{345}$ lie in the open half-space $\mathbb{S} \{\beta_6 > 0\}$, while the vertex $v_{245}$ lies in the opposite half-space $\mathbb{S} \{\beta_6 < 0\}$.
	Using the explicit forms of the first five functionals $\beta_1, \ldots, \beta_5$, we obtain $v_{134} = [(0,1,0,0)]$, $v_{123} = [(0,0,0,1)]$, $v_{345} = [(1,1,0,0)]$ and $v_{245} = [(1,0,1,0)]$.
	Substituting these vertices, we find that $d$ determines an element of $\widetilde{\mathcal{E}} (\mathcal{G})$ if and only if $d \in (1, \infty)$.
	
	For each $d\in(1,\infty)$, put $\beta^d:=(\beta_1,\ldots,\beta_5,\beta_6^d)$, where
	\begin{align*}
		\beta_i &= e_i^* \quad \text{for } i=1,2,3,4, \\
		\beta_5 &= -e_1^*+e_2^*+e_3^*+e_4^*, \\
		\beta_6^d &= -e_1^*+d e_2^*+e_4^*.
	\end{align*}
	We show that the parameter $d$ is uniquely determined by the $\mathbb G$-orbit of $\beta^d$.
	Suppose that $d,d'\in(1,\infty)$ and that $g=(A,c_1,\ldots,c_6)\in\mathbb G$ satisfies
	\[
	g\cdot\beta^d=\beta^{d'}.
	\]
	Comparing the first four components gives $A=\text{diag}(c_1^{-1},c_2^{-1},c_3^{-1},c_4^{-1})$.
	Comparing the fifth components gives $c_1=c_2=c_3=c_4=c_5$.
	Since $A\in\text{SL}_{\pm}(V)$ and all the numbers $c_i$ are positive, it follows that $c_1=\cdots=c_5=1$ and $A=\text{Id}$.
	Comparing the sixth components then gives $c_6=1$ and $d=d'$.
	Thus distinct values of $d$ determine distinct points of $\mathcal E(\mathcal G)$.
	
	The preceding normalization shows that every point of $\mathcal E(\mathcal G)$ is represented by $\beta^d$ for some $d\in(1,\infty)$.
	Consequently, the map
	\[
	(1,\infty) \to \mathcal E(\mathcal G),
	\quad
	d \mapsto [\beta^d],
	\]
	is a homeomorphism.
	
	\subsubsection{The restricted deformation spaces and the global deformation space}\label{example-1-fibers}
	
	We compute the restricted deformation spaces for the same labeled combinatorial 3-polytope $\mathcal{G}$.
	Recall from Section \ref{Section 6.2} that the restricted deformation spaces in $\mathcal{C} (\mathcal{G})$ can be identified with the nonempty fibers of the map $\mathcal{C} (\mathcal{G}) \to \mathcal{RS} (\mathcal{G})$.
	In the commuting diagram of Figure \ref{commuting diagram}, we see that they can also be identified with the nonempty fibers of the map $\overline{\rho}: \mathcal{D} (\mathcal{G}) \to \mathcal{E} (\mathcal{G})$.
	Moreover, since $\eta: \widetilde{\mathcal{D}} (\mathcal{G}) \to \mathcal{D} (\mathcal{G})$ maps each nonempty fiber of $\rho: \widetilde{\mathcal{D}} (\mathcal{G}) \to \widetilde{\mathcal{E}} (\mathcal{G})$ homeomorphically onto a nonempty fiber of $\overline{\rho}$, it suffices to compute the nonempty fibers of $\rho$.
	
	Here is the outline of the computation.
	Let $d \in (1, \infty)$, and let $\beta := (\beta_1, \ldots, \beta_5, \beta_6^d) \in \widetilde{\mathcal{E}} (\mathcal{G})$ as described in the previous computation.
	An element $w = (w_1, \ldots, w_6) \in V^6$ lies in $\rho^{-1} (\beta)$ if and only if $(\beta, w)$ satisfies the conditions \hyperref[V1]{(V1)} through \hyperref[V6]{(V6)}. 
	The condition \hyperref[V6]{(V6)} is already given since $\beta \in \widetilde{\mathcal{E}} (\mathcal{G})$.
	We first find such $w_1, \ldots, w_6$ successively by solving the linear equations arising from the conditions \hyperref[V1]{(V1)}, \hyperref[V3]{(V3)}, and \hyperref[V4]{(V4)}.
	Then we sort out the solutions $w_1, \ldots, w_6$ satisfying the condition \hyperref[V2]{(V2)}.
	Checking condition \hyperref[V5]{(V5)} separately is unnecessary, since the restrictions of conditions \hyperref[V1]{(V1)}--\hyperref[V4]{(V4)} to adjacent pairs, together with condition \hyperref[V6]{(V6)}, imply condition \hyperref[V5]{(V5)} by the argument at the end of the proof of Lemma~\ref{standard chart lemma}.
	In the process, we will also determine all values of $d$ for which the fiber $\rho^{-1} (\beta)$ is nonempty.
	
	We first find the vector $w_1$.
	By conditions \hyperref[V1]{(V1)} and \hyperref[V3]{(V3)}, it is uniquely determined by the equalities
	\begin{align*}
		\beta_1 (w_1) &= 2, \\
		\beta_i (w_1) &= 0 \quad \text{for} \quad i = 2,3,4.
	\end{align*}
	Hence, we obtain $w_1 = (2, 0, 0, 0)$.
	
	Similarly, the equations $\beta_2 (w_2) = 2$, $\beta_i (w_2) = 0$ for $i = 1, 3$, and $\beta_6^d (w_2) = 0$ yield $w_2 = (0, 2, 0, -2d)$.
	
	The third vector $w_3$ is not uniquely determined.
	It is given by only the three equations $\beta_3 (w_3) = 2$ and $\beta_i (w_3) = 0$ for $i = 1,2$.
	Thus, $w_3 = (0, 0, 2, x)$ for some $x \in \mathbb{R}$.
	Since $\beta_4 = e_4^*$, we have $\beta_4(w_3)=x$.
	The edge $\mathcal{G}_3 \cap \mathcal{G}_4$ has order $3$, so condition \hyperref[V4]{(V4)} gives
	\[
	\beta_3(w_4)\beta_4(w_3) = 4\cos^2\left(\frac{\pi}{3}\right) = 1.
	\]
	Hence, $x\neq 0$.
	More generally, for any $w\in\rho^{-1}(\beta)$ and any distinct indices $i,j$ with $m_{i,j}\neq 2$, condition \hyperref[V2]{(V2)}, together with \hyperref[V4]{(V4)} when $\mathcal{G}_i$ and $\mathcal{G}_j$ are adjacent and \hyperref[V5]{(V5)} otherwise, implies that $\beta_i(w_j)<0$.
	The remaining inequalities determine the admissible range of $x$.
	
	For $w_4$, we need an additional condition \hyperref[V4]{(V4)}, and the resulting system of equations is
	\begin{align*}
		\beta_4 (w_4) &= 2, \\
		\beta_i (w_4) &= 0 \quad \text{for} \quad i = 1, 6, \\
		\beta_3 (w_4) &= \frac{4 \cos^2 \left( \frac{\pi}{3} \right) }{\beta_4 (w_3)} = \frac{1}{x}.
	\end{align*}
	Thus, we obtain $w_4 = \left(0, - \frac{2}{d}, \frac{1}{x}, 2 \right)$, which depends uniquely on $d$ and $x$.
	
	Solving the corresponding system of four linear equations arising from conditions \hyperref[V1]{(V1)}, \hyperref[V3]{(V3)}, and \hyperref[V4]{(V4)} gives
	\[
	w_5 = \left( \frac{d\left(\frac{1}{x+2}-2\right)}{d-1} +\frac{1}{-\frac{2}{d}+\frac{1}{x}+2}, \frac{\frac{1}{x+2}-2}{d-1}, \frac{1}{x+2}, \frac{1}{-\frac{2}{d}+\frac{1}{x}+2} \right).
	\]
	The corresponding system for $w_6$ gives $w_6=(-2,0,-2,0)$.
	
	To summarize, we have obtained
	\begin{align*}
		w_1 &= (2, 0, 0, 0), \quad w_2 = (0, 2, 0, -2d), \quad w_3 = (0, 0, 2, x), \quad w_4 = \left(0, - \frac{2}{d}, \frac{1}{x}, 2 \right), \\
		w_5 &= \left( \frac{d \left( \frac{1}{x + 2} - 2 \right)}{d - 1} + \frac{1}{- \frac{2}{d} + \frac{1}{x} + 2}, \frac{\frac{1}{x + 2} - 2}{d - 1}, \frac{1}{x + 2}, \frac{1}{- \frac{2}{d} + \frac{1}{x} + 2} \right), \\
		w_6 &= (-2, 0, -2, 0).
	\end{align*}
	The realization space is parametrized by $d$, and the restricted deformation space corresponding to the parameter $d$ is then parametrized by $x$.
	
	Next, we find all $d$ and $x$ such that $(\beta, w)$ satisfies condition \hyperref[V2]{(V2)}.
	Since the cases $\beta_i (w_j) = 0$ are already covered in the construction of $w_1, \ldots, w_6$, we only need to consider the pairs $i, j$ such that $m_{i, j} \geq 3$.
	By \cite[Proposition~14]{MR0302779}, once condition \hyperref[V2]{(V2)} has been verified for adjacent facets, condition \hyperref[V2]{(V2)} for nonadjacent facets and condition \hyperref[V5]{(V5)} follow from the linear relations among the supporting functionals.
	It is therefore enough to check the adjacent edges of order at least $3$.
	
	It follows that the values of $d$ and $x$ satisfying condition \hyperref[V2]{(V2)} must satisfy the inequalities
	\[
	\beta_3 (w_4) < 0, \quad \beta_3 (w_5) < 0, \quad \beta_4 (w_5) < 0.
	\]
	These inequalities, combined with the condition $d \in (1, \infty)$, yield 
	\[
	1 < d < \frac{4}{3}, \quad \frac{d}{2 - 2d} < x < - 2.
	\]
	
	We conclude that the restricted deformation space $\rho^{-1} (\beta)$ over the point $\beta = (\beta_1, \ldots, \beta_5, \beta_6^d) \in \widetilde{\mathcal{E}} (\mathcal{G})$ is nonempty if and only if $1 < d < \frac{4}{3}$.
	Moreover, each nonempty fiber $\rho^{-1} (\beta)$, parametrized by $x$, is homeomorphic to an open interval.
	
	From these observations, we conclude that the deformation space $\mathcal{C} (\mathcal{G})$ is homeomorphic to an open 2-cell $\mathbb{R}^2$, which provides additional information to the statement of Theorem \ref{main theorem}.
	
	\subsection{A forgetful map that is neither surjective nor a fiber bundle}\label{counterexample}
	
	Let $\mathcal{G}$ be the labeled combinatorial $3$-polytope shown in Figure~\ref{example 2}.
	Its underlying polytope is a truncation polytope, and \cite[Theorem~3.16]{MR2660566} shows that $\mathcal{C}(\mathcal{G})$ is homeomorphic to $\mathbb{R}^4 \sqcup \mathbb{R}^4$.
	The purpose of this example is to show that, unlike the map in Section~\ref{Example 7.1}, the forgetful map $\kappa$ need not be a fiber bundle over its image.
	In Section~\ref{example-2-realization}, we parameterize $\mathcal{E}(\mathcal{G})$ and choose a path $\gamma(s)$ of realizations.
	In Section~\ref{example-2-fibers}, we compute the fibers along this path, show that their number of connected components changes from one to two, and deduce that $\kappa$ is not a fiber bundle over its image.
	The occurrence of empty fibers along the same path also shows that $\kappa$ is not surjective.
	We omit computational details that are identical to those in Section~\ref{Example 7.1}.

\begin{figure}[t]
	\centering
	\begin{tikzpicture}
		\draw (0, 0.8) -- (0.69282, -0.4);
		\draw (0.69282, -0.4) -- (-0.69282, -0.4);
		\draw (-0.69282, -0.4) -- (0, 0.8);
		
		\draw (-0.4330125, 1.75) -- (0.4330125, 1.75);
		\draw (0.4330125, 1.75) -- (0, 1.2);
		\draw (0, 1.2) -- (-0.4330125, 1.75);
		
		\draw (0.4330125, 1.75) -- (2.1650625, -1.25);
		\draw (2.1650625, -1.25) -- (-2.1650625, -1.25);
		\draw (-2.1650625, -1.25) -- (-0.4330125, 1.75);
		
		\draw (0, 0.8) -- (0, 1.2);
		\draw (0.69282, -0.4) -- (2.1650625, -1.25);
		\draw (-0.69282, -0.4) -- (-2.1650625, -1.25);
		
		\node at (-1,0) {\footnotesize $\mathcal{G}_1$};
		\node at (1,0) {\footnotesize $\mathcal{G}_2$};
		\node at (0,-0.9) {\footnotesize $\mathcal{G}_3$};
		\node at (0,0) {\footnotesize $\mathcal{G}_4$};
		\node at (2.5,0) {\footnotesize $\mathcal{G}_5$};
		\node at (0,1.55) {\footnotesize $\mathcal{G}_6$};
		
		\node at (-0.1,1) {\tiny 3};
		\node at (-0.44641, 0.3) {\tiny 2};
		\node at (-1.42894125, -0.725) {\tiny 3};
		\node at (-1.3990375, 0.35) {\tiny 3};
		\node at (-0.31650625, 1.375) {\tiny 2};
		
		\node at (0.31650625, 1.375) {\tiny 2};
		\node at (1.3990375, 0.35) {\tiny 3};
		\node at (1.42894125, -0.725) {\tiny 4};
		\node at (0.44641, 0.3) {\tiny 2};
		
		\node at (0, -0.5) {\tiny 3};
		\node at (0, -1.35) {\tiny 2};
		
		\node at (0, 1.85) {\tiny 3};
	\end{tikzpicture}
	\caption{A labeled combinatorial $3$-polytope for which the forgetful map is neither surjective nor a fiber bundle; see Section~\ref{counterexample}.
	}
	\label{example 2}
\end{figure}

\subsubsection{The realization space and a path of realizations}\label{example-2-realization}

	We compute $\widetilde{\mathcal{E}}(\mathcal{G})$ by the same normalization as in Section~\ref{Example 7.1}.
	Specifically, we consider the unique representative $\beta = (\beta_1, \ldots, \beta_6)$ in each $\mathbb{G}$-orbit such that the first five functionals determine the triangular prism $T$ as in Section \ref{Example 7.1}.
	This prism is defined by the five inequalities $\beta_1 \geq 0, \ldots, \beta_5 \geq 0$, where $\beta_i = e_i^*$ for $i = 1,2,3,4$ (the standard basis elements of $V^*$), and $\beta_5 = e_1^* + e_2^* + e_3^* - e_4^*$.
	
	Note that the positions of $\mathcal{G}_1$ and $\mathcal{G}_4$ are switched compared to Section \ref{Example 7.1}.
	This adjustment ensures that the orderability conditions are satisfied.
	
	For each $i \in \{1, \ldots, 5\}$, let $T_i$ denote the facet of $T$ supported by $\beta_i$, and let $v_{125}$ denote the vertex $T_1 \cap T_2 \cap T_5$.
	For the sixth element $\beta_6 \in V^*$, the set $\mathbb{S} (\{x \in V \ | \ \beta_i (x) \geq 0 \ \text{for} \ i = 1, \ldots, 6\})$ forms a convex 3-polytope combinatorially equivalent to $\mathcal{G}$ if and only if the vertex $v_{125}$ lies in the open half-space $\mathbb{S} (\{\beta_6 < 0\})$ and the remaining five vertices of $T$ lie in the opposite half-space $\mathbb{S} (\{\beta_6 > 0\})$.
	
	These conditions, along with normalization by positive scalars, yield the parametrization of $\mathcal{E} (\mathcal{G})$ given by
	\begin{align*}
		\beta_i &= e_i^* \quad \text{for} \quad i = 1,2,3,4, \\
		\beta_5 &= e_1^* + e_2^* + e_3^* - e_4^*, \\
		\beta_6 &= d_1 e_1^* + d_2 e_2^* + d_3 e_3^* - e_4^*,
	\end{align*}
	with $(d_1, d_2, d_3) \in (1, \infty) \times (1, \infty) \times (0, 1)$.
	The number of parameters also coincides with Theorem \ref{realization spaces} since we have $e = 12$.
	
	Next, we consider a specific path inside $\mathcal{E}(\mathcal{G})$ using the parametrization.
	To detect a change in the topology of the fibers, consider the path $\gamma(s)\in\widetilde{\mathcal{E}}(\mathcal{G})$, $s\in(1,2]$, defined by
	\[
	d_1=s,\qquad d_2=2,\qquad d_3=\frac{1}{2}.
	\]
	
	\subsubsection{The restricted deformation spaces along the path}\label{example-2-fibers}
	
	We examine the fiber $\rho^{-1} (\gamma (s))$ over the point $\gamma (s)$, i.e., the restricted deformation space associated with the polytope determined by $\gamma (s)$.
	Solving conditions \hyperref[V1]{(V1)}, \hyperref[V3]{(V3)}, and \hyperref[V4]{(V4)} recursively gives
	\begin{align*}
		w_1 &= \left( 2, t, -4 s - 4 t, 0 \right), \\
		w_2 &= \left( \frac{1}{t}, 2, -\frac{2(s + 4t)}{t}, 0 \right), \\
		w_3 &= \left( -\frac{1}{2(2s + 2t)}, -\frac{t}{s + 4t}, 2, -\frac{\frac{s}{2} - 4s^2 + 2t - 18st - 14t^2}{(2s + 2t)(s + 4t)} \right), \\
		w_4 &= \left( 0, 0, \frac{(2s + 2t)(s + 4t)}{-\frac{s}{2} + 4s^2 - 2t + 18st + 14t^2}, 2 \right),
	\end{align*}
	\[
	w_5=
	\begin{pmatrix}
		-\dfrac{1}{2\left(-1+2s+\frac{3t}{2}\right)}\\[2mm]
		\dfrac{t}{2\left(\frac12-s-3t\right)}\\[2mm]
		0\\[2mm]
		-\dfrac{\frac34-\frac{7s}{2}+4s^2-\frac{13t}{2}+16st+\frac{39t^2}{4}}
		{\left(\frac12-s-3t\right)\left(1-2s-\frac{3t}{2}\right)}
	\end{pmatrix}^T,
	\]
	\[
	w_6
	= \begin{pmatrix}
		0 \\[2mm]
		0 \\[2mm]
		\frac{2\left(-1 + \frac{9s}{2} - 5s^2 + \frac{29t}{4} - \frac{35st}{2} - 12t^2\right)}{\frac{3}{4} - \frac{13s}{4} + \frac{7s^2}{2} - \frac{11t}{2} + \frac{25st}{2} + \frac{33t^2}{4}} \\[2mm]
		\frac{2\left(-\frac{5}{4} + \frac{11s}{2} - 6s^2 + \frac{73t}{8} - \frac{85st}{4} - \frac{57t^2}{4}\right)}{\frac{3}{4} - \frac{13s}{4} + \frac{7s^2}{2} - \frac{11t}{2} + \frac{25st}{2} + \frac{33t^2}{4}}
	\end{pmatrix}^T.
	\]
	Here $t\in\mathbb{R}$ is the remaining parameter after solving the equality conditions \hyperref[V1]{(V1)}, \hyperref[V3]{(V3)}, and \hyperref[V4]{(V4)}.
	The remaining inequalities determine the admissible values of $t$, the nonemptiness of the fiber, and its number of connected components.
	
	Since $\gamma(s)\in\widetilde{\mathcal{E}}(\mathcal{G})$, condition \hyperref[V6]{(V6)} is already satisfied, and the displayed vectors satisfy conditions \hyperref[V1]{(V1)}, \hyperref[V3]{(V3)}, and \hyperref[V4]{(V4)} by construction.
	As in Section \ref{Example 7.1}, condition \hyperref[V5]{(V5)} does not need to be checked separately.
	Thus, the remaining nontrivial requirements in condition \hyperref[V2]{(V2)} are exactly the following inequalities:
	\begin{align*}
		\beta_3 (w_4) &< 0, \quad \beta_2 (w_5) < 0, \quad \beta_1 (w_5) < 0, \quad \beta_1 (w_2) < 0, \\
		\beta_2 (w_3) &< 0, \quad \beta_1 (w_3) < 0, \quad \beta_5 (w_6) < 0.
	\end{align*}

For $1<s\leq 2$, the five inequalities
\[
\beta_2(w_5)<0,\qquad
\beta_1(w_5)<0,\qquad
\beta_1(w_2)<0,\qquad
\beta_2(w_3)<0,\qquad
\beta_1(w_3)<0
\]
are equivalent to
\[
L_s<t<U_s,
\qquad
L_s:=\frac{2-4s}{3},
\qquad
U_s:=-\frac{s}{4}.
\]
On this interval, the remaining two expressions can be written as
\[
\beta_3(w_4)
=
\frac{4(s+t)(s+4t)}{q_4(s,t)}
\]
and
\[
\beta_5(w_6)
=
\frac{(2s+6t-1)(4s+3t-2)}{q_6(s,t)},
\]
where
\[
q_4(s,t):=28t^2+(36s-4)t+8s^2-s
\]
and
\[
q_6(s,t):=33t^2+(50s-22)t+14s^2-13s+3.
\]
For $L_s<t<U_s$, we have
\[
s+t>0,\qquad
s+4t<0,\qquad
2s+6t-1<0,\qquad
4s+3t-2>0.
\]
It follows that the inequalities $\beta_3(w_4)<0$ and $\beta_5(w_6)<0$ are equivalent to
\[
q_4(s,t)>0
\qquad\text{and}\qquad
q_6(s,t)>0,
\]
respectively.

Set
\[
R_4(s):=\frac{1-9s-\sqrt{1-11s+25s^2}}{14}
\]
and
\[
R_6(s):=\frac{11-25s+\sqrt{22-121s+163s^2}}{33}.
\]
These are the relevant roots of $q_4(s,t)$ and $q_6(s,t)$, respectively.
Also set
\[
s_0:=\frac{193+3\sqrt{697}}{176},
\qquad
s_1:=\frac{4}{19}\left(5+\sqrt{6}\right).
\]
A comparison of the roots with $L_s$ and $U_s$ gives the complete admissible range of $t$ for $1<s\leq2$:
\[
t\in
\begin{cases}
	\varnothing,
	& 1<s\leq s_0,\\[2mm]
	\bigl(L_s,R_4(s)\bigr),
	& s_0<s\leq s_1,\\[2mm]
	\bigl(L_s,R_4(s)\bigr)\cup\bigl(R_6(s),U_s\bigr),
	& s_1<s\leq2.
\end{cases}
\]

It follows that the fiber $\rho^{-1}(\gamma(s))$ is empty for $1<s\leq s_0$, connected for $s_0<s\leq s_1$, and has two connected components for $s_1<s\leq2$.
Let $\overline{\gamma}(s)\in\mathcal{RS}(\mathcal{G})$ denote the realization corresponding to $\pi(\gamma(s))\in\mathcal{E}(\mathcal{G})$.
By Lemma \ref{fiber identification}, whenever the fiber is nonempty, the fiber of $\kappa$ over $\overline{\gamma}(s)$ is homeomorphic to $\rho^{-1}(\gamma(s))$.
If $\kappa$ were a fiber bundle over its image, then its pullback along the path $\overline{\gamma}|_{(s_0,2]}$ would be a fiber bundle over the connected interval $(s_0,2]$, and all its fibers would be homeomorphic.
This is impossible since the number of connected components of the fibers changes from one to two at $s=s_1$.
Therefore, $\kappa$ is not a fiber bundle over its image.
The empty fibers for $1<s\leq s_0$ also show that $\kappa$ is not surjective.

\appendix

\section{Invariant projective subspaces for normal-type Coxeter $3$-polytopes}
\label{appendix:invariant-subspaces}

As noted in Remark~\ref{rem:geomded}, the proof of \cite[Proposition~2]{MR2247648} contains gaps in the argument excluding holonomy-invariant projective subspaces.
The proof of Theorem~\ref{main theorem} uses only the restricted-deformation-space theorem stated as Theorem~\ref{restricted}, so the correction below is not a separate step in our main argument.
We include it to complete the cited proof.
Specifically, Proposition~\ref{proposition:appendix-choi-correction} and its proof replace the part of \cite[Proposition~2]{MR2247648} concerning an invariant projective line, a disjoint union of two projective lines, and an invariant projective plane.

Let $(P,r_1,\ldots,r_f)$ be a Coxeter $3$-polytope realizing a labeled combinatorial $3$-polytope $\mathcal{G}$, and write
\[
P=\mathbb{S}(K),
\qquad
K=\{x\in\mathbb{R}^4\mid \alpha_i(x)\geq 0\text{ for }i=1,\ldots,f\},
\qquad
r_i=\text{Id}-\alpha_i\otimes v_i,
\]
where $\alpha_i(v_i)=2$.  Denote by $P_i$ the facet supported by $\alpha_i$, and put
\[
H_i:=\mathbb{S}(\ker\alpha_i)\subset\mathbb{S}^3.
\]
The unordered antipodal pair
\[
\mathfrak{p}_i:=\mathbb{S}(\mathbb{R}v_i)=\{[v_i],[-v_i]\}
\]
is called the \emph{polar $0$-sphere} of $r_i$; the line $\mathbb{R}v_i$ is the $(-1)$-eigenspace of $r_i$.  This is the object called the set of ``antipodally fixed points'' in \cite{MR2247648}.  A \emph{projective $k$-subspace} of $\mathbb{S}^3$ means a subset $\mathbb{S}(W)$ for a $(k+1)$-dimensional linear subspace $W\subset\mathbb{R}^4$.  
Throughout this appendix, the word \emph{invariant} is used in the setwise sense.
Thus, if a group $\Lambda$ acts on a set $Y$, a subset $X\subset Y$ is $\Lambda$-invariant if $\lambda(X)=X$ for every $\lambda\in\Lambda$; this does not require $\Lambda$ to fix the points of $X$ individually.

Let
\[
\Gamma:=\langle r_1,\ldots,r_f\rangle,
\qquad
\widetilde{\Omega}:=\text{Int}\!\left(\bigcup_{\gamma\in\Gamma}\gamma K\right) \subset \mathbb{R}^4,
\qquad
\Omega:=\mathbb{S}(\widetilde{\Omega}) \subset \mathbb{S}^3.
\]  
We shall repeatedly use the following consequences of Vinberg's relations:
\begin{equation}
	\alpha_i(v_j)\leq 0\quad(i\neq j),
	\qquad
	\alpha_i(v_j)=0\Longleftrightarrow \alpha_j(v_i)=0.
	\label{eq:appendix-vinberg-zero}
\end{equation}
Moreover, if the equivalent equalities in \eqref{eq:appendix-vinberg-zero} hold, then $r_i$ and $r_j$ commute.  Since the presentation in Theorem~\ref{Vinberg 2} identifies the order of $r_i r_j$ with the Coxeter order prescribed by $\mathcal{G}$, the facets $P_i$ and $P_j$ are then adjacent and their common edge has order $2$.

\begin{proposition}
	\label{proposition:appendix-choi-correction}
	Let $\mathcal{G}$ be a labeled combinatorial $3$-polytope of normal type, and let $(P,r_1,\ldots,r_f)$ be a Coxeter $3$-polytope realizing $\mathcal{G}$.  Then the reflection group $\Gamma=\langle r_1,\ldots,r_f\rangle$ preserves none of the following subsets of $\mathbb{S}^3$:
	\begin{enumerate}[(i)]
		\item a projective line;
		\item a disjoint union of two projective lines;
		\item a projective plane.
	\end{enumerate}
\end{proposition}

This is the invariant-subspace assertion that replaces the corresponding part of the proof of \cite[Proposition~2]{MR2247648}.
We first establish two elementary lemmas used in the proof.

\begin{lemma}
	\label{lemma:appendix-invariant-subspace}
	Let $W\subset\mathbb{R}^4$ be a linear subspace invariant under $r_i$.  Then either
	\[
	W\subset\ker\alpha_i
	\qquad\text{or}\qquad
	v_i\in W.
	\]
\end{lemma}

\begin{proof}
	If $W\not\subset\ker\alpha_i$, choose $w\in W$ with $\alpha_i(w)\neq 0$.  Since both $w$ and $r_i(w)$ belong to $W$, we have
	\[
	v_i=\frac{w-r_i(w)}{\alpha_i(w)}\in W.
	\]
\end{proof}

\begin{lemma}
	\label{lemma:appendix-no-invariant-plane}
	Suppose that $\mathcal{G}$ is of normal type.  Then $\Gamma$ does not preserve a projective plane in $\mathbb{S}^3$.
\end{lemma}

\begin{proof}
	Suppose, to the contrary, that $S=\mathbb S(W)$ is $\Gamma$-invariant for a $3$-dimensional linear subspace $W\subset\mathbb R^4$.
	By Lemma~\ref{lemma:appendix-invariant-subspace}, for every $i$, either $H_i=S$ or $v_i\in W$.
	At most one mirror can equal $S$ since two distinct facets of $P$ cannot have the same supporting projective plane.
	
	Assume first that $H_i=S$ for one index $i$.  For every $j\neq i$ we then have $v_j\in W=\ker\alpha_i$, and hence $\alpha_i(v_j)=\alpha_j(v_i)=0$.
	Consequently, $P_i$ is adjacent to every other facet $P_j$, and every edge of $P_i$ has order $2$.  In addition, since $\alpha_i(v_i)=2 > 0$ and $\alpha_j(v_i)=0$ for $j \ne i$, the point $[v_i]$ belongs to every facet except $P_i$.
	
	The point $[v_i]$ must be a vertex of $P$.
	Otherwise, it would lie in the relative interior of an edge or a facet of $P$.
	In the former case $[v_i]$ lies in precisely two facets of $P$, and in the latter case $[v_i]$ lies in precisely one facet of $P$.
	Both are impossible since $P_i$ is adjacent to at least three facets.
	Therefore
	\[
	P=\text{conv}\bigl(P_i\cup\{[v_i]\}\bigr)
	\]
	is a cone over the polygon $P_i$, and all edges of its base polygon have order $2$.  This is precisely case~\textup{(i)} excluded in Definition~\ref{normal type}, a contradiction.
	
	We may therefore assume that no mirror equals $S$.  Lemma~\ref{lemma:appendix-invariant-subspace} gives
	\begin{equation}
		v_i\in W\qquad\text{for every }i.
		\label{eq:appendix-all-polars-in-plane}
	\end{equation}
	Choose a nonzero linear functional $\beta\in(\mathbb{R}^4)^*$ with $W=\ker\beta$.  Equation~\eqref{eq:appendix-all-polars-in-plane} implies
	\[
	\beta\circ r_i=\beta
	\qquad\text{for every }i.
	\]
	If $S\cap\Omega=\varnothing$, the convex domain $\Omega$ is contained in one of the two affine charts
	\[
	\{[x]\in\mathbb{S}^3\mid \beta(x)>0\},
	\qquad
	\{[x]\in\mathbb{S}^3\mid \beta(x)<0\}.
	\]
	That chart is preserved by all $r_i$.  Thus $\mathcal{G}$ admits an affine Coxeter group representation, which is case~\textup{(iv)} excluded in Definition~\ref{normal type}.
	
	It remains to consider $S\cap\Omega\neq\varnothing$.
	Since each element of $\Gamma$ preserves $S$ and no mirror equals $S$, we have $S \ne \gamma H_i$ for all $\gamma \in \Gamma$ and $i$.
	Since $\Gamma$ is countable and each intersection
	$S\cap\gamma H_i$ is a proper projective subspace of $S$, their union
	has empty interior in $S$.  We may therefore choose a point of
	$S\cap\Omega$ that does not belong to any $\gamma H_i$.
	Such a point lies in $\gamma\text{Int}(P)$ for some
	$\gamma\in\Gamma$.
	Since $S$ is $\Gamma$-invariant, applying $\gamma^{-1}$ gives a point $[x]\in S\cap\text{Int}(P)$.
	Choose its lift $x\in W\cap\text{Int}(K)$, so that $\alpha_i(x)>0$ for every $i$.
	
	The section $Q:=P\cap S$ is a properly convex polygon in the projective plane $S$.  
	In fact, every facet of $P$ cuts out a facet (hence an edge) of $Q$.  
	More precisely, for each $i$ we set $x_i:=x-\frac{\alpha_i(x)}{2}v_i\in W$.
	Then $\alpha_i(x_i)=0$, while for $j\neq i$ Vinberg's inequalities give
	\[
	\alpha_j(x_i)
	=\alpha_j(x)-\frac{\alpha_i(x)}{2}\alpha_j(v_i)
	\geq \alpha_j(x)>0.
	\]
	The point $[x_i]$ lies in the relative interior of $P_i\cap S$, so $Q_i:=P_i\cap S$ is a side of $Q$.
	The restrictions $r_i|_W$ are reflections in these sides, and the images of $Q$ under the elements of $\Gamma$ tessellate $S\cap\Omega$.
	The induced action of $\Gamma$ on $S\cap\Omega$ is faithful.
	Indeed, an element acting trivially on $S$ fixes $[x]\in\text{Int}(P)$, whereas the interiors of $P$ and the images of $P$ under any nontrivial elements of $\Gamma$ are disjoint.
	It follows that $Q$ is a fundamental polygon for the action of $\Gamma$ on $S\cap\Omega$, and this action is generated by the reflections $r_i|_W$.

	We now show that every edge of $P$ also meets $S$.  Let $P_i\cap P_j$ be an edge of order $m_{i,j}$.
	By Vinberg's relation,
	\[
	\alpha_i (v_j) \alpha_j (v_i) = 4\cos^2\!\left(\frac{\pi}{m_{i,j}}\right)<4.
	\]
	It follows that there are unique numbers $s,t$ satisfying
	\begin{equation}
		\begin{pmatrix}
			2&\alpha_i (v_j)\\
			\alpha_j (v_i) & 2
		\end{pmatrix}
		\begin{pmatrix}s\\t\end{pmatrix}
		=
		\begin{pmatrix}\alpha_i(x)\\\alpha_j(x)\end{pmatrix}.
		\label{eq:appendix-edge-system}
	\end{equation}
	We have
	\[
	\begin{pmatrix}
		2&\alpha_i (v_j)\\
		\alpha_j (v_i)&2
	\end{pmatrix}^{-1}
	=\frac{1}{4-\alpha_i (v_j) \alpha_j (v_i)}
	\begin{pmatrix}
		2&-\alpha_i (v_j)\\
		-\alpha_j (v_i)&2
	\end{pmatrix}.
	\]
	Since $[x]\in\text{Int}(P)$, both entries of the vector on the right-hand side of \eqref{eq:appendix-edge-system} are positive.
	The inverse matrix displayed above has nonnegative entries and positive diagonal entries.
	Therefore, $s,t>0$.
	Set
	\[
	y:=x-sv_i-tv_j\in W.
	\]
	The defining equations of $s,t$ give $\alpha_i(y)=\alpha_j(y)=0$.  For $k\notin\{i,j\}$,
	\[
	\alpha_k(y)
	=\alpha_k(x)-s\alpha_k(v_i)-t\alpha_k(v_j)
	\geq\alpha_k(x)>0.
	\]
	Thus $[y]$ lies in the relative interior of the edge $P_i\cap P_j$ and is a vertex of the polygon $Q$, at which precisely the two sides $Q_i$ and $Q_j$ meet.
	
	We fix an index $i$.  
	The facet $P_i$ is a convex polygon and therefore has at least three distinct edges.  
	Applying the preceding construction to those edges produces at least three distinct vertices of the single side $Q_i$.  
	The vertices are distinct since the point obtained from an edge $P_i\cap P_j$ lies on no facet other than $P_i$ and $P_j$.  
	This is impossible, since a side of a convex polygon has exactly two vertices.  
	We have obtained the final contradiction.
\end{proof}

\begin{proof}[Proof of Proposition~\ref{proposition:appendix-choi-correction}]
	Part~\textup{(iii)} is Lemma~\ref{lemma:appendix-no-invariant-plane}.  
	We first prove part~\textup{(i)}.  
	Suppose that $l=\mathbb{S}(L)$ is a $\Gamma$-invariant projective line, where $\dim L=2$.  
	We call a facet $P_i$ \emph{parallel to $l$} if $l\subset H_i$.
	In terms of the corresponding linear subspaces, this means that $L\subset\ker\alpha_i$.
	A facet that is not parallel to $l$ is called \emph{transverse to $l$}.
	By Lemma~\ref{lemma:appendix-invariant-subspace},
	\begin{equation}
		v_i\in L
		\qquad\text{for every facet transverse to }l.
		\label{eq:appendix-transverse-polar}
	\end{equation}
	
	There are at most two facets parallel to $l$.
	Indeed, let
	\[
	K^*:=\{\lambda\in(\mathbb{R}^4)^*\mid \lambda(x)\geq 0\text{ for every }x\in K\}
	\]
	be the dual cone.
	We call a convex cone \emph{sharp} if it contains no affine line.
	Since $K$ has nonempty interior, the cone $K^*$ is sharp.
	Indeed, if $\lambda_0+\mathbb{R}\mu\subset K^*$, then $\lambda_0(x)+t\mu(x)\geq 0$ for every $x\in K$ and every $t\in\mathbb{R}$, so $\mu$ vanishes on $K$.
	The nonempty interior of $K$ then gives $\mu=0$.
	Set
	\[
	L^\perp:=\{\lambda\in(\mathbb{R}^4)^*\mid \lambda|_L=0\}.
	\]
	The cone $K^*\cap L^\perp$ is sharp as a subcone of $K^*$ and has dimension at most $2$.
	The functionals supporting facets parallel to $l$ belong to this cone.
	Among any three rays in a sharp convex cone of dimension at most $2$, one lies in the positive span of the other two.
	The corresponding inequality would then be implied by the other two inequalities, contrary to the assumed nonredundancy of the facet inequalities defining $P$.
	
	If no facet is parallel to $l$, equation~\eqref{eq:appendix-transverse-polar} shows that every $v_i$ lies in $L$.  Hence every $3$-dimensional linear subspace containing $L$ is invariant under all $r_i$, producing a $\Gamma$-invariant projective plane.  This contradicts Lemma~\ref{lemma:appendix-no-invariant-plane}.
	
	Suppose that exactly one facet, say $P_i$, is parallel to $l$.  
	Since $\alpha_i(v_i)=2$ and $\alpha_i|_L=0$, we have $v_i\notin L$.  
	The 3-dimensional subspace $W := L+\mathbb{R}v_i$ is invariant under every generator.
	The polar of $r_i$ lies in $W$, and the polar of every other generator lies in $L$ by \eqref{eq:appendix-transverse-polar}.
	Again, $\mathbb{S}(W)$ is an invariant projective plane, a contradiction.
	
	Finally, suppose that exactly two facets, say $P_i$ and $P_j$, are parallel to $l$.  Put
	\[
	W:=L+\text{Span}\{v_i,v_j\}.
	\]
	This subspace is invariant under every generator.  
	Since $v_i\notin L$, we have $\dim W\geq 3$.  
	If $\dim W=3$, then $\mathbb{S}(W)$ is an invariant projective plane, contradicting Lemma~\ref{lemma:appendix-no-invariant-plane}.  
	We may therefore assume $W=\mathbb{R}^4$.  
	Then
	\[
	U:=\text{Span}\{v_i,v_j\}
	\]
	has dimension $2$ and $\mathbb{R}^4=L\oplus U$.
	
	For every transverse facet $P_k$, we have $v_k\in L$ and hence
	\[
	\alpha_i(v_k)=\alpha_j(v_k)=0.
	\]
	By \eqref{eq:appendix-vinberg-zero}, this gives
	\[
	\alpha_k(v_i)=\alpha_k(v_j)=0,
	\]
	so $\alpha_k\in U^\perp := \{\lambda \in (\mathbb{R}^4)^* \ | \ \lambda|_U = 0\}$.  
	The cone $K^*\cap U^\perp$ is sharp as a subcone of the sharp cone $K^*$.
	It has dimension at most $2$ since $\dim U^\perp=2$.
	Among any three rays in this cone, one lies in the positive span of the other two.
	If three of these rays arose from nonredundant facet inequalities, the inequality corresponding to that ray would be implied by the other two.
	At most two rays in $K^*\cap U^\perp$ can therefore arise from facet-defining functionals.
	It follows that there are at most two transverse facets.
	In particular, $P$ has at most four facets.
	A $3$-polytope has at least four facets, so $P$ is a tetrahedron with exactly two parallel and two transverse facets.
	
	Every pair of facets of a tetrahedron is adjacent.
	For each transverse facet $P_k$, the equalities $\alpha_i(v_k)=\alpha_k(v_i)=0$ and $\alpha_j(v_k)=\alpha_k(v_j)=0$ imply that both $r_i$ and $r_j$ commute with $r_k$, by the consequence of Vinberg's relations stated after \eqref{eq:appendix-vinberg-zero}.
	Moreover, the subgroup generated by the two parallel reflections and the subgroup generated by the two transverse reflections are both finite dihedral groups.
	Since these two subgroups commute and together generate $\Gamma$, the group $\Gamma$ is finite.
	This is case~\textup{(iii)} excluded in Definition~\ref{normal type}.
	This proves part~\textup{(i)}.
	
	For part~\textup{(ii)}, suppose that $\Gamma$ preserves a disjoint union $l_1\sqcup l_2$ of projective lines.  
	Every projective line meets the mirror $H_i$ of a reflection $r_i$, since a $2$-dimensional linear subspace of $\mathbb{R}^4$ meets the $3$-dimensional subspace $\ker\alpha_i$ nontrivially.  
	If $r_i$ exchanged $l_1$ and $l_2$, a point of $l_1\cap H_i$, fixed by $r_i$, would also belong to $l_2$, contradicting the disjointness.  
	Hence every $r_i$ preserves each line separately.  
	In particular, $l_1$ is $\Gamma$-invariant, contradicting part~\textup{(i)}.
\end{proof}

\bibliographystyle{plain}
\bibliography{References}

\end{document}